\documentclass[11pt]{article}

\usepackage{amsmath}
\usepackage{amsfonts}
\usepackage{amssymb}
\usepackage{amsthm}
\usepackage{graphicx}
\usepackage{cite}
\usepackage{caption}
\usepackage{url}

\usepackage{algorithm}
\usepackage{algpseudocode}
\usepackage{mathtools}
\usepackage{multirow}
\usepackage{bbm}
\usepackage[table]{xcolor}
\usepackage{xr}

\usepackage{subfigure}
\usepackage{empheq}
\usepackage[figuresright]{rotating}

\usepackage{wrapfig}
\usepackage{comment}

\makeatletter
\newcommand*{\addFileDependency}[1]{% argument=file name and extension
  \typeout{(#1)}
  \@addtofilelist{#1}
  \IfFileExists{#1}{}{\typeout{No file #1.}}
}
\makeatother

\newtheorem{remark}{Remark}

\algblock{ParFor}{EndParFor}
\algnewcommand\algorithmicparfor{\textbf{parallel for}}
\algnewcommand\algorithmicpardo{\textbf{do}}
\algnewcommand\algorithmicendparfor{\textbf{end\ parallel for}}
\algrenewtext{ParFor}[1]{\algorithmicparfor\ #1\ \algorithmicpardo}
\algrenewtext{EndParFor}{\algorithmicendparfor}

\definecolor{Reviewer1}{rgb}{0,0,0}
\definecolor{Reviewer2}{rgb}{0,0,0}
\definecolor{Reviewer4}{rgb}{0,0,0}
\definecolor{Reviewer5}{rgb}{0,0,0}
\definecolor{Reviewer2_new}{rgb}{0.6,0,0}
\definecolor{Reviewer4_new}{rgb}{0,0.6,0}

\title{High order-accurate solution of scattering integral equations\\
  with unbounded solutions at corners} \author{Constantine Sideris$^*$, Davit Aslanyan\footnote{Electrical and Computer Engineering, University of Southern California, Los Angeles, CA 90089, USA}
  \ and Oscar P. Bruno\footnote{Computing and Mathematical Sciences,
    Caltech, Pasadena, CA 91125, USA}}

\begin{document}
\date{}
\maketitle
% \tableofcontents
\begin{abstract}
  Although high-order Maxwell integral equation solvers provide
  significant advantages in terms of speed and accuracy over
  corresponding low-order integral methods, their performance
  significantly degrades in presence of non-smooth geometries---owing
  to field enhancement and singularities that arise at sharp edges and
  corners which, if left untreated, give rise to significant accuracy
  losses. The problem is particularly challenging in cases in which
  the ``density'' (i.e., the solution of the integral equation) tends
  to infinity at corners and edges---a difficulty that can be bypassed
  for 2D configurations, but which is unavoidable in 3D Maxwell
  integral formulations, wherein the component tangential to an edge
  of the electrical-current integral density vector tends to infinity
  at the edge. In order to tackle the problem this paper restricts
  attention to the simplest context in which the unbounded-density
  difficulty arises, namely, integral formulations in 2D space whose
  integral density blows up at corners; the strategies proposed,
  however, generalize directly to the 3D context. The novel
  methodologies presented in this paper yield high-order convergence
  for such challenging equations and achieve highly accurate solutions
  (even near edges and corners) without requiring a priori analysis of
  the geometry or use of singular bases.
\end{abstract}
\vspace{0.5 cm}
\maketitle

\section{Introduction}

The computational solution of electromagnetic scattering problems
plays fundamental roles in many areas of science and engineering,
ranging from photonics to remote sensing, and including
communications, optics and electromagnetic compatibility among many
others.  Computational electromagnetic scattering solvers based on
boundary-integral equation (BIE) offer a number of advantages, such as
reliance on discretizations of lower dimensionality; an ability to
directly utilize CAD geometries as input without significant geometry
processing, and associated faithfulness to the geometries provided; a
lack of dispersion errors of significance; and, when used in
conjunction with suitable acceleration
methods~\cite{bauinger_bruno,bles,br-k1,rokhlinetal}, high efficiency
in terms of computing times and parallel performance. The
Chebyshev-based ``rectangular-polar'' Boundary Integral Equation (RP)
method~\cite{2,3} is a computationally efficient, high-order Nystr\"om
method for the discretization and solution of scattering and radiation
problems in 2D and 3D which, based on use of representation of surface
currents in terms of local Chebyshev polynomial expansions, achieves
high-order accuracy for smooth geometries. The method's convergence
rate, however, can drop to even less than linear for practical
geometries that include geometric singularities such as sharp edges
and corners, thus eliminating, in such cases, one of the significant
benefits offered by the methodology. The difficulty arises from the
singular behavior of field components and current densities (and/or
their derivatives) near geometric singularities. This behavior leads
to poor numerical approximations and significant integration errors in
these regions. The challenge is further exacerbated in cases where the
solution of the integral equation (the ``density'') diverges at
corners and edges. While such issues can often be circumvented in 2D
configurations~\cite{0}, they are unavoidable in 3D Maxwell integral
formulations, where the tangential component of the current-density
vector along an edge becomes infinite. While such singular fields near
an edge can be represented by an asymptotic expansion of quasistatic
singular fields, which can in turn be incorporated in numerical
methods in the form of singular basis functions, unfortunately,
algorithms based on such approaches may suffer from significantly
increased complexity in view of the geometrical dependence of the
singular exponents involved---which, in particular, are not even known
in closed form for corners in 3D space and would thus require
numerical evaluation, which in itself amounts to a challenging
computational problem.

In order to tackle the problem at hand, this paper restricts attention
to the simplest context in which the unbounded-density difficulty
arises, namely, an integral formulation in 2D space whose integral
density blows up at corners. The novel methodologies presented in what
follows recover high-order convergence for such challenging
geometries, and they achieve highly accurate solutions (even near
edges and corners) without requiring a priori analysis of the geometry
or use of singular bases. In detail, this paper considers a
prototypical integral-equation (for the 2D Neumann problem for the
Helmhotz equation) whose density solution is unbounded at
corners. Alternative integral equation formulations can be used for
the 2D Neumann-Helmholtz problem considered here: for example, the
surface currents associated with the Direct Regularized Combined Field
Integral Equations (DCFIE-R) for the Neumann problem are bounded near
edges and corners~\cite{0}, and they are thus amenable to effective
treatment by means of well known regularization methods. Our selection
of integral equation with unbounded density was made, precisely, to
enable the development of strategies capable to accurately deal with
infinite densities, which could then be generalized and applied in the
corresponding 3D cases for which, as mentioned above, infinite
densities at edges and corners exist regardless of the formulation
used. As discused below, the novel methodologies presented in this
paper yield high-order convergence for such challenging geometries and
achieve highly accurate solutions (even near edges and corners)
without requiring a priori analysis of the geometry or use of singular
bases.

The following discussion provides a brief review of the literature,
beginning with contributions on finite element methods (FEM) for
partial differential equations (PDEs) in domains with geometric
singularities. While FEM results in sparse systems of linear
equations, these methods are not always competitive with integral
equation methods. The latter, leveraging lower-dimensional
discretizations, can benefit from fast solvers~\cite{br-k1,bles,r3}
and, when applicable, significantly outperform FEM
approaches. Nonetheless, FEM-based techniques for boundary value
problems in domains with singular boundaries are crucial in many
fields and have inspired a rich body of research. For instance,
spatially refined meshes near regions of geometric singularities have
been explored
in~\cite{babushkaNM1979,cox_fixCMA1984,babuska_94,Babuska:The,Demkowitz:On,Devloo:Recursive}. An
alternative approach explicitly incorporates the \textit{known}
singular behavior into the Galerkin basis, as discussed
in~\cite{sternIJNME1979,fixJCP1973,Lin:Singular,Hughes:Techniques}. The
``DtN finite element method,'' detailed
in~\cite{kellerIJNME1992,Wu:Discrete,Wu:FE}, identifies neighborhoods
of corner singularities and poses a new boundary value problem on the
complement of these neighborhoods using artificial boundary conditions
derived from Dirichlet-to-Neumann maps.

Green function-based integral equation methods, in turn, particularly
high-order techniques for solving two- and three-dimensional problems
in domains with smooth boundaries, are
well-established~\cite{kress_ie_book,br-k1,ganesh_graham_ho,strain}. Significant
advancements have also been achieved for non-smooth domains, as
demonstrated by works such
as~\cite{maichak_stephan_97,atkinson_book}. For instance, the approach
in~\cite{maichak_stephan_97} utilizes theoretically robust first-kind
(singular or hypersingular) integral equations combined with
high-order Galerkin boundary element methods. These techniques are
effective for both Dirichlet and Neumann problems; however, they
exhibit limited accuracy in scenarios where the integral equation
solutions, as considered in this study, become
unbounded~\cite{heuer_mellado_stephan}.

Nystr\"om-based integral-equation methodologies, which align more
closely with the approach proposed here, have also been applied with
notable success~\cite{8,Atkinson87,Graham88}; a comprehensive
discussion of this literature can be found
in~\cite[Chap. 8]{atkinson_book}. These methodologies use special
graded-mesh quadratures to achieve high-order accuracy for solving the
{\em Dirichlet problem for the Laplace equation} in two-dimensional
domains with corners via second-kind integral equations. However, the
integral equations addressed in these contributions do not involve
hypersingular operators, which are central to the present study.  As
demonstrated in Section~\ref{comparisonagainstRP}, directly extending
these methods to Maxwell problems which involve hypersingular
operators as well as unbounded integral densities does not generally
yield accurate solutions. Additional Nystr\"om methodologies including
both hypersingular operators and unbounded solutions are presented
in~\cite{bremer2012fast,helsing2013solving} and references
therein. Reference~\cite{bremer2012fast} builds upon the ``augmented''
integral equation method~\cite{yaghjian1981augmented}, in that it
relies on use of integral operators supported in curves in the
interior of the obstacle---which may pose some challenges for
realistic 3D configurations. The method demonstrates high accuracies
of the computed fields, albeit at significant distances from the
corner points. Reference~\cite{helsing2013solving}, in turn, presents
the recursive ``RCIP'' method, which recursively subdivides regions
near corners. That reference presents highly accurate results for
a domain with a $90^\circ$ corner, and it indicates that the opening
angle should be in the interval $\pi/3\leq \theta\leq 5\pi/3$.

This paper presents a numerical methodology that recovers high-order
convergence in presence of unbounded densities at corners, and which
produces highly accurate solutions even at extremely close proximity
from the singularity and for extremely sharp corners (such as, e.g.,
near machine precision accuracy at a distance $d=10^{-8}$ for a corner
with an interior angle of $0.01$ radians $\approx 0.57^\circ$); it is
expected that forthcoming 3D versions of the proposed algorithm will
prove similarly effective.  The proposed approach utilizes (i)~A {\em
  regularized-operator} integral formulation
(Section~\ref{regularized}) that eliminates the deleterious spectral
character arising from hypersingular operators inherent in the
integral formulation; and, (ii)~A regularization procedure and
associated regular density unknown $\psi$ that incorporates a
corner-regularization Jacobian factor (Section~\ref{corn_reg});
together with, (iii)~Specialized near-corner and self-interaction
quadrature techniques that effectively manage the extreme proximity of
discretization points, and associated cancellation errors, that are
introduced by the corner regularization methods
(Sections~\ref{interpolation_header_section}
and~\ref{precomputationswithcov}). A variety of numerical examples
provided in this paper illustrate the character of the proposed
approach. Results for the wave-equation solution $u$ with near
machine-precision accuracy are obtained arbitrarily close to the
corner: such accuracies are indeed consistently demonstrated in this
paper at a nominal extremely-close distance of $d = 10^{-8}$ away from
corners (see Remark~\ref{rem-11-dig} for additional details in this
regard).

This paper is organized as follows. Section~\ref{int_eq} presents
necessary background concerning integral equations (which give rise to
infinite densities at corners), and it reviews the Chebyshev-based
(Nystr\"om) rectangular polar discretization method~\cite{2,3}.
Section~\ref{op-corn_regul} then presents the proposed strategy for
high-order solution of the infinite-density problem, which involves
use of a change of integration variables and associated graded meshes
as well as a change of unknown, and which requires in particular
detailed treatment of arithmetic near-limits of the form $0\cdot\infty$ in
finite-precision arithmetic. A variety of numerical results
illustrating the character of the proposed methods are presented in
Section~\ref{numer}, and a few concluding remarks, finally, are
provided in Section~\ref{conc}.

\section{Integral Equation Formulations for Smooth Surfaces\label{int_eq}}

We consider the problem of scattering of Transverse Magnetic (TM)
Electromagnetic waves by a closed perfect electrical conductor (PEC)
scatterer $\Omega$ with boundary $\Gamma$ in 2D space, wherein the
spatial variable $r = (x,y)$ is used, under an incident wave with
transverse magnetic-field component
$u^\mathrm{inc} = H_z^\mathrm{inc}(r)$ such as, e.g.
$H_z^\mathrm{inc}(r) = e^{ikx}$.  The resulting electromagnetic
scattered field in the domain $\Omega^c$ exterior to $\Gamma$ can be
determined on the basis of the transverse scattered magnetic field
component $u(r) = H_z^\mathrm{scat}(r)$, which is a solution of the
problem
\begin{empheq}[left=\empheqlbrace\ ]{align}
& \Delta u + k^2 u = 0  \quad \mbox{in}\quad \Omega^c,\label{helm}\\
&\frac{\partial u}{\partial n} = - \frac{\partial u^\mathrm{inc}}{\partial n}  \quad \mbox{on}\quad \Gamma.\label{Neumann}
\end{empheq}

\subsection{Magnetic Field Integral Equation}\label{MFIESection}
Using the single-layer representation
\begin{equation}\label{SL}
    u(r) = \int_\Gamma{G_{k}(r,r')\phi d\ell'},\quad r\in\Omega^c,
  \end{equation}
  where
\begin{equation}
    G_{k}(r,r') = \frac{i}{4}H_0^1(k|r-r'|)
\end{equation}
denotes the outgoing Helmholtz free-space Green's function, the
imposition of the Neumann boundary condition~\eqref{Neumann} results
in the Magnetic Field Integral Equation (MFIE)~\cite{7} 
\begin{equation}\label{MFIE}
  -\frac{1}{2}\phi(r) + \frac{\partial}{\partial n(r)}\int_{\Gamma}G_{k}(r,r')\phi(r')d\ell(r') = -\frac{\partial u^\mathrm{inc}}{\partial n(r)}, \quad r\in\Gamma.
\end{equation}
Once a solution $\phi$ of this equation has been produced and
substituted into~\eqref{SL}, the solution $u$
of~\eqref{helm}-\eqref{Neumann} can be obtained by straightforward
numerical quadrature.

The MFIE~\eqref{MFIE} is not uniquely solvable at a
scatterer-dependent countable sequence of frequencies
$k = k_1,k_2,\dots$, see e.g.~\cite{7}, around each one of which the
associated numerical implementations become ill-conditioned and
inaccurate. This problem can be eliminated by utilizing instead a
certain combined field integral equation formulation, which is
considered in Section \ref{sec:CFIE}. Unfortunately, however, that
formulation incorporates a hypersingular operator whose practical
implementation presents certain difficulties, most notably in terms of
the associated conditioning characteristics and large iteration
numbers it requires if iterative solvers are used for its
solution. Following~\cite{1}, Section \ref{regularized} introduces a
second version of the combined equation, where the hypersingular
operator is regularized using a smoothing operator. This approach
resolves the difficulties associated with the hypersingular nature of
the problem, at least for smooth surfaces $\Gamma$. However, new
complications emerge when this formulation is applied to domains with
corners, since, as is the case for the regular CFIE formulation, its
solutions $\phi$ are unbounded at corner points. An additional
reformulation of the problem, that introduces the new
corner-regularized unknown and associated equation, which is one of
the key contributions in this paper, is then presented in Section
\ref{corn_reg}. Sections \ref{interpolation_header_section} and
\ref{precomputationswithcov} detail additional components of the
proposed solver designed to address finite precision issues associated
with the use of fine singularity-resolving meshes and fine-scale
pre-computation problems. The remainder of the present
Section~\ref{int_eq}, in turn, introduces the aforementioned
smooth-surface uniquely solvable CFIE equation and associated
parametrization and differentiation methods, which are needed for
reference in the subsequent Section~\ref{op-corn_regul}.

\subsection{Unique Solvability: Combined Field Integral Equation\label{sec:CFIE}}
Using now the representation
\begin{equation}
    u = -i\eta\int_\Gamma{G_{k}(r,r')\phi d\ell'} +\int_\Gamma{\frac{\partial G_{k}(r,r')}{\partial{n}(r')}\phi d\ell'},\quad r\in\Omega^c,
\end{equation} 
instead of~\eqref{SL} (see~\cite{1} and references therein),
enforcement of the Neumann boundary conditions~\eqref{Neumann} yields
the Combined Field Integral Equation (CFIE) formulation
\begin{equation}\label{CFIE}
    \frac{i\eta}{2}\phi(r) - i\eta\int_{\Gamma}\frac{\partial G_{k}(r,r')}{\partial n(r)}\phi(r')d\ell(r') + \frac{\partial}{\partial n(r)}\int_{\Gamma}\frac{\partial G_{k}(r,r')}{\partial n(r')}\phi(r')d\ell(r') 
    = -\frac{\partial u^\mathrm{inc}}{\partial n(r)},
\end{equation}
where $\eta \neq 0$ is a coupling parameter. Following ~\cite{1}, throughout this paper the value $\eta = 1$ is used. In equation ~\eqref{CFIE} and throughout this paper,  $\frac{\partial}{\partial n(r)}$ (resp. $\frac{\partial}{\partial n(r')}$) denotes differentiation with respect to $r$ in the direction of $n(r)$  (resp. differentiation with respect to $r'$ in the direction of $n(r')$). %The normal derivative of the double-layer operator kernel is hypersingular, and takes the form
%\begin{equation}
%    \frac{\partial ^2G_{k}(r,r')}{\partial n(r) \partial n(r')} = \frac{i}{4}{k}^2H_2^1(k|r-r'|)\cdot\frac{n(r')\cdot(r-r')}{|r-r'|} \cdot\frac{n(r)\cdot(r-r')}{|r-r'|} -\frac{i}{4}kH_1^1(k|r-r'|)\cdot\frac{n(r')\cdot n(r)}{|r-r'|}.
%\end{equation}
The last operator on the left-hand side of~\eqref{CFIE}, namely, the normal derivative of the double-layer operator, is a hypersingular operator---whose computation by direct implementation of normal derivatives can be rather challenging. Fortunately, however, the hypersingular operator may be expressed in the form 
\begin{equation}\label{tangential}
    \frac{\partial}{\partial n(r)}\int_{\Gamma}\frac{\partial G_{k}(r,r')}{\partial n(r')}\phi(r')d\ell(r') = {k}^2 \int_{\Gamma} G_{k}(r,r') (n(r) \cdot n(r'))\phi(r') \, d\ell(r') +  \partial_{\tau(r)}\int_{\Gamma}  G_{k}(r,r') \partial_{\tau(r')} \phi(r') \, d\ell(r');
\end{equation}
(see~\cite{1} and references therein), where $\partial_{\tau(r)}$ and
$\partial_{\tau(r')}$ denote the operators of differentiation with respect to the unit tangent vectors $\tau(r)$ and $\tau(r')$
at the respective ``target" and ``source" points $r$ and $r'$ on $\Gamma$; clearly equation~\eqref{tangential}  does not require evaluation of derivatives in directions other than those tangential to $\Gamma$, and can therefore be computed effectively by means of the discretization methods used in this paper. (Note that the identity~\eqref{tangential} holds for closed curves, provided the tangential vector $\tau(r)$ is defined in a consistent manner throughout the scatterer's boundary, e.g., in counterclockwise fashion.) In sum, defining the operators
\begin{equation}\label{T1}
 T_1[\phi](r) = \int_{\Gamma}\frac{\partial G_{k}(r,r')}{\partial n(r)}\phi(r')d\ell(r'),
\end{equation}
 
\begin{equation}\label{T2}
T_2[\phi](r) = \int_{\Gamma} G_{k}(r,r') (n(r) \cdot n(r'))\phi(r') \, d\ell(r')\quad\mbox{and}
\end{equation}
 
\begin{equation}\label{T3}
T_3[\phi](r) = \partial_{\tau(r)}\int_{\Gamma}  G_{k}(r,r') \partial_{\tau(r')} \phi(r') \, d\ell(r'),
\end{equation}
the CFIE equation \eqref{CFIE} takes the form
\begin{equation}\label{CFIE_subdivided}
  \frac{i\eta}{2}\phi(r) - i\eta T_1[\phi](r) + k^2 T_2[\phi](r) + T_3[\phi](r) 
  = -\frac{\partial u^\mathrm{inc}}{\partial n(r)}.
\end{equation}
Note that the integrals required by the operators $T_1$ and $T_2$ may
be expressed in the forms
  \begin{equation}\label{gen_int}
\mathcal{I}[H,f](r) = \int_\Gamma H(r,r') f(r')d\ell(r'),
\end{equation}
while $T_3$ may be presented in the form
\begin{equation}\label{gen_int_dif}
\partial_{\tau(r)}\mathcal{I}[H,f](r) = \partial_{\tau(r)}\int_\Gamma H(r,r') f(r')d\ell(r'),
\end{equation}
for certain kernels $H$ and functions $f$: we have
$H(r,r') = \frac{\partial G_{k}(r,r')}{\partial n(r)}$ and
$f(r') = \phi(r')$ for the operator $T_1$;
$H(r,r') = G_{k}(r,r') (n(r) \cdot n(r'))$ and $f(r') = \phi(r')$ for
the operator $T_2$; and $H(r,r') = G_{k}(r,r')$ and
$f(r') = \partial_{\tau(r')}\phi(r')$ for the opeator $T_3$. Clearly,
the operator $T_3$ requires an additional differentiation of the
unknown $\phi$---which, as discussed in.~\cite{1}, impacts
significantly upon the accuracy and conditioning of the resulting
numerical implementation.

In preparation for the introduction of the main operator- and
corner-regularization methods in Section~\ref{op-corn_regul}, the
following two sections describe, in the context of smooth surfaces,
the particular parametrization and discretization methods used, and
the associated methods employed for differentiation of surface
variables and evaluation of Green function-based integral operators.

\subsection{Chebyshev-Based Parametrization and Rectangular-Polar
  Integration: smooth-scatterer
  case\label{cheby_old}}
Throughout this paper the discretizations of integral operators such
as those considered in the previous sections and elsewhere are
implemented on the basis of of the Chebyshev-based Rectangular Polar
(RP) Nystr\"om methods introduced in~\cite{2,3}. This section briefly
reviews the RP method as it applies in the context of smooth surfaces
and for operators that do not contain tangential differentiation
(e.g., the operators~\eqref{T1} and~\eqref{T2} above); additions
presented in Section \ref{numerical_derivative} then enable the
discretization of the operator~\eqref{T3} with similar efficiency and
accuracy.

\begin{figure}[!ht]
    \centering
    \includegraphics[width=0.4\linewidth]{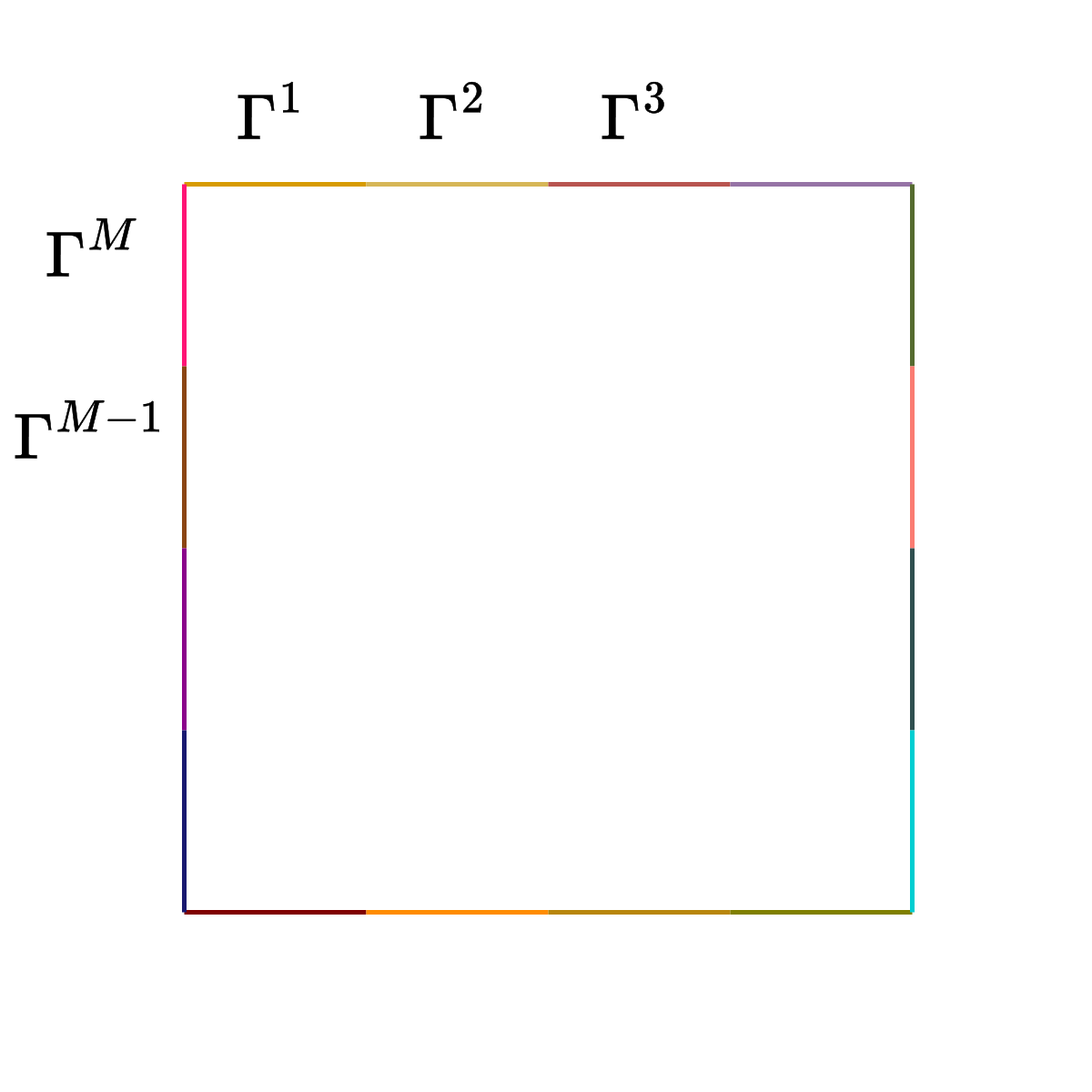}
    \vspace{-1cm}
    \caption{Partitioning of the domain boundary $\Gamma$ as the
      union~\eqref{partition} of nonoverlapping patches $\Gamma^q$
      ($1\leq q\leq M$).}
    \label{fig:patch_splitting}
\end{figure}

The RP method relies on use of a partition of curve   (surface in 3D) $\Gamma$ into a set a set of $M$ non-overlapping patches $\Gamma^q$
\begin{equation}\label{partition}
\Gamma = \bigcup_{q=1}^M \Gamma^q,
\end{equation}
as depicted in Figure~\ref{fig:patch_splitting}. Without loss of generality it is assumed that each curve $\Gamma^q$ is parametrized by a vector function with scalar argument
\begin{equation}\label{param}
r_q:[0,1] \to \Gamma^q \in \mathbb{R}^2.
\end{equation}
(Note that, while previous implementations of the RP method use patch
parametrizations defined in the interval $[-1,1]$---given the reliance
of the RP method on Chebyshev expansions---, throughout this paper the
parameter segment $[0,1]$ is utilized instead, on account of the
``zero-centering'' technique introduced per Sections~\ref{corn_reg}
and~\ref{interpolation_header_section}, which requires the use of such
a parametrization interval.)
On the basis of the parametrizations $r_q$, a function $f(r)$ defined
on $\Gamma^q$ corresponds to the function $f(r_q(s))$ defined for
$s\in [0,1]$, which can then be discretized by means of the $Q$-point
Chebyshev grid given by~\cite{10}
\begin{equation}\label{zeroonenodes}
s_j = \frac{\zeta_j+1}{2}  \quad j = 0, \ldots, Q - 1; \quad \omega_j = \frac{W_j}{2}, 
\end{equation}
where
\begin{equation}\label{negoneonenodes}
\zeta_j = \cos \left( \pi \frac{2j + 1}{2Q} \right), \quad j = 0, \ldots, Q - 1; 
\quad W_j = \frac{2}{Q} \left( 1 - 2 \sum_{k=1}^{Q/2} \frac{1}{4k^2 - 1} \cos \left( k\pi \frac{2j + 1}{Q} \right) \right).
\end{equation}
Clearly, using $Q$ points per patch for  a total $M$, the total number of discretization points used is
\begin{equation}
  \label{eq:N}
  N = MQ.
\end{equation}

The implementation of the integral operators at hand requires
consideration of the line element
\begin{equation}\label{LE}
    d\ell_{q'} = L_{q'}(s')ds,\quad \mbox{where}\quad L_{q'}(s')\mbox{ = }\sqrt{\left(\frac{dx_{q'}}{ds}(s')\right)^2 + \left(\frac{dy_{q'}}{ds}(s')\right)^2},
\end{equation} 
as well as the unit tangential and normal basis vectors which, at a
point $r_{q'}(s') = \left(x_{q'}(s'),y_{q'}(s')\right)\in\Gamma^{q'}$, are given
by
\begin{equation}
  {\tau}(r_{q'}(s')) = \left(\frac{dx_{q'}}{ds}(s'),\frac{dy_{q'}}{ds}(s')\right)\bigg/L_{q'}(s'),\quad
  n(r_{q'}(s')) = \left(\frac{dx_{q'}}{ds}(s'),\frac{dy_{q'}}{ds}(s')\right)\bigg/L_{q'}(s')
\end{equation}
In order to achieve high computational efficiency and accuracy the RP
quadrature algorithm proceeds by separately treating singular,
near-singular, and regular interactions. In view of~\eqref{gen_int}
and~\eqref{partition}, the integral $\mathcal{I}[H,f]$ may be
expressed in the form
\begin{equation}
  \label{eq:sum-int}
  \mathcal{I}[H,f](r) = \sum_{{q'}=1}^M \mathcal{I}^{q'}[H,f](r),
\end{equation}
where
\begin{equation}
\mathcal{I}^{q'}[H,f](r) =
\int_{\Gamma^{q'}} H(r,r')f(r')d\ell(r') = \int_0^1 H(r,r_{q'}(s'))f(r_{q'}(s')) L_{q'}(s') ds'.  \\[8pt]
\end{equation}

% \begin{equation}
% I^q(r) =
% \begin{cases} 
% \int_{\Gamma^q} H(r,r')f(r')d\ell(r'),  \\[8pt]
% \int_{\Gamma^q} H(r,r')\frac{\partial \phi(r')}{\partial \tau(r')}d\ell(r'), 
% \end{cases}
% \end{equation} with
% \begin{equation}
% H^q(r,r') = 
% \begin{cases} 
% \frac{\partial G_k(r,r')}{\partial n(r)}, & \text{for $T_1$ \eqref{T1}}, \\[8pt]
% G_k(r,r')n(r)\cdot n(r'), & \text{for $T_2$ \eqref{T2}}, \\[8pt]
% G_k(r,r') \quad \mbox{for $T_3$
% \end{cases}
% \end{equation}

For target points $r$ away from the source patch $\Gamma^{q'}$ (as
determined by a certain proximity parameter $\delta^\mathrm{prox}>0$)
the kernels of the integral operators are smooth over the source
patch, and, at least in the smooth-surface case considered in this
section, so are the density functions $f$; thus, for such far
interactions (on smooth domains) the method directly utilizes
Fej\'er's first quadrature rule for integration, and the approximation
\begin{equation}\label{far_nocorner}
  \mathcal{I}^{{q'}} [H,f](r) \approx \sum_{j=0}^{N-1} H (r, r_{{q'}}(s_j))f(r_{{q'}}(s_j)) L_{{q'}} (s_j) \omega_j , \quad \mbox{for } \quad r\quad \mbox{away from} \quad \Gamma^{{q'}}
\end{equation}
results, where $s_j$ and $\omega_j$ denote the integration points and
weights introduced in equation \eqref{zeroonenodes}.

When the target point $r$ is close to or within the source patch, on
the other hand, the kernels of the integral operators are near
singular or singular, respectively, and accurate evaluation of the
corresponding integrals requires use of specialized integration
techniques. To produce such integrals, the RP method relies on the
Chebyshev approximation
\begin{equation}\label{density_cheby_expansion}
  f(r_{q'}(s')) \approx \sum_{m=0}^{N-1} a_m^{q'} T_m
  (2s'-1)
\end{equation}
which can be obtained for the necessary functions $f$ (equal either to
$\phi$ or its tangential derivative, see equation~\eqref{gen_int_dif})
from the available values of the density $\phi$ on the Chebyshev mesh.
Using this expansion, the operator values in singular and
near-singular cases can be efficiently and accurately computed on the
basis of the numerically precomputed integration weights
\begin{equation}\label{precomp_def} \beta_m^{q'}[H](r) =
\int_{0}^{1}H(r,r_{q'}(s')) \, L_{q'}(s') \, T_m(2s'-1) \, ds',
\end{equation}
These weights can be precomputed by integrating the kernel against
each Chebyshev polynomial (e.g., as described in~\cite{2,3}) and
stored for repeated use. Using these weights we obtain the accurate
approximation
\begin{equation}\label{precomp_sum}
  \mathcal{I}^{q'}[H,f](r) \approx \sum_{m=0}^{N-1} a_m^{q'} \beta_m^{q'}[H](r), \quad \mbox{for} \quad r\quad \mbox{on or near} \quad \Gamma^{q'}.
\end{equation}
This approach applies to each one of the kernels $H$ listed near the
end of Section~\ref{sec:CFIE} and, together with the aforementioned
smooth integration methods for $r$ far from $\Gamma^{q'}$, it enables
the necessary integrals over $\Gamma$ to be computed efficiently via
equation~\eqref{eq:sum-int}.  Note that, per the discussion above in
this section, precomputations are necessary for target discretization
points $r = r_{q'}(s_j)$ on the source patch $\Gamma^{q'}$, as well as
for near singular target discretization points $r = r_{q'}(s_j)$ on
the nearby patches $\Gamma^{q'}$ for which the distance between the
discretization points $r = r_{q'}(s_j)$ and $\Gamma_{q'}$ is less than
$\delta^\mathrm{prox}$.
\subsection{Differentiation}\label{numerical_derivative}
The evaluation of the operator $T_3$ in~\eqref{T3}, which contains two tangential derivatives, is discussed in what follows. In view of~\eqref{partition} and writing

\begin{equation}\label{T33}
\left .\partial_{\tau(r)}\int_{\Gamma^{q'}}  G_{k}(r,r') \partial_{\tau(r')} \phi(r') \, d\ell(r')\right|_{r=r_q(s)} = \frac{1}{L_q(s)}\frac{\partial}{\partial s}\int_{0}^{1}   G_k(r_q(s),r_{q'}(s')) \frac{\partial \phi(r_{q'}(s'))}{\partial s'} ds', 
\end{equation}
% \begin{equation}\label{T3}
% T_3[\phi] = \partial_{\tau(r)}\int_{\Gamma}  G_{k}(r,r') \partial_{\tau(r')} \phi(r') \, d\ell(r')
% \end{equation}
% \begin{equation}\label{oldchebydiff}
% \int_{\Gamma^{q'}}  \frac{\partial G_k(r_q(s),r_{q'}
% (s'))}{\partial \tau(r_q(s))} \frac{\partial \phi(r_{q'}(s'))}{\partial \tau(r_{q'}(s'))} L(s') ds' = \frac{1}{L(s)}\frac{\partial}{\partial s}\int_{0}^{1}   G_k(r_q(s),r_{q'}(s')) \frac{\partial \phi(r_{q'}(s'))}{\partial s'} ds', 
% \end{equation}
(since $\partial_{\tau(r_q(s))} = \frac{1}{L_q(s)}\partial/\partial s$
and
$\partial_{\tau(r_q'(s'))} = \frac{1}{L_{q'}(s')}\partial/\partial
s'$), it suffices to evaluate the right-hand side expression in this
equation for all $q$ and $q'$ with $1\leq q,q' \leq M$ and then sum
over $q'$ for each $q$.  To do this we leverage the Chebyshev
structure inherent in the RP algorithm described in
Section~\ref{cheby_old} together with the numerically stable Chebyshev
differentiation approach~\cite{10} (which is briefly reviewed below),
all of which results in the following procedure.
\begin{enumerate}
\item For each $q'$ evaluate the derivative $\left .\frac{\partial \phi(r_{q'}(s'))}{\partial s'}\right|_{s' = s_j}$ for $j=0,\dots,Q-1$ using Chebyshev differentiation.
\item For each $q$ and each $q'$ evaluate the integral
  $S^{q'}(r_q(s)) = \int_{0}^{1} G_k(r_q(s),r_{q'}(s')) \frac{\partial
    \phi(r_{q'}(s'))}{\partial s'} ds'$ by means of the RP methods
  described in Section~\ref{cheby_old} with $H(r,r') = G_k(r,r')$ and
  on the basis of the quantity
  $f(r_{q'}(s_j)) = \left .\frac{\partial \phi(r_{q'}(s'))}{\partial s'}\right|_{s' = s_j}$
  obtained per point 1, and then obtain the quantity
  $S(r_q(s)) = \sum_{q=1}^M S^{q'}(r_q(s))$.
\item Calculate $T_3[\phi](r_q(s)) = \frac{1}{L_q(s)}\frac{\partial S(r_q(s))}{\partial s}$ using once again Chebyshev differentiation.
\end{enumerate}
The aforementioned Chebyshev-based differentiation algorithm~\cite{10}
numerically obtains the derivative of a given function $g(s)$ defined
in the interval $[0,1]$ as follows:
\begin{itemize}
\item[(a)] Obtain the Chebyshev
  expansion $g(s)  = \sum_{i=0}^{Q-1} c_iT_i(2s-1)$ of order $Q-1$ (throughout this paper we have used the value $Q = 10$ for differentiation purposes);
\item[(b)] Obtain the derivative coefficients $c'_i$ by means of the recurrence relation
  $$c_{i-1}' = c_{i+1}' + 2ic_i \quad (i = Q-1, Q-2, \ldots, 1); \quad
  c_{Q}' = c_{Q-1}' = 0, $$
\item[(c)] Evaluate the desired derivative values via the inverse Chebyshev transform $g'(s)  = \sum_{i=0}^{Q-1} c'_iT_i(2s-1)$.
\end{itemize}

\section{Operator-Regularized and Corner-Regularized Combined Field
  Integral equations\label{op-corn_regul}}
\subsection{Operator-Regularized Combined Field Integral Equation: A Well-Conditioned Smooth-Boundary  Formulation\label{regularized}}
Computational experiments presented in~\cite{1} illustrate the fact
that, while the CFIE formulation is uniquely solvable for all
wavenumbers, solving it accurately can give rise to certain
difficulties. The presence of two tangential derivatives in the
hypersingular operator~\eqref{tangential} associated with the
formulation~\eqref{CFIE} induces arbitrarily large eigenvalues as the
discretizations are refined. As a result, the conditioning of the
system suffers and the numbers of iterations for convergence to a
required tolerance increases significantly.  As shown in~\cite{1}
these difficulties can be effectively resolved, for smooth surfaces
$\Gamma$, by representing the scattered fields in the alternative form
\begin{equation}\label{reg_int_form}
    u = - i\eta\int_\Gamma{G_{k}(r,r')\phi(r') d\ell'} +\int_\Gamma{\frac{\partial G_{k}(r,r')}{\partial{n}(r')}R[\phi(r')] d\ell'}
\end{equation} where 
\begin{equation}\label{continuous_regularization}
    R[\phi](r') = \int_\Gamma{G_{K}(r',r'')\phi(r'') d\ell(r'')}
\end{equation}
is a regularizing operator. The wavenumber $K$ can take any real or imaginary values and does not necessarily need to depend on the wavenumeber of the problem $k$, but~\cite{1}  recommends the use of the value
$$ K = ik.$$ 
Enforcement of the Neumann condition~\eqref{Neumann} then gives rise
to the (Operator) Regularized Combined Field Integral Equation (CFIE-R)
\begin{multline}\label{operator_regularized}
    \frac{i\eta}{2} \phi(r) - i\eta \frac{\partial}{\partial n(r)} \int_{\Gamma} G_{k}(r,r') \phi(r') \, d\ell(r')
+ k^2 \int_{\Gamma} G_{k}(r,r') (n(r) \cdot n(r')) \left( \int_\Gamma{G_{K}(r',r'')\phi(r'') d\ell(r'')} \right)\, d\ell(r') \\
+ \frac{\partial }{\partial \tau(r)}\int_{\Gamma}  G_k(r,r') \frac{\partial}{\partial \tau(r')}\left( \int_\Gamma{G_{K}(r',r'')\phi(r'') d\ell(r'')} \right)\  d\ell(r') 
= -\frac{\partial u^\mathrm{inc}(r)}{\partial n(r)}.
\end{multline}

As illustrated in~\cite{1}, at least for smooth surfaces, this
formulation leads to well conditioned numerical methods which
additionally, if used in conjunction with iterative solvers, require
reduced iteration numbers. But the equation suffers from severe
accuracy and conditioning problems for surfaces $\Gamma$ containing
corners. A ``Corner-Regularized'' discretization is introduced in
Section \ref{corn_reg}; unlike the operator regularization introduced
in the present section, the corner regularization is introduced on the
basis of the parameterizations of the various smooth curves that make
up the piecewise-smooth scatterer $\Gamma$. Since the RP
discretization method~\cite{2,3} utilized in this paper otherwise
relies on use of certain parametrizations and discretizations of
$\Gamma$ in terms of cosine transforms and Chebyshev polynomials, for
conciseness Section~\ref{corn_reg} introduces the
corner-regularization method on the basis of such parametrizations. We
should note, however, that, (i)~The corner regularization method
presented in Section~\ref{corn_reg} could also be implemented in
connection with other Nystr\"om or Boundary-Element discretizations;
and, (ii)~It is expected that, in conjunction with the electromagnetic
operator-regularized equations presented in~\cite{11} (3D Maxwell
analogous of the ones described in Section~\ref{sec:CFIE}), the
methods inherent in the corner-regularized version of of the RP method
presented in Section~\ref{corn_reg} should provide an excellent basis
for a 3D-Maxwell corner-and-edge regularization methodology.

\subsection{Corner-Regularizing Change of Unknown: Regularization of Unbounded Densities \label{corn_reg}}

The proposed corner-regularization approach partitions the domain into
a set of $M$ non-overlapping patches, parameterized by the functions
$r_q = (x_q(s),y_q(s))$, $q = 1, 2, \ldots, M$, each one of which
contains either one corner or no corners, and it is assumed that, by
construction, if the $q$-th patch contains a corner then the corner
point equals $r_q(0)$---that is, the corner point corresponds to the
parameter value $s=0$---a feature that is utilized to eliminate
catastrophic cancellation errors in the context of
Section~\ref{interpolation_header_section}. Note that, in view of this
prescription, the tangential vectors obtained by direct
differentiation of the parametrization $r_q(s)$ on some of the patches
$\Gamma^q$ need to be flipped (multiplied by $-1$) in order to obtain
a consistent prescription (say, counterclockwise) of the tangential
vector field around $\Gamma$.

The regularization method used seeks to yield an integral formulation
whose unknown density is bounded and sufficiently smooth at
corners. To do this, the algorithm relies on the use of a 1-1 change
of variables (CoV) $s = s(\theta)$, $s:[0,1]\to [0,1]$ for each patch
containing a corner~\cite{7,8,9}. The corner CoV gives rise to an
integration Jacobian that vanishes at the corner together with its
derivatives up to a prescribed order, and it maps a uniform mesh in
the $\theta$ variable to a mesh in the $s$ variable that is graded
toward the corner. The CoV Jacobians, which appear in the integrand
multiplying the unknown density, are such that the product of Jacobian
and density has a number of bounded derivatives at the corner
point---thus resulting in a ``regularized'' quantity which can be
approximated closely by the Chebyshev expansions used in the RP
method.

Regularizing CoVs have previously been used for regularization of {\em
  bounded} (Lipschitz continuous) integral densities~\cite{7,8,9}; the
regularization strategy embodied in equations~\eqref{R-CFIE}
and~\eqref{absorbed_ricfie} below, which is applicable to integral
equations whose density solutions are unbounded, is one of the main
enabling elements introduced in this paper. Other essential
algorithmic elements in the proposed approach are introduced in
Sections~\ref{interpolation_header_section}
and~\ref{precomputationswithcov}.

\begin{figure}[ht]
  \begin{center}
    \includegraphics[width=0.4\textwidth]{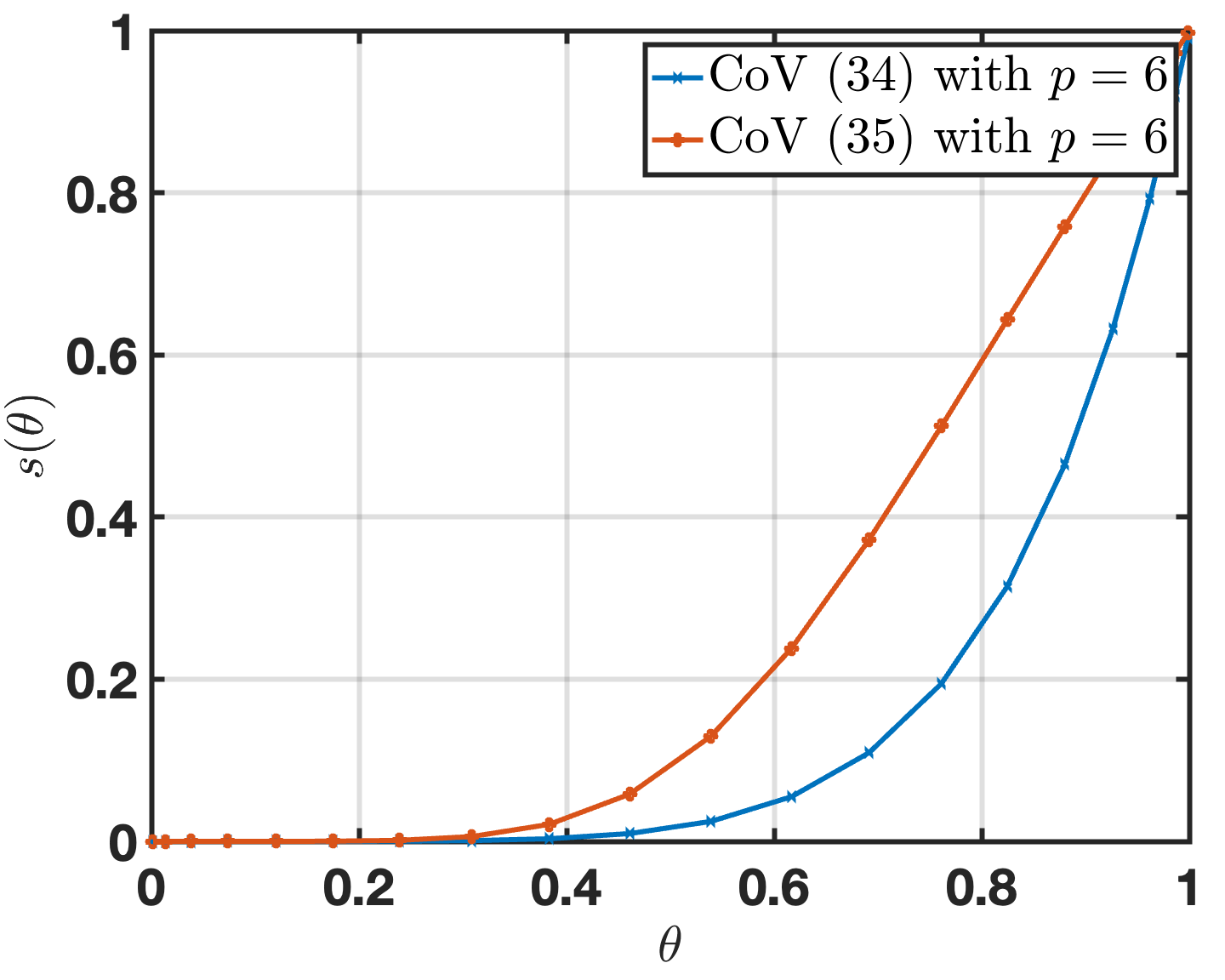}
  \end{center}
  \vspace{-.5cm}
  \caption{Illustration of the graded-mesh Changes of Variables (CoV)  in equations~\eqref{eq:CoVCK} and~\eqref{eq:power_cov} of order $p=6$. \label{fig:ck_cov}}
  \end{figure}

Per equation~\eqref{cov} below the changes of variables utilized in
this paper are based on use of the CoV function~\cite{7,8}
\begin{equation}
  \label{eq:CoVCK}
  w(\theta) = 2\pi \frac{[v(\theta)]^p}{[v(\theta)]^p + [v(2\pi -
    \theta)]^p}, \quad 0 \leq \theta \leq 2\pi,
\end{equation}
which is depicted in
Figure~\ref{fig:ck_cov}, where
$$v(\theta) = \left(\frac{1}{p} - \frac{1}{2}\right) \left(\frac{\pi -
    \theta}{\pi}\right)^3 + \frac{1}{p} \left(\frac{\theta -
    \pi}{\pi}\right) + \frac{1}{2}.$$ The derivatives of the function
$w(\theta)$ of orders less than or equal to $p-1$ vanish at
$\theta =0$ and $\theta =2\pi$---so that, in particular, we have
\begin{equation}
  \label{eq:w-asym}
  w(\theta)\sim \theta^p\quad \mbox{as}\quad  \theta\to 0.
\end{equation}
While simpler changes of
variables such as the monomial CoV~\cite{9}
\begin{equation}
  \label{eq:power_cov}
  s(\theta) = \theta^p,
\end{equation}
$p-1$ of whose derivatives vanish at $\theta =0$ can also be
effective, the CoV~\eqref{eq:CoVCK}, which is used throughout this
paper (with exception of the comparison results presented in
Figure~\ref{fig:results_teardrop} right), gives rise to
significantly better performance than the simpler
version~\eqref{eq:power_cov} on account of the rather evenly split
distribution it enjoys of dicretization points near to and away from
the endpoints of the interval $0 \leq \theta \leq
2\pi$~\cite{7}. Results based on the monomial change of
variables~\eqref{eq:power_cov} presented in
Figure~\ref{fig:results_teardrop} right provide an indication of the
improvements that result from the use of the CoV~\eqref{eq:CoVCK}.

Recalling that, as indicated at the beginning of this section, corners
are assumed to occur only at $s=0$ points, on any one of the patches
$\Gamma^q$ the CoV is induced by the change of
variables~\eqref{eq:CoVCK} in accordance with the expression
% \begin{empheq}[left=\empheqlbrace\ ]{align}
% &  s_q(\theta) = \theta,  \quad \mbox{if $\Gamma^q$ has no corners},\label{cov}\\
% &s_q(\theta) = \frac{1}{\pi} w(\pi
%   \theta),  \quad \mbox{if $\Gamma^q$ has corner at $\theta=0$.}\label{eq:CoV}
% \end{empheq}
\begin{equation}\label{cov}
  \left\{
   \begin{aligned}
&  s_q(\theta) = \theta,  \quad \mbox{if $\Gamma^q$ has no corners},\\
&s_q(\theta) = \frac{1}{\pi} w(\pi
  \theta),  \quad \mbox{if $\Gamma^q$ has a corner at $\theta=0$.}
   \end{aligned}
   \right .
 \end{equation}
Note that the first $p-1$ derivatives of $s(\theta)$ vanish at
$\theta = 0$ for patches containing a corner and
$\frac{ds}{d\theta} = 1$ for patches not containing a corner.

Incorporating the CoV~\eqref{cov} and defining
\begin{equation}\label{elem_th} 
 \widetilde{r}_q(\theta) = r_q(s_q(\theta))\quad\mbox{and}\quad \widetilde{L}_{q'}(\theta') = L_{q'}(s(\theta'))  \frac{ds_{q'}(\theta')}{d\theta'},
\end{equation}
the adjoint double-layer
integral in~\eqref{operator_regularized} (that is, the first integral
on the left-hand side in that equation) with
$r=\widetilde{r}_{q}(\theta)$ becomes
\begin{equation}\label{MFIE_kernel} 
  \int_{0}^{1} \frac{\partial G_{k}(\widetilde{r}_{q}(\theta), \widetilde{r}_{q'}(\theta')}{\partial n(\widetilde{r}_q(\theta))} \phi_{q'}(\widetilde{r}_{q'}(\theta')) \widetilde{L}_{q'}(\theta') \, d\theta'.
\end{equation}
(It is important to recall that, in view of the definitions and
assumptions introduced above in this section, for patches $\Gamma^{q'}$
containing a corner $C$, the corner point is necessarily given by
$C=\widetilde{r}_{q}(0)$.)  The regularization
operator~\eqref{continuous_regularization}, in turn, takes the form
  \begin{equation}\label{regularization_withoutcov}
R[\phi](\widetilde{r}_{q'}(\theta')) =  \sum_{q''=1}^{M}\int_{0}^{1} G_{ik}(\widetilde{r}_{q'}(\theta')),\widetilde{r}_{q''}(\theta'')) \phi_{q''}(\widetilde{r}_{q''}(\theta'')) \widetilde{L}_{q''}(\theta'') \, d\theta''.
\end{equation}
Similarly, for the last two integrals
in~\eqref{operator_regularized} we obtain the expressions
\begin{equation}\label{EFIE_firstkernel}
\int_{0}^{1} G_{k}(\widetilde{r}_{q}(\theta),\widetilde{r}_{q'}(\theta') (n(\widetilde{r}_{q}(\theta)) \cdot n(\widetilde{r}_{q'}(\theta')) \widetilde{L}_{q'}(\theta') R[\phi](\widetilde{r}_{q'}(\theta'))\, d\theta' 
\end{equation}
and
\begin{equation}\label{EFIE_secondkernel}
\frac{1}{\widetilde{L}_q(\theta)}\frac{\partial}{\partial \theta}\int_{0}^{1}G_{k}(\widetilde{r}_{q}(\theta),\widetilde{r}_{q'}(\theta')) \frac{\partial }{\partial \theta'} R[\phi](\widetilde{r}_{q'}(\theta'))\, d\theta',
\end{equation}
respectively.

In all, incorporating  the CoV~\eqref{cov} in~\eqref{operator_regularized}, using the expressions~\eqref{MFIE_kernel}
through~\eqref{EFIE_secondkernel}, and
multiplying the result by $\widetilde{L}_{q}(\theta)$, we obtain
\begin{multline}\label{R-CFIE}
    \frac{i\eta}{2} \phi(\widetilde{r}_{q}(\theta))\widetilde{L}_{q}(\theta) - i\eta \widetilde{L}_{q}(\theta) \sum_{q'=1}^{M} \int_{0}^{1} \frac{dG_{k}(\widetilde{r}_{q}(\theta), \widetilde{r}_{q'}(\theta')}{\partial n(\widetilde{r}_{q}(\theta))} \phi(\widetilde{r}_{q'}(\theta') \widetilde{L}_{q'}(\theta') \, d\theta'  \\
+ {k}^2  \widetilde{L}_{q}(\theta)\sum_{q'=1}^{M}\int_{0}^{1} G_{k}(\widetilde{r}_{q}(\theta),\widetilde{r}_{q'}(\theta')) (n(\widetilde{r}_{q}(\theta)) \cdot n(\widetilde{r}_{q'}(\theta')) \widetilde{L}_{q'}(\theta') R[\phi](\widetilde{r}_{q'}(\theta'))\, d\theta' \\
+ \frac{\partial }{\partial \theta}\sum_{q'=1}^{M}\int_{0}^{1} G_{k}(\widetilde{r}_{q}(\theta),\widetilde{r}_{q'}(\theta')) \frac{\partial }{\partial \theta'} R[\phi](\widetilde{r}_{q'}(\theta'))\, d\theta'
= -\widetilde{L}_{q}(\theta)\frac{\partial u^{\text{inc}}(\widetilde{r}_{q}(\theta))}{\partial n(\widetilde{r}_{q}(\theta))}.
\end{multline}
Thus, defining the new unknown
\begin{equation}
  \label{eq:ch_of_unk}
  \psi_q(\theta) = \phi(s_q(\theta)) \widetilde{L}_{q}(\theta),
\end{equation}
and expressing the regularization operator \eqref{regularization_withoutcov} in the form
\begin{equation}\label{regularization_withcov}
\widetilde{R}[\psi_q](r) =  \sum_{q'=1}^{M}\int_{0}^{1} G_{ik}(r,\widetilde{r}_{q'}(\theta'))\, \psi_q(\theta') \, d\theta',
\end{equation}
we finally obtain the {\em change-of-unknown, Corner-Regularized} version
\begin{multline}\label{absorbed_ricfie}
  \frac{i\eta}{2} \psi_q(\theta) - i\eta \widetilde{L}_q(\theta)\sum_{q'=1}^{M}\int_{0}^{1} \frac{\partial G_{k}(\widetilde{r}_{q}(\theta), \widetilde{r}_{q'}(\theta')}{\partial n(\widetilde{r}_{q}(\theta))} \psi_q(\theta') \, d\theta' \\
  + {k}^2  \widetilde{L}_q(\theta)\sum_{q'=1}^{M}\int_{0}^{1} G_{k}(\widetilde{r}_{q}(\theta),\widetilde{r}_{q'}(\theta')) (n(\widetilde{r}_{q}(\theta)) \cdot n(\widetilde{r}_{q'}(\theta'))\widetilde{L}_{q'}(\theta') R[\psi](\widetilde{r}_{q'}(\theta))\, d\theta' \\
  + \frac{\partial }{\partial \theta}\sum_{q'=1}^{M}\int_{0}^{1} G_{k}(\widetilde{r}_{q}(\theta),\widetilde{r}_{q'}(\theta')) \frac{\partial }{\partial \theta'} R[\psi](\widetilde{r}_{q'}(\theta'))\, d\theta' = -\frac{\partial u^{\text{inc}}}{\partial n(\widetilde{r}_{q}(\theta))}
   \widetilde{L}_q(\theta)
\end{multline}
of the Operator Regularized R-CFIE
equation~\eqref{operator_regularized}.  It is important to emphasize
that $\psi(\theta)$ is a much smoother function than $\phi$ (and,
thus, significantly easier to approximate by polynomials) on account
of the multiplicative line element $\widetilde{L}_q(\theta)$ factor
contained in $\psi$, whose value together with that of $p-1$ of its
derivatives, vanishes at the corner point $\theta=0$ where the field
singularity is located. This is another key aspect of the proposed
methodology for scattering by domains with corners.
 
For the purpose of comparison in the numerical results section we
additionally consider the following change-of-unknown,
corner-regularized version of the MFIE equation~\eqref{MFIE}:
\begin{equation}\label{absorbed_mfie}
  -\frac{1}{2} \psi_q(\theta) + \widetilde{L}_q(\theta)\sum_{q'=1}^{M}\int_{0}^{1} \frac{\partial G_{k}(\widetilde{r}_{q}(\theta), \widetilde{r}_{q'}(\theta')}{\partial n(\widetilde{r}_{q}(\theta))} \psi_{q'}(\theta') \, d\theta = -\frac{\partial u^{\text{inc}}}{\partial n(\widetilde{r}_{q}(\theta))}
  \widetilde{L}_q(\theta).
\end{equation}

Note that in light of the introduction of the new unknown
$\psi(\theta)$, which includes a the $\widetilde{L}_q(\theta)$ factor,
the rectangular polar discretization method needs to be modified as follows:
\begin{enumerate}
\item Substitute the Fej\'er evaluation of the far
  interactions~\eqref{far_nocorner} by the corresponding Fej\'er
  approximation
\begin{equation}\label{far_nocorner2}
  \mathcal{I}^{q'} [H,\widetilde{f}](r) \approx \sum_{j=0}^{N-1} H (r, \widetilde{r}_{q'}(\theta_j))\widetilde{f}((\theta_j)) \omega_j , \quad \mbox{for } \quad r\quad \mbox{away from} \quad \Gamma^{q'}.
\end{equation}
\item Substitute~\eqref{density_cheby_expansion} by the corresponding expansion
  \begin{equation}\label{density_cheby_expansion2}
  \widetilde{f}_{q'}(\theta') \approx \sum_{m=0}^{N-1} \widetilde{a}_m^{q'} T_m
  (2\theta'-1)
\end{equation}
(where, in analogy to Section~\ref{sec:CFIE}, here we set either
$\widetilde{f}_q(\theta')=\psi_q(\theta')$ or
$\widetilde{f}_q(\theta')=\frac{\partial \psi_q}{\partial
  \theta}(\theta')$).
\item Replace the precomputation integrals~\eqref{precomp_def} by
\begin{equation}\label{precomp_def2}
  \widetilde{\beta}_m^{q'}[H](r) = \int_{0}^{1}H(r,\widetilde{r}_{q'}(\theta')) \, \, T_m(2\theta'-1) \, d\theta',
\end{equation}
for the kernels $H$ listed near the end of Section~\ref{sec:CFIE}, and
for target discretization points $r = \widetilde{r}_{q'}(\theta_j)$ on
the source patch $\Gamma^{q'}$, as well as for near singular
target discretization points $r = \widetilde{r}_{q}(\theta_j)$ on
nearby patches $\Gamma^{q}$. (Note that, in particular, the
$\widetilde{L}_{q'}(\theta')$ factor, which is contained in the new
unknown $\psi(\theta')$, does not appear in the integrand of
equation~\eqref{precomp_def2}.)  In view of the singularities introduced by the corner
geometries under consideration, certain significant variations of the
rectangular polar method~\cite{2,3} are needed, which are described in
Sections~\ref{interpolation_header_section} and~\ref{precomputationswithcov}.
\item Using the precomputed weights~\eqref{precomp_def2} evaluate the
  integral
\begin{equation}\label{precomp_sum2}
\mathcal{I}^{q'}[H,\widetilde{f}](r) = \sum_{m=0}^{N-1} \widetilde{a}_m^{q'} \widetilde{\beta}_m^{q'}[H](r)
\end{equation}
in terms of the coefficients $\widetilde{a}^{q'}_m$
in~\eqref{density_cheby_expansion2}, and evaluate the needed overall
integrals by a sum analogous to~\eqref{eq:sum-int}. The tangential
differentiations needed for the evaluation of the operator $T_3$
in~\eqref{gen_int_dif} are produced by a procedure analogous to the
one described in Section~\ref{numerical_derivative} applied in the
$\theta'$ variable (instead of the $s'$ variable used in that
section).
\end{enumerate}

\subsection{Finite precision evaluation of Green-function source-target interactions}\label{interpolation_header_section}
The introduction of the CoV and associated change of unknown, which
effectively enables the approximation of the new integral densities in
terms of corresponding Chebyshev expansions, can also give rise to
significant numerical cancellation errors unless suitable numerical
procedures are used, on account of the finite precision inherent in
computer arithmetic. The strategies we use to eliminate cancellation
errors for ``near-corner interactions'' and ``self interactions'' are
described in what follows. Note that, in view of the overall strategy
used, source-target near interactions only arise as part of the
precomputation stage described in
Section~\ref{precomputationswithcov}.

Let us consider a corner point $C$ located at the intersection of the
patches $\Gamma^q$ and $\Gamma^{q'}$, $C = \Gamma^q\cap\Gamma^{q'}$,
for which, in accordance with the conventions put forth in
Section~\ref{corn_reg}, we have
$C = \widetilde{r}_{q}(0) = \widetilde{r}_{q'}(0)$. Clearly, the CoV
introduced in that section generally induces graded meshes that can
result in marked clustering of points near $C$ on both the source and
target patches ($\Gamma^{q'}$ and $\Gamma^q$ respectively). The
resulting extreme proximity of source and target points (generically)
gives rise to catastrophic cancellation in the evaluation of certain
source-target differences
$|\widetilde{r}_{q}(\theta) - \widetilde{r}_{q'}(\theta')|$ that are
required for evaluation of Green-function values. In detail,
catastrophic cancellations occur for source and target points
$\widetilde{r}_{q}(\theta)$ and $\widetilde{r}_{q'}(\theta')$ near $C$
if at least one of the coordinates of
$C = (\widetilde{x}_{q}(0),\widetilde{y}_{q}(0)) =
(\widetilde{x}_{q'}(0),\widetilde{y}_{q'}(0))$ is an $\mathcal{O}(1)$
quantity. Indeed, if a source discretization point
$\widetilde{r}_{q}(\theta) =
(\widetilde{x}_{q}(\theta),\widetilde{y}_{q}(\theta))\in\Gamma_{q}$
and a target discretization point
$\widetilde{r}_{q'} (\theta')=
(\widetilde{x}_{q'}(\theta'),\widetilde{y}_{q'}(\theta'))\in\Gamma_{q'}$
are each within a distance $\delta$ from the corner, where
$\delta \ll|\widetilde{r}_{q}(0)| = |\widetilde{r}_{q'}(0)|$, the
resulting computation of the difference
$|\widetilde{r}_{q}(\theta)-\widetilde{r}_{q'}(\theta')|$ would suffer
from a significant cancellation error: a loss of a number
$\approx |\log_{10}\delta|$ of digits would ensue from such a
calculation. Taking into account the small $\delta$ values that may
arise in this context (e.g., the numerical examples presented in this
paper feature $\delta$ values in the range
$ 10^{-14}\leq\delta\leq 10^{-10}$ for $p=4$ and
$ 10^{-24}\leq\delta\leq 10^{-15}$ for $p=6$ with
$|r_{q}(\theta)| =\mathcal {O}(1)$,
$|r_{q'}(\theta')| =\mathcal {O}(1)$) the accuracy loss can be highly
detrimental.  Accurately computing these small quantities (and
associated large Hankel-function values) is crucial to produce the
necessary and correct $\mathcal{O}(1)$ contributions to the values of
the solution both near and away from the corner point.

Large source-target cancellation errors arise across smooth portions
of the surface $\Gamma$ as well, in view of the need for evaluating
the normal derivative $\frac{\partial G_k(r, r')}{\partial n(r)}$ of
the Green function $G_k$ in such regions. Cancellation errors also
affect the Green function itself, which also appears in
equations~\eqref{T2} and~\eqref{T3}, but the cancellation-error
difficulty is particularly problematic in the context of the adjoint
double-layer operator~\eqref{T1} since the normal derivative contains
the term
\begin{equation}
  \label{factor}
  F(r,r') =  \frac{(r-r')\cdot
    n(r)}{|r-r'|^2}.
\end{equation}
The issue stems from the quadratic decay of the numerator
$n(r) \cdot (r - r')$ in~\eqref{factor}, whose magnitude decreases
like $\mathcal{O}(|r - r'|^2)$ as $|r - r'| \to 0$, at points
$r \in \Gamma$ where the surface $\Gamma$ is smooth, and whose
resulting (large) relative errors are then amplified upon division by
the small denominator $|r-r'|^2$. Fortunately such errors can be
effectively treated by using a Taylor expansion of the the
parametrization $r = \widetilde{r}_{q'}(\theta')$ around
$\theta' = \theta$. Neighboring pairs of points
$\widetilde{r}_{q}(\theta)$ and $\widetilde{r}_{q'}(\theta')$ that are
additionally near corners, whether they are on the same patch ($q=q'$)
or on different patches ($q\ne q'$), require a somewhat different
treament, in view of the corner changes of variables and associated
fine refinements utilized by the proposed algorithm---which may
compound the degree of proximity, but which would be cumbersome to
treat via Taylor expansions on account of the multiple-derivative
vanishing at corners that is induced by the corner changes of
variables parametrization. To tackle the near corner cancellation
problem around a corner point $C$ we consider the parametrizations
$\widetilde{r}_{q}(\theta)$ and $\widetilde{r}_{q'}(\theta')$ of the
two patches that contain the point $C$ (possibly with $q=q'$)---so
that, in accordance with the prescriptions in Section~\ref{corn_reg}
we have
$C = r_q(0)=r_{q'}(0)=\widetilde{r}_{q}(0)=\widetilde{r}_{q'}(0)$. In
view of the first equation in~\eqref{elem_th} our algorithm
incorporates the lowest order Taylor expansions (with remainder)
$r_q(s) = r_q(0) + s R_q(s)$ and $r_{q'}(s) = r_{q'}(0) + s R_q(s)$ so
that
$\widetilde{r}_{q}(\theta)-\widetilde{r}_{q'}(\theta') =
r_q(s_q(\theta)) - r_{q'}(s_{q'}(\theta')) =
s_q(\theta)R_q(s_q(\theta)) -
s_{q'}(\theta')R_{q'}(s_{q'}(\theta'))$. Since in
view~\eqref{eq:w-asym} and~\eqref{cov} of the $s_q$ changes of
variables satisfy $s_q(\theta)\sim \theta^p$ and
$s_{q'}(\theta')\sim (\theta')^p$, this procedure produces
$\widetilde{r}_{q}(\theta)-\widetilde{r}_{q'}(\theta')$ as the
difference of two suitably small numbers thereby bypassing the
near-corner cancellation problem. (Note that in the case $q=q'$ the
algorithm additional incoporates the aforementioned Taylor expansion
in $\theta'$ around $\theta' =\theta$, thereby leading to accurate
evaluation of small differences throughout $\Gamma$.)

\begin{figure}[!ht] % example dataset
  \centering% <-- superfluous in this example as commented by Zarko
  \includegraphics[width=0.4\textwidth]
  {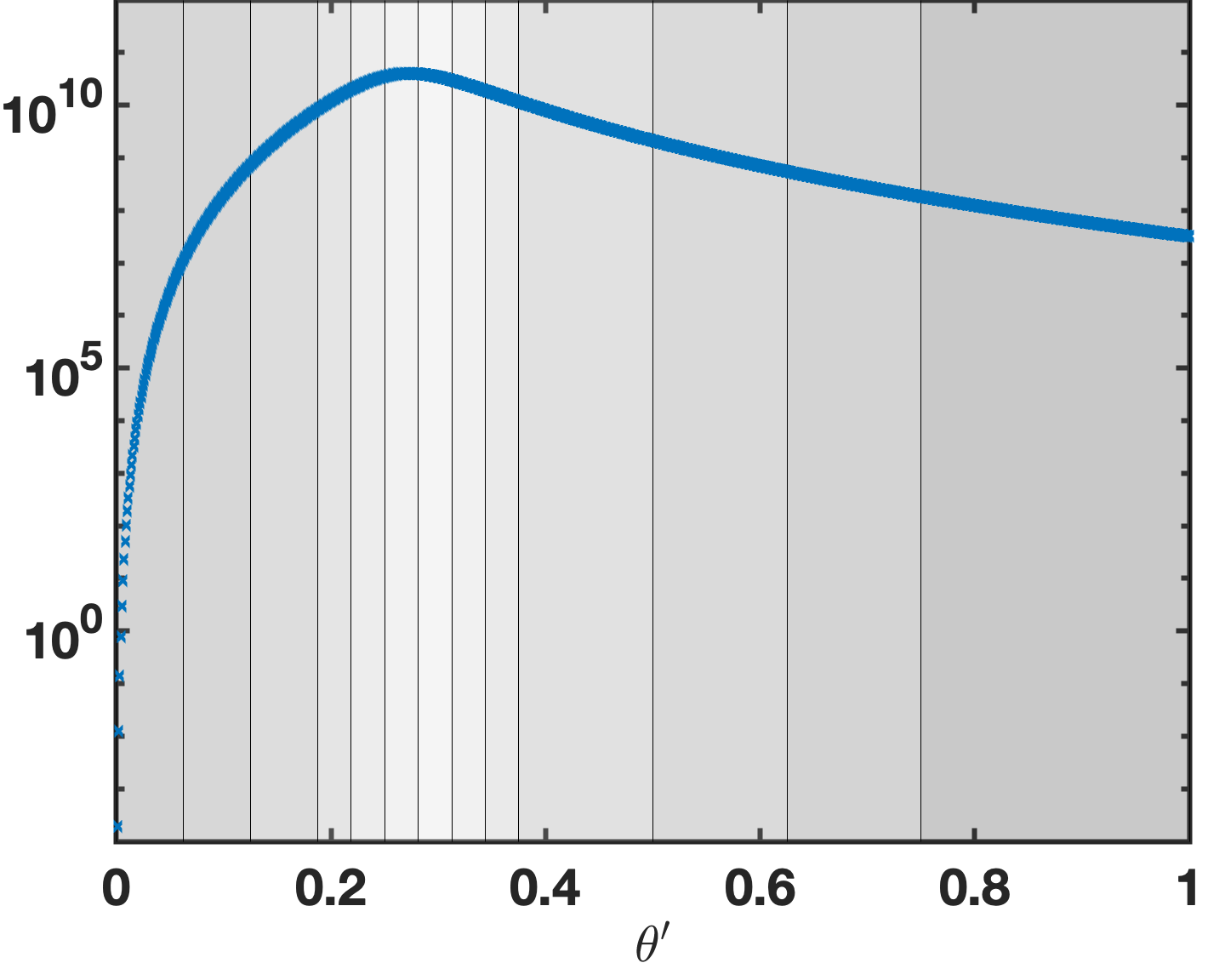}\hspace{0.8cm}
  \includegraphics[width=0.4\textwidth]
  {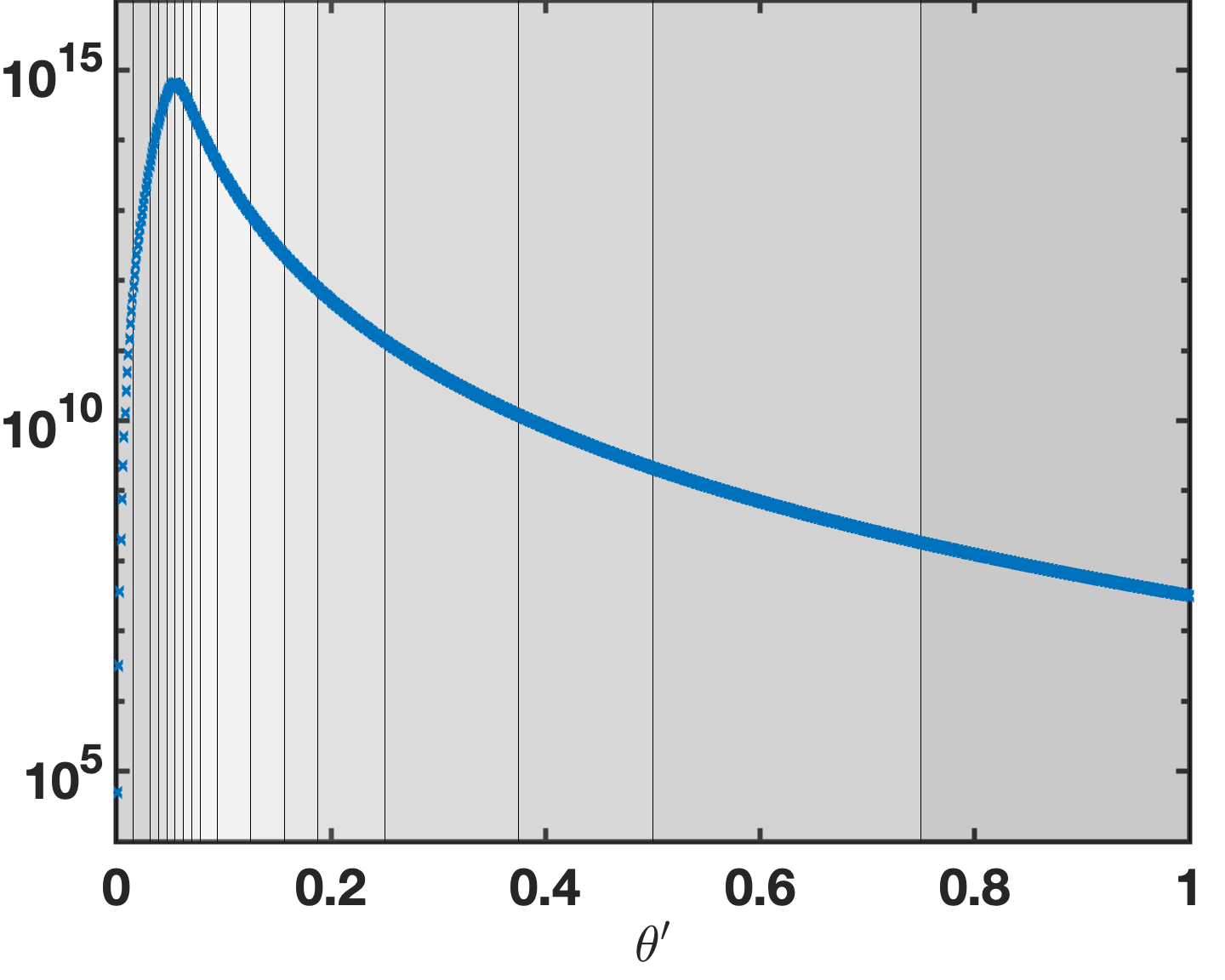}% <-- No space here
  \caption{The blue crosses display the magnitude of the kernel
    $\frac{\partial G_{k}(\widetilde{r}_{q}(\theta),
      \widetilde{r}_{q'}(\theta')}{\partial
      n(\widetilde{r}_{q}(\theta))}$ in the first integral
    in~\eqref{absorbed_ricfie}, as a function of $\theta'$, for values
    of $\theta$ corresponding to target points at distances
    $1.98\times10^{-12}$ (left) and $1.23\times10^{-16}$ (right) away
    from the corner (target point and source patch are on different
    sides of a corner).  The vertical lines correspond to abscissas
    selected by the Gauss-Kronrod adaptive quadrature algorithm. }
\label{fig:mfiekernels}.  
\end{figure}

\begin{figure}[!ht] % example dataset
  \centering% <-- superfluous in this example as commented by Zarko
  \begin{center}
    \includegraphics[width=0.4\textwidth]
    {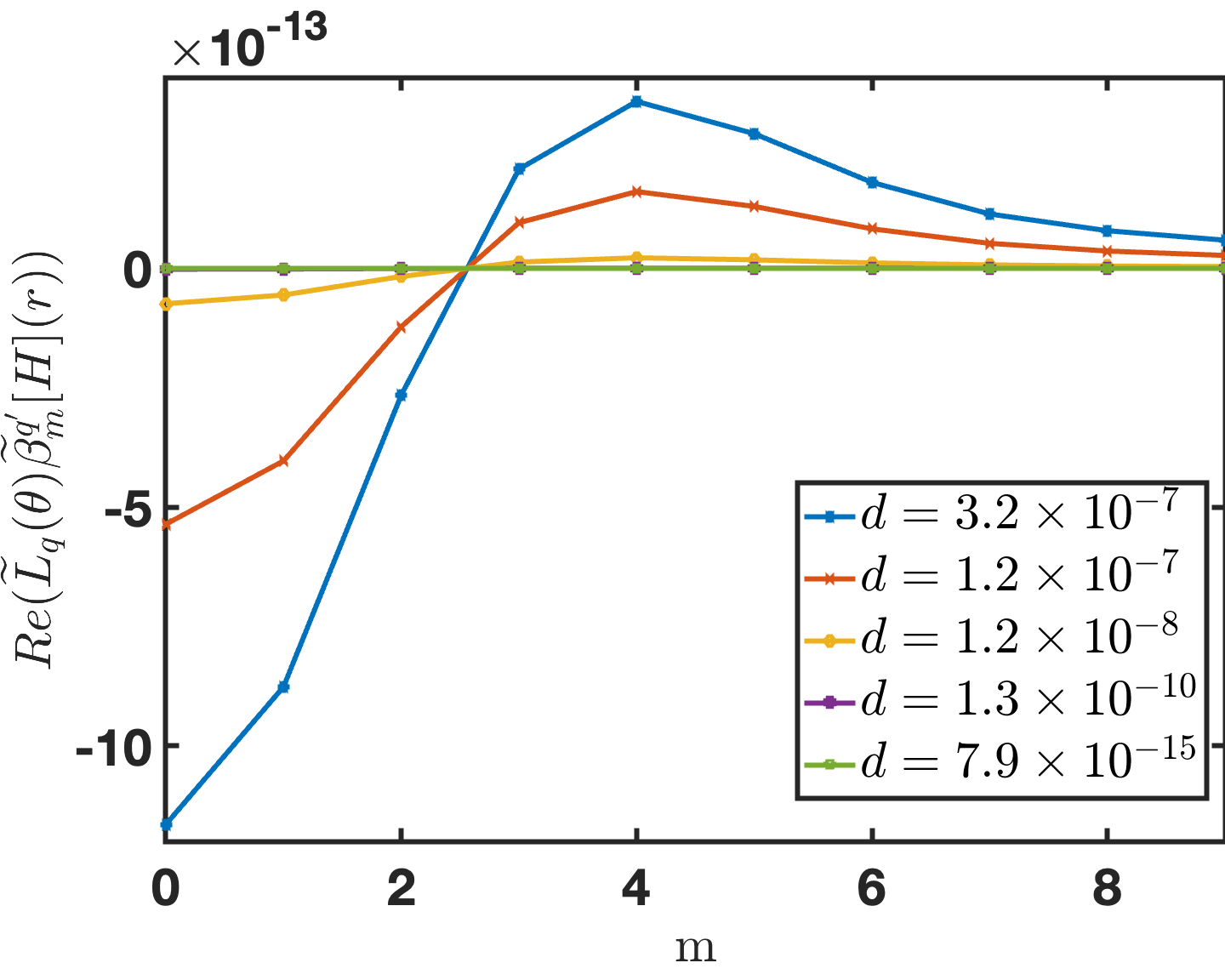}\hspace{0.8cm}
  \includegraphics[width=0.4\textwidth]
  {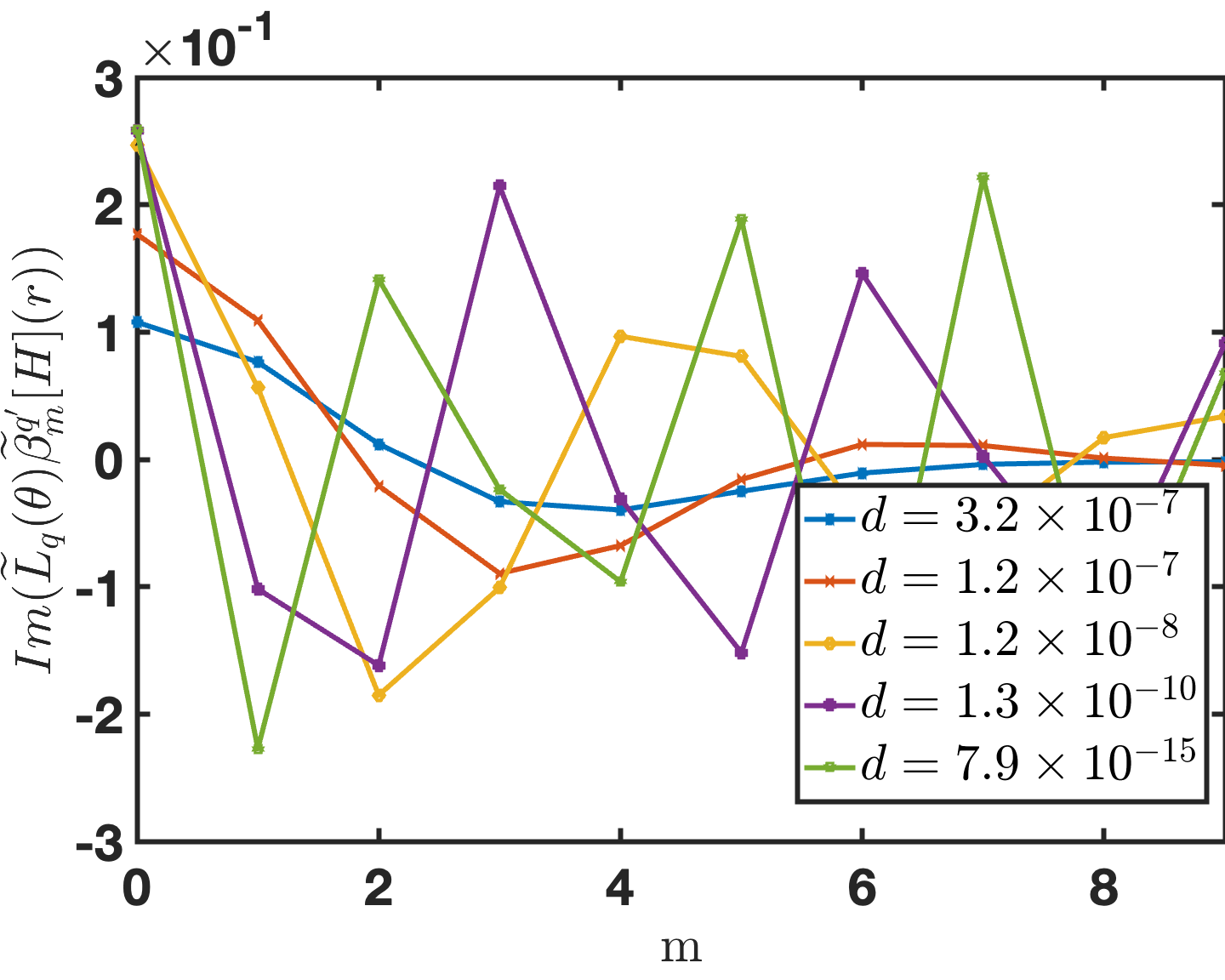}
\end{center}
\vspace{-0.5cm}
\caption{  The real and imaginary parts of the quantities $\widetilde{L}_q(\theta)\widetilde{\beta}_m^{q'}[H](r)$ are shown, at points $ r $ located at various distances $ d $ from a corner $ C $ of a square scatterer, where these quantities are obtained as the product of the line element $\widetilde{L}_q(\theta)$ and the precomputed weights~\eqref{precomp_def2} for the kernel $ H = \frac{\partial G_k(r, r')}{\partial n_t(r)} $ associated with the corner-regularized formulation. These values are relevant for implementing the integral formulation~\eqref{absorbed_ricfie}.  The results are shown for the $ q' $-th source patch and for target points $ r $ on the $ q $-th observation patch ($ q \neq q' $), where the patches are positioned on opposite sides of the corner $ C $. As illustrated, the precomputed values remain bounded throughout the $ q $-th patch for all Chebyshev polynomial orders $ m $.}
\label{fig:precomputations}
\end{figure}
\subsection{Numerical treatment of pre-computations}\label{precomputationswithcov}

The rectangular-polar method introduced in Section~\ref{cheby_old}
requires the precomputation of the integration weights
$\beta^{q'}_{m}$ in equation~\eqref{precomp_def} for each one of the
kernels $H(r,r') = \frac{\partial G_{k}(r,r')}{\partial n(r)}$,
$H(r,r') = G_{k}(r,r') (n(r) \cdot n(r'))$ and
$H(r,r') = G_{k}(r,r')$, and for the singular and near-singular
integration cases considered in Section~\ref{cheby_old}, namely, for
target points $r$ that are either within or close to the integration
patch, respectively---and for which, in particular, the algorithm
relies on the techniques introduced in the
Section~\ref{interpolation_header_section} for accurate computation of
source-target distances. In previous implementations~\cite{2,3} a
change of variables was utilized as a singularity cancellation
technique, upon which a high-order Fej\'er quadrature integration rule
may be employed to precompute the integration weights $\beta^{q'}_{m}$
with high accuracy. In the particular situation at hand, however, the
kernel $H(r,r') = \frac{\partial G_{k}(r,r')}{\partial n(r)}$, which
is actually bounded for $r'$ at the corner and $r$ near the corner,
still tends to infinity like $1/|r-r'|$, as illustrated in
Figure~\ref{fig:mfiekernels}, as both $r$ and $r'$ approach the corner
while $(r-r')$ maintains an angle bounded away from $\pi/2$ with the
normal $n(r)$---a difficulty that is not shared by the kernels
$H(r,r') = G_{k}(r,r') (n(r) \cdot n(r'))$ and
$H(r,r') = G_{k}(r,r')$, which tend to infinity only logarithmically
and are thus natrually resolved by the rapid vanishing of the density
$\psi(r')$ as $r'$ tends to the corner. (In contrast, the
$1/|r-r'|\sim (\theta')^{-p}$ asymptotics of the normal-derivative
kernel are not cancelled by the near vanishing of the density
$\psi_q(\theta') \sim (\theta')^{(p-1) - \nu p}$, where $\nu$ denotes
the corner singularity exponent: $\psi(r)\sim d^{-\nu}$ as the
distance $d$ from $r$ to the corner tends to zero.) As illustrated in
Figure~\ref{fig:mfiekernels} (and as it is easy to appreciate by
considering the geometry at hand), this peak corresponds to a
different point in each case related both to the distance of the
target point to the corner as well as the corner angle
itself. Although, in principle, the location of this peak could be
computed analytically and a similar additional singularity
annihilating CoV could be incorporated to allow the successful use of
Fej\'er quadrature for these scenarios as well, the implementation
proposed in this paper relies instead on the adaptive Gauss-Kronrod
(GK) quadrature scheme~\cite{davis2007methods}. Since the GK scheme
determines the optimal placement of integration points adaptively,
this approach offers the advantage of simplicity (at the cost of a
slight potential computational overhead and, in the current
implementation, an error floor dictated by Matlab's $10^{-12}$ GK
error limit; see Remark~\ref{rem-11-dig}). Figure~\ref{fig:precomputations} displays values of precomputed integration weights; clearly, as suggested by the figure, the precomputed values are bounded quantities throughout the observation patch for all Chebyshev polynomial orders $m$.

\section{Numerical Results\label{numer}}

This section illustrates the character and performance of the proposed
algorithm by means of a variety of numerical results, including,
1)~Quantitative evaluation of the algorithm's performance as each one
of the new key innovations is omitted
(Section~\ref{comparisonagainstRP}); 2)~Evaluation of the conditioning
of the linear systems and the numbers of GMRES iterations resulting
from the proposed regularized R-CFIE formulation and comparison with
the corresponding metrics resulting from the (unregularized) MFIE
formulation (Section~\ref{conditioning}); 3)~Convergence studies for
computed field values at points extremely close to corners
(Section~\ref{accuracy}); 4)~Convergence studies for the integral
density $\psi$ both near and away from corners
(Section~\ref{exponent}); and, 5)~Evaluation of the singular exponent
of the numerically computed solution density and comparison with the
exponents based on theoretical considerations
(Section~\ref{exponent}). For uniformity, results of all such studies
are presented for a cylindrical geometry with a square
cross-section. To demonstrate the generality of the method, field
convergence studies were conducted using scatterers of varying
geometries, including parallelograms and teardrop-shaped objects with
a wide range of interior angles. Specifically, the parallelograms
included shapes with interior angles of $45^\circ$ and
$0.57^\circ (=0.01$ radians), while the teardrop-shaped scatterers feature curved
boundaries with interior corner angles of $90^\circ$ and $9^\circ$,
respectively.  All of the numerical errors reported in this section
where obtained by means of the relative error expression
\begin{equation}\label{epsilon}
  \varepsilon = |u - u_\mathrm{ref}|/|u_\mathrm{max}|,
\end{equation}
where $|u_{\mathrm{max}}|$ denotes the maximum magnitude of the total
field $u$ outside the scatterer. In all cases the reference solution
$u_\mathrm{ref}$ was obtained as the numerical solution on a fine
discretization, containing, in each case, at least twice as many
discretization points as the finest discretization
shown. Additionally, comparisons with solutions produced by a
commercial solver are presented in Section~\ref{accuracy}
(Figure~\ref{fig:lumericalvsbie}). All the results presented were
obtained on the basis of Matlab implementation of the proposed
algorithm running on an Apple Mac Mini M2 2023 personal
computer. Computing times are not reported on account of the
significant overhead associated with the Matlab implementation
used. However, in view of the character of the algorithm it is
expected that the computational cost required by an efficient
implementation for a given number of discretization points should be
similar to that obtained for related unaccelerated Nystr\"om
implementations~\cite{2,3} for the same number of points, and they
could of course be further accelerated by means of $N\log N$-type
algorithmic acceleration of the type considered
in~\cite{bauinger_bruno} and references therein.

\subsection{Quantitative Evaluation with Key Improvements Omitted}\label{comparisonagainstRP}
The scattering-solver algorithm for domains with corners described in
previous sections relies on four main elements, namely, (i)~The
regularized integral formulation~\eqref{reg_int_form} and associated
regularized integral equation~\eqref{operator_regularized} (which,
while uniquely solvable, gives rise to unbounded densities $\phi$ and
to integral kernels that are not weakly singular around corner
points); (ii)~The regularizing changes of variables~\eqref{eq:CoVCK},
as described in Section~\ref{corn_reg}; as well as, (iii)~The change
of unknown~\eqref{eq:ch_of_unk}; and, (iv)~The precision-preserving
methodologies described in
Section~\ref{interpolation_header_section}. In order to highlight the
beneficial effects provided by these strategies, in this section we
compare the performance of the proposed solver to the performance that
results as one or more of the proposed techniques (ii), (iii) and~(iv)
are excluded from the proposed algorithm; the benefits arising from
point~(i) are then explored in Section~\ref{conditioning}.

\begin{figure}[!ht]
  \begin{center}
    \includegraphics[width=0.3\textwidth]{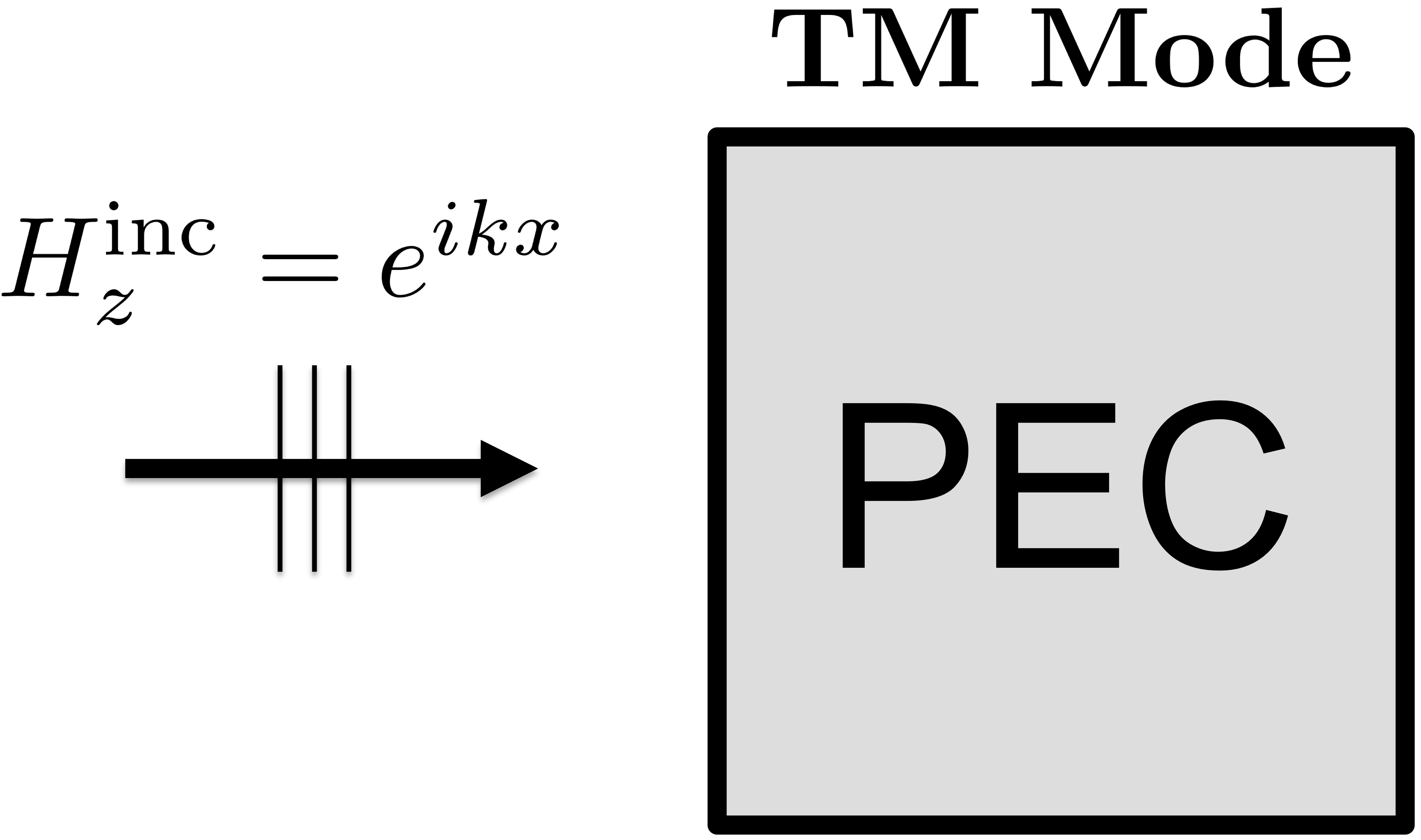}
  \end{center}
  \caption{Plane wave incident upon a PEC cylinder of square cross-section.}
  \label{fig:excitation}
\end{figure}

\begin{figure}[!ht]
  \begin{center}
    \centering% <-- superfluous in this example as commented by Zarko
      \includegraphics[width=0.43\textwidth]
      {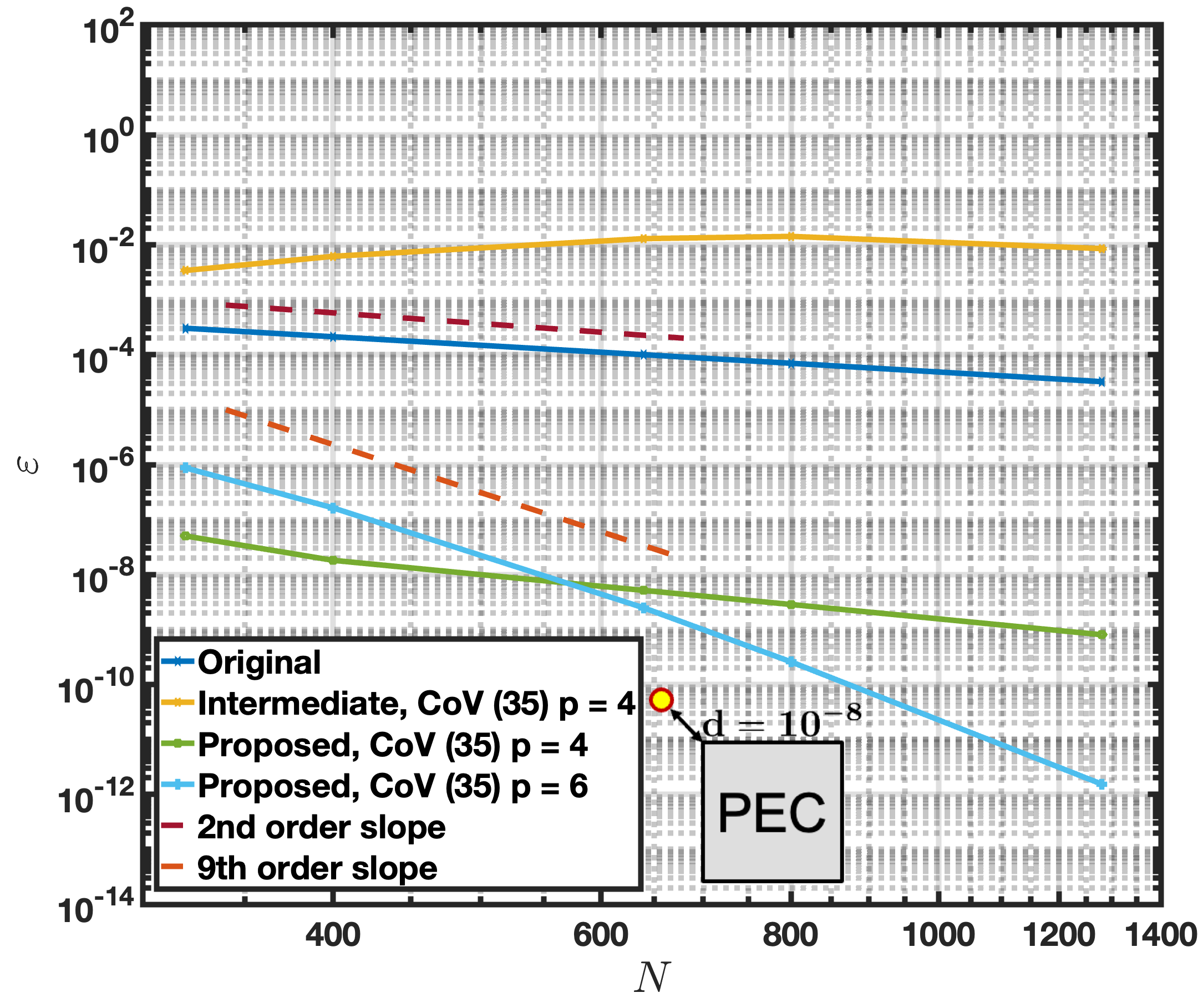} \hspace{0.6cm}
      \includegraphics[width=0.43\textwidth]
      {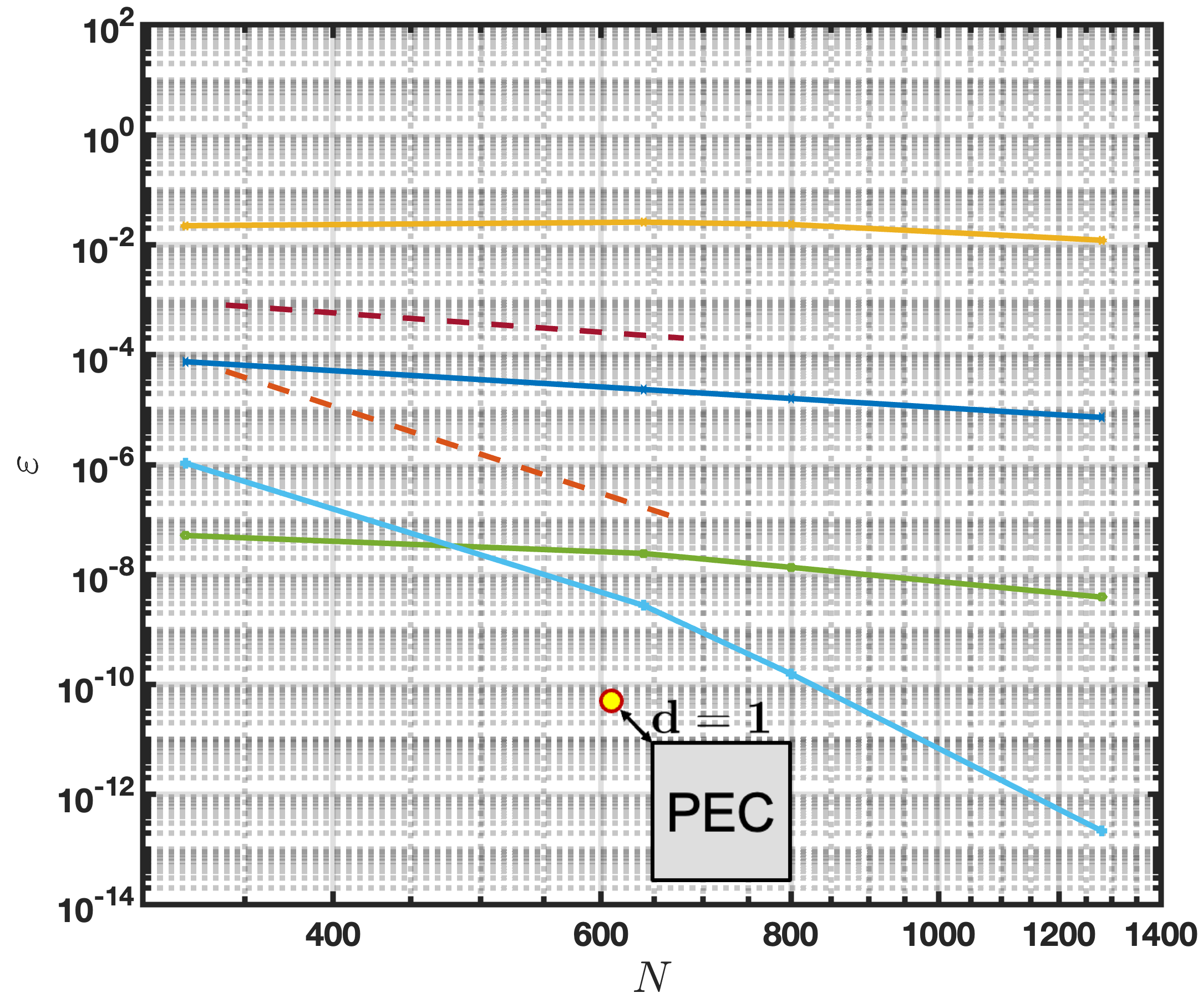}% <-- No space here
    \end{center}
    \vspace{-0.5cm}
    \caption{Error $\varepsilon$ (equation~\eqref{epsilon}) as a
      function of $N$ (equation~\eqref{eq:N}) in the solution $u$ at
      distances of $d=10^{-8}$ (left) and $d=1$ (right) from a corner
      point, obtained by means of implementations with key
      improvements omitted, as described in
      Section~\ref{comparisonagainstRP}, compared to results produced
      by the full proposed method.}
  \label{fig:square_allmethodcomparison}
\end{figure}

In detail this section present numerical results obtained on the basis
of the following formulations:
\begin{itemize}
\item[--] Original formulation: R-CFIE with all additional key
  improvements (ii), (iii) and (iv) omitted;
%\item[--] Intermediate~1 formulation: R-CFIE formulation with the
  %addition of the regularizing CoV mentioned in point~(ii) above
  %(without change of unknown);
\item[--] Intermediate formulation: R-CFIE formulation with the addition of the regularizing CoV mentioned in point~(ii) above (without change of unknown) and incorporating additionally the precision-preserving methodologies mentioned in points~(ii) and~(iv);
%\item[--] Intermediate~3 formulation: R-CFIE formulation with the
  %addition of the regularizing CoV and change of unknown methodologies
  %mentioned in points~(ii) and~(iii), without use of the
  %precision-preserving methodologies in point~(iv);
\item[--] Proposed formulation: R-CFIE formulation with the addition of
  the regularizing CoV, and the change of unknown and
  precision-preserving methodologies mentioned in points~(ii)--(iv).
\end{itemize}
Figure~\ref{fig:square_allmethodcomparison} presents numerical
convergence results obtained by means of each one of these
formulations for the problem of scattering by a PEC cylinder with a
square cross section of side $a = 2$ that is illuminated by a plane
wave excitation travelling in the positive $x$ direction, as depicted
in Figure~\ref{fig:excitation}, with $k = 10$ (that is
$\lambda = 0.628$ and $a \approx 3.2\lambda$). The scatterer is
discretized into a finite number of patches, each one of which
contains 10 discretization points. Error values at two points
diagonally away from the top left corner of the scatterer are
displayed in the figure, at distances $d=d_1= 1\times10^{-8}$ (left)
and $d=d_2=1$ (right) away from the corner, respectively, as indicated
in the figure. The figure suggests that, as verified by means of
observation points at varied distances, including much closer
distances than the one cited above, the method is capable of producing
high accuracies with similar numbers of discretization points both
near and away from the corner point.
\begin{remark}\label{rem-11-dig}
  The 12-digit accuracy limit depicted in
  Figure~\ref{fig:square_allmethodcomparison} and other error plots
  presented in this paper originates from the precision of the
  Gauss-Kronrod integration algorithm employed near corners, which is
  constrained to 12-digit accuracy in the Matlab implementation we have
  used. It is expected that use of a higher accuracy adaptive
  quadrature rule near corners would eliminate this limitation.
\end{remark}

Figure~\ref{fig:square_allmethodcomparison} highlights the impact of
each one the ``key improvements'' mentioned above in this
section. Indeed, the figure illustrates the character of the proposed
approach which, incorporating the Operator Regularized formuation,
CoV, change of unknown, and precision preserving techniques (items~(i)
through~(iv) above), achieves high-order convergence for any CoV
degree---as illustrated in the figure for the CoV degrees $p=4$ and
$p=6$. In particular, the $p=6$ CoV results in an observed convergence
rate faster than 9th order. The faster-than-expected convergence
results from the high spectral accuracy achieved in patches away from
corners. Notably, while these patches carry the dominant error in the
examples considered, they still exhibit faster convergence compared to
the error produced under the CoV around corners---as clearly
illustrated in Figure~\ref{fig:square_errors} left in the
numerical-results section. This observation suggests the potential for
designing an error-balancing strategy that mitigates these two sources
of error while reducing computational costs for a given error
tolerance. However, for the sake of simplicity, it may be reasonable
to accept a limited degree of suboptimality in exchange for a more
straightforward algorithmic design.

Next, we examine the algorithms that arise when some or all of the key
improvements are not included. The original method (including
point~(i) but excluding point~(ii), and therefore excluding
points~(iii) and~(iv) which are rendered irrelevant in absence of
point~(ii)) only achieves second order convergence, on account of the
corner singular character of the integral density. The Intermediate algorithm, which omits the change of unknown
strategy (point~(iii)), performs worse than the original approach,
emphasizing the significant impact inherent in this important
algorithmic strategy. The precision-preserving strategies mentioned in point~(iv) are important for any approaches using a CoV, since the CoV causes the discretization points to cluster so close to the corners that values of the source-target difference $|x-x'|$ exactly vanishes numerically for certain $x\ne x'$, thus giving rise to runtime errors at significant numbers of discretization points.
%The Intermediate~3 curve, which incorporates
%both the $p=4$ CoV and the change of unknown (item~(ii) with $p=4$ and
%item~(iii)), but which excludes the %precision preserving techniques
%(point~(iv)), performs as well as the full proposed method for the
%case with CoV order $p=4$.  Note that
%Figure~\ref{fig:square_allmethodcomparison} does not include an error
%curve excluding point~(iv) for the CoV degree $p=6$. The reason for
%this is that for the larger CoV degree the discretization points
%cluster so close to the corners that values of the source-target
%difference $|x-x'|$ exactly vanishes numerically for certain
%$x\ne x'$, thus giving rise to runtime errors at significant numbers
%of discretization points.  This underscores the importance of  the
%precision-preserving strategies mentioned in point~(iv) for higher CoV
%degrees.

\subsection{Conditioning and Number of GMRES Iterations}\label{conditioning}

\begin{figure}[!ht] % example dataset
  \centering% <-- superfluous in this example as commented by Zarko
  \includegraphics[width=0.4\textwidth] {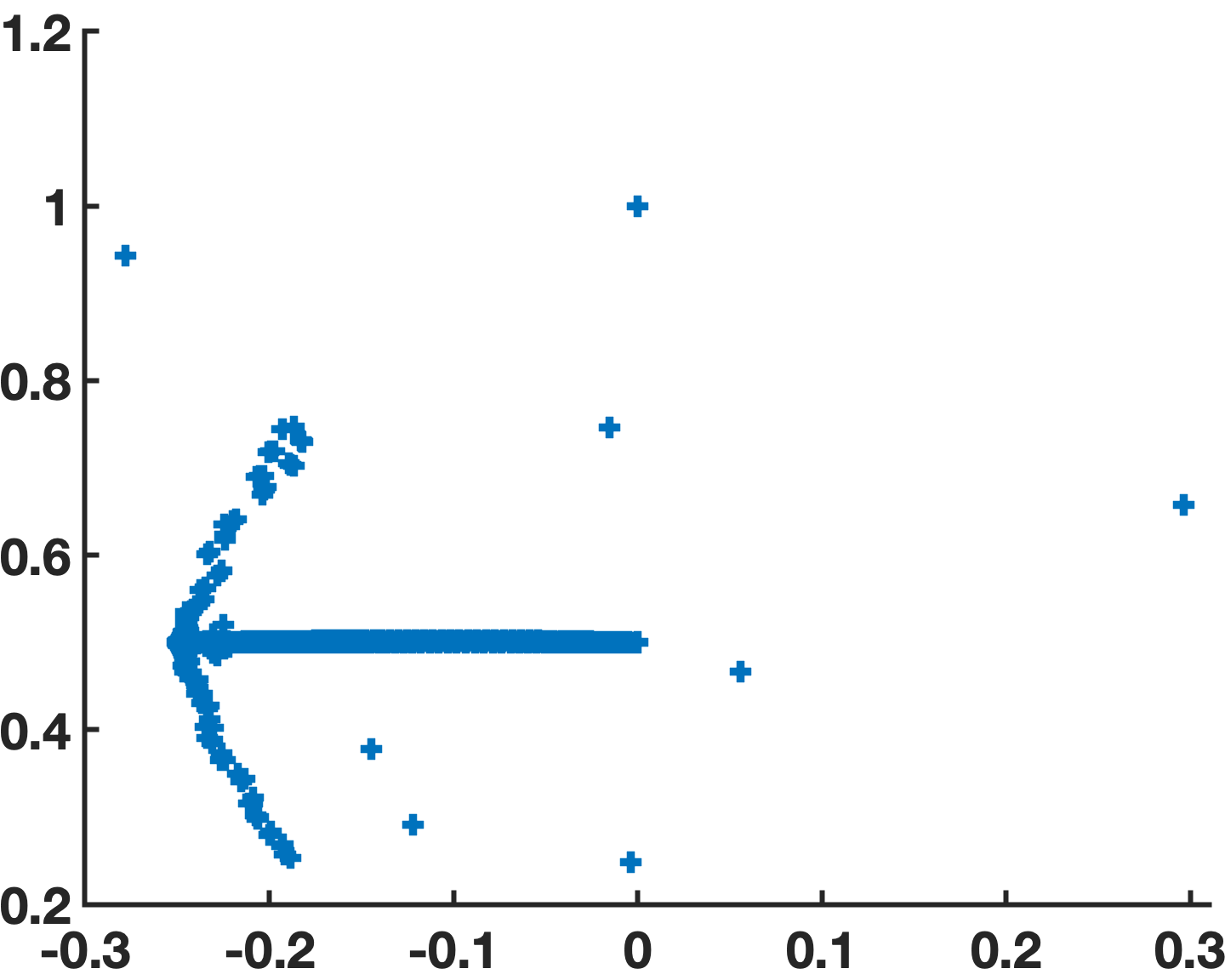}\hspace{0.5cm}
  \includegraphics[width=0.4\textwidth]
  {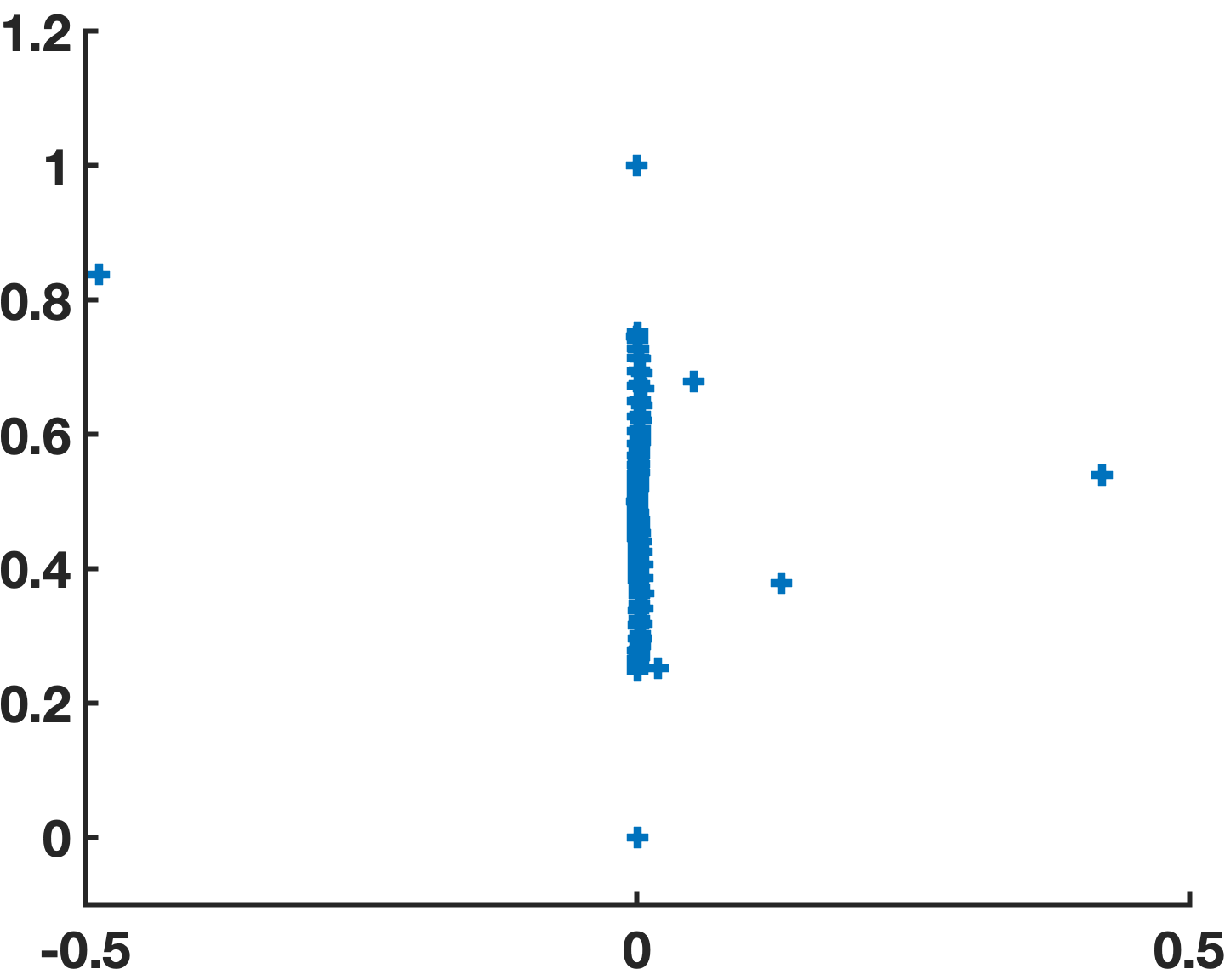}% <-- No space here
  \caption{Eigenvalue distributions for the left hand-side operators in the ``change-of-unknown'' regularized equations~\eqref{absorbed_ricfie}  and~\eqref{absorbed_mfie} (left and right panels, respectively), with $k =2.2214$, and  on a PEC cylinder with square
    cross-section.}
  \label{fig:eigenvalues_square_resonance}
\end{figure}

\begin{figure}[!ht] % example dataset
    \centering% <-- superfluous in this example as commented by Zarko
        \includegraphics[width=0.4\textwidth]
        {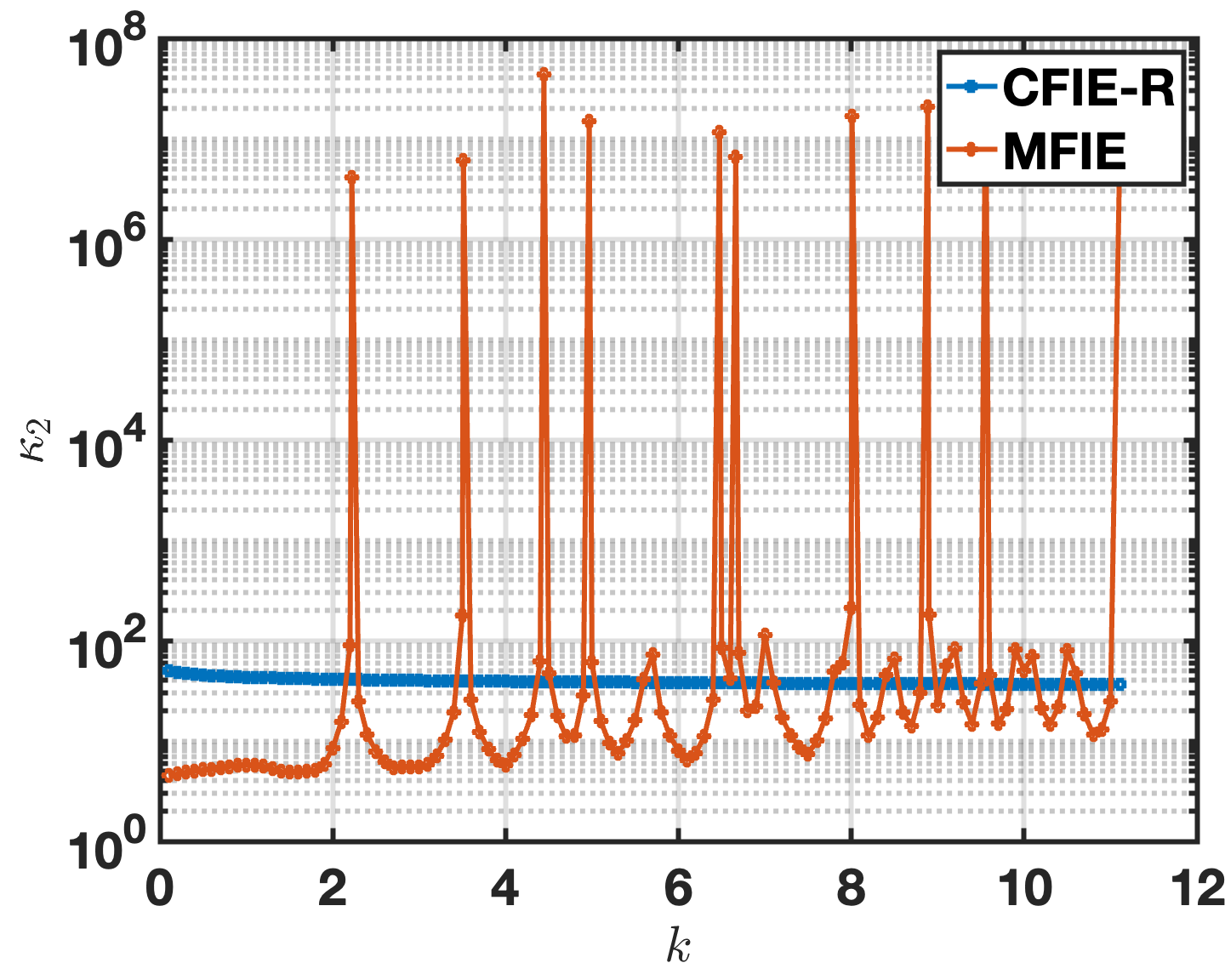}\hspace{0.5cm}
%   <-- No space here
        \includegraphics[width=0.4\textwidth]
        {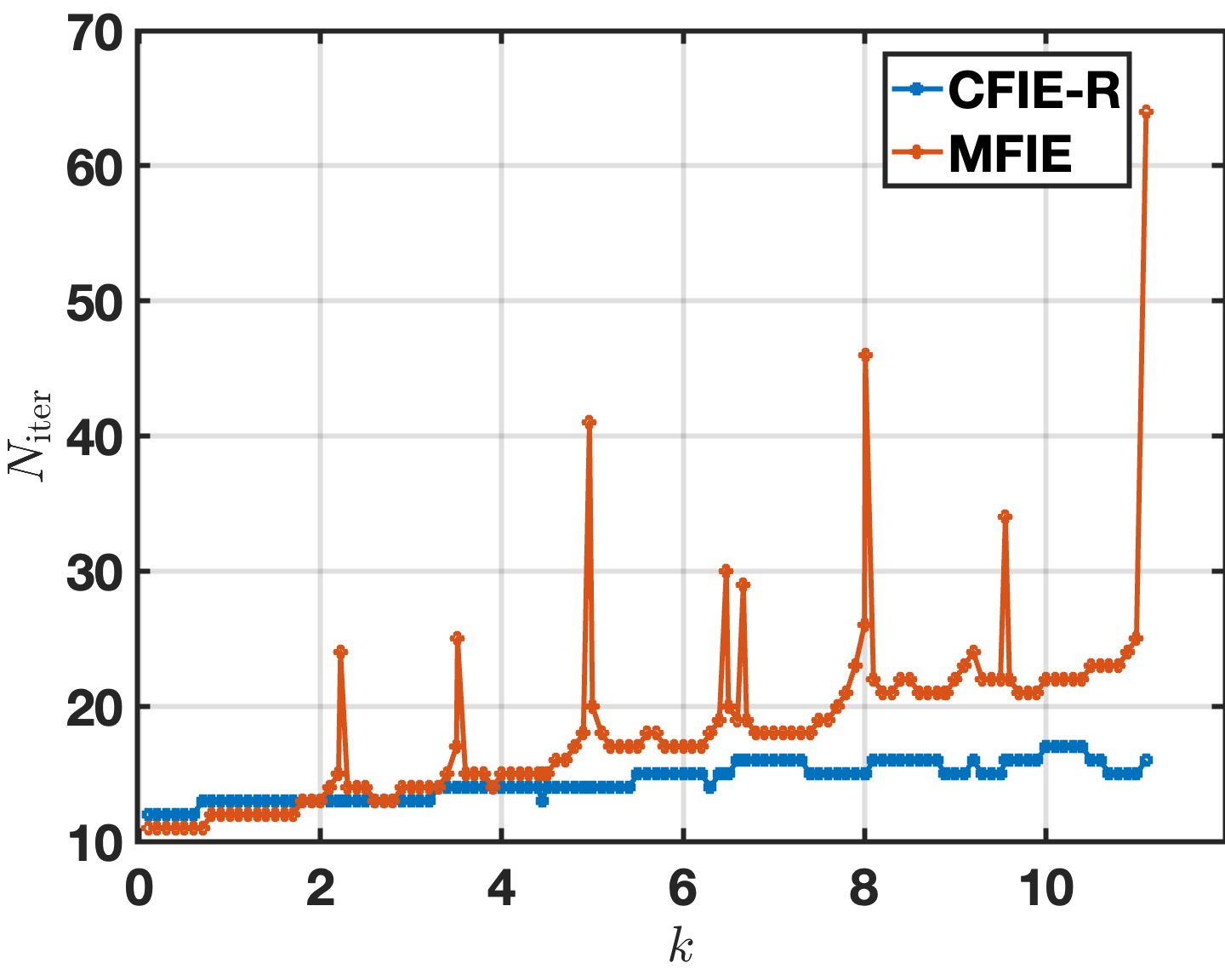}%      <-- No space here
\caption{Condition number $\kappa_2$ (left) and number of GMRES iterations required to reach $10^{-5}$ residual (right) as functions of the wavenumber $k$ for a PEC cylinder with square cross-section of side 2.}
\label{fig:square_condnumber_gmres}
\end{figure}

In order to highlight the attractive linear-algebra conditioning
properties resulting from the proposed corner-capable CFIE-R
regularized combined formulation, we compare it in this section with
the corresponding MFIE formulation---using, in both cases, all of the
key improvements~(i) through~(iv) in
Section~\eqref{comparisonagainstRP}. In order to do this in a
controlled manner, we consider once again the problem of scattering by
a square cylinder of side $a = 2$, for which the interior Laplace
resonances can be computed in closed form:
$$k_c = \frac{\pi}{a}\sqrt{ m^2 + n^2},
$$
for all positive integer values of $m$ and $n$. We thus consider the
BIE systems for $k = \frac{\pi\sqrt{2}}{2}\approx 2.2214$ and, to
visualize the singularity/non-singularity of the associated
MFIE/CFIE-R matrix equations, we present in Figure
\ref{fig:eigenvalues_square_resonance} (left for MFIE and right for CFIE-R) the corresponding eigenvalue
distributions. In particular, we see that the MFIE system is singular:
$0$ is an eigenvalue. At the same wavenumber, in contrast, the
eigenvalues of the CFIE-R system are all non-zero. Figure
\ref{fig:square_condnumber_gmres} left, which displays the condition
numbers vs. wavenumber for the MFIE and CFIE-R formulations, shows
significant condition-number spikes for the MFIE formulation at each
of the resonances corresponding to the cavity modes of the square
cylinder while remaining bounded for the CFIE-R formuilation---thus
indirectly demonstrating the resonance-free/singular character of the
CFIE-R/MFIE formulations considered.

Figure~\ref{fig:square_condnumber_gmres} right presents the number
of iterations required by the GMRES iterative algorithm to achieve a
residual of $10^{-5}$ vs. wavenumber for both the MFIE and CFIE-R. It
can be seen that the CFIE-R system achieves the prescribed tolerance
for all wavenumbers considered in less than 15 iterations, growing
only moderately with the wavenumber, while the number of iterations
required to solve the MFIE spikes near each resonance and grows
steadily even away from the resonances. This highlights the favorable
iteration numbers and resonance-free properties of the proposed CFIE-R
corner regularized formulation. We expect that the favorable number of
iterations of the CFIE-R formulation will be crucial in the 3D case
where iterative solvers become a necessity for large problems.

\subsection{Field evaluation accuracy}\label{accuracy}
\begin{figure}[!ht]
  \begin{center}
    \includegraphics[width=0.4\textwidth]{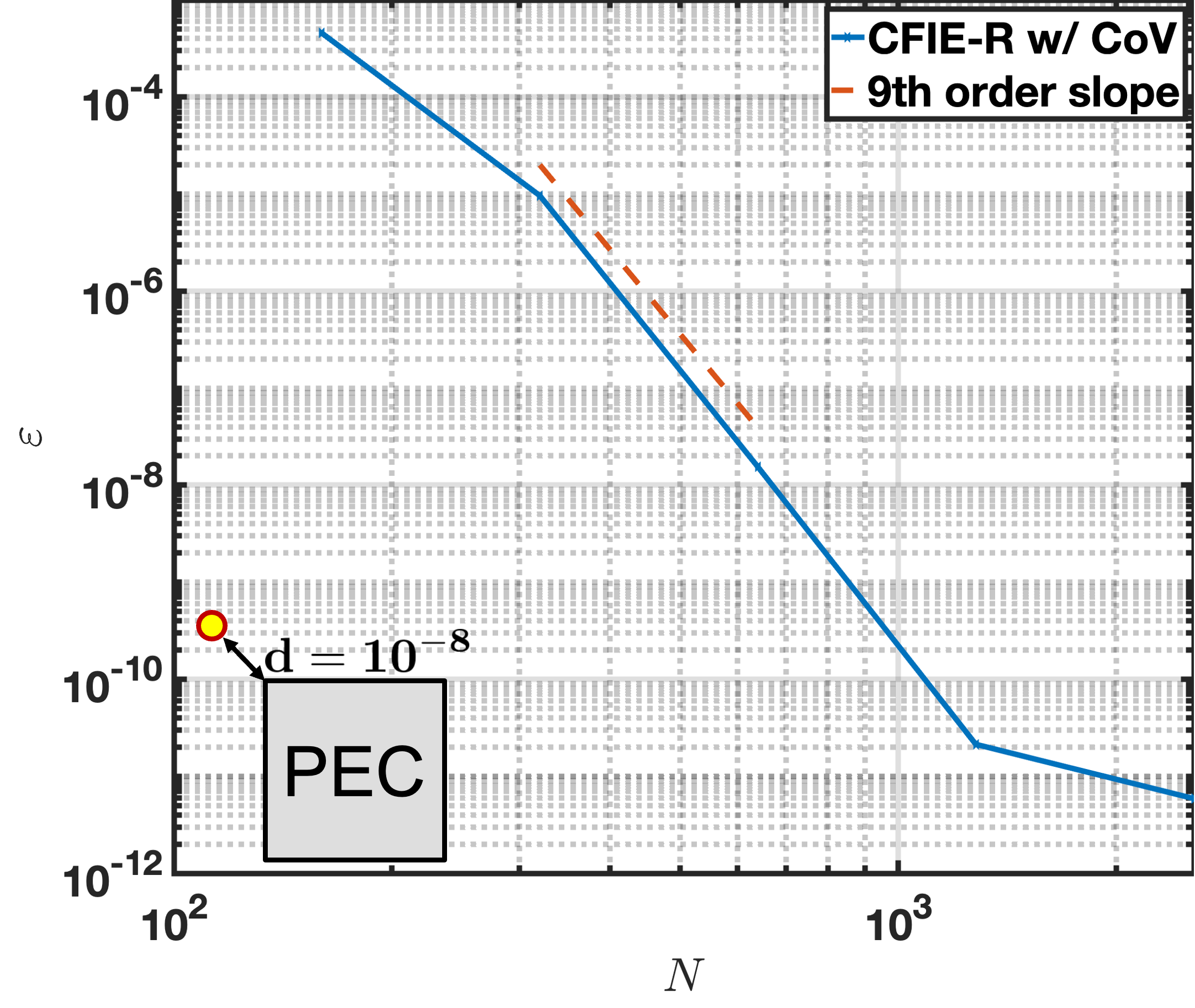}
  \end{center}
  \vspace{-0.5cm}
  \caption{Convergence of the field $u$  at a distance $10^{-8}$
    away from a corner of a PEC cylinder of square
    cross-section with point-source excitation.}
\label{fig:square_monopole}
\end{figure}

In this section, the accuracy of the proposed method with all of the
key improvements incorporated is evaluated for several different
geometries and excitations. The first example considered is a PEC
cylinder with square cross-section of sidelength $a=2$ with a monopole
point source excitation $u^{\mathrm{inc}}(r) = H_0^{(1)}(k_0 |r|)$
located within the scatterer. Due to the PEC nature of the object
(under the boundary conditions~\eqref{Neumann}), the resulting
scattered field outside the scatterer equals the negative of the
incident monopole field. Thus, the incident field can be used to
determine the accuracy of the scattered field. Figure
\ref{fig:square_monopole} plots the convergence of the
corner-regularized CFIE-R with CoV order $p=6$ vs. number of unknowns
of evaluating the field at a point that is a distance $10^{-8}$ away
diagonally from the top left corner of the square cylinder. Errors
smaller than  $10^{-11}$ are obtained with a 9th-order convergence
slope.

\begin{figure}[!ht] % example dataset
    \centering% <-- superfluous in this example as commented by Zarko
        \includegraphics[width=0.4\textwidth]
        {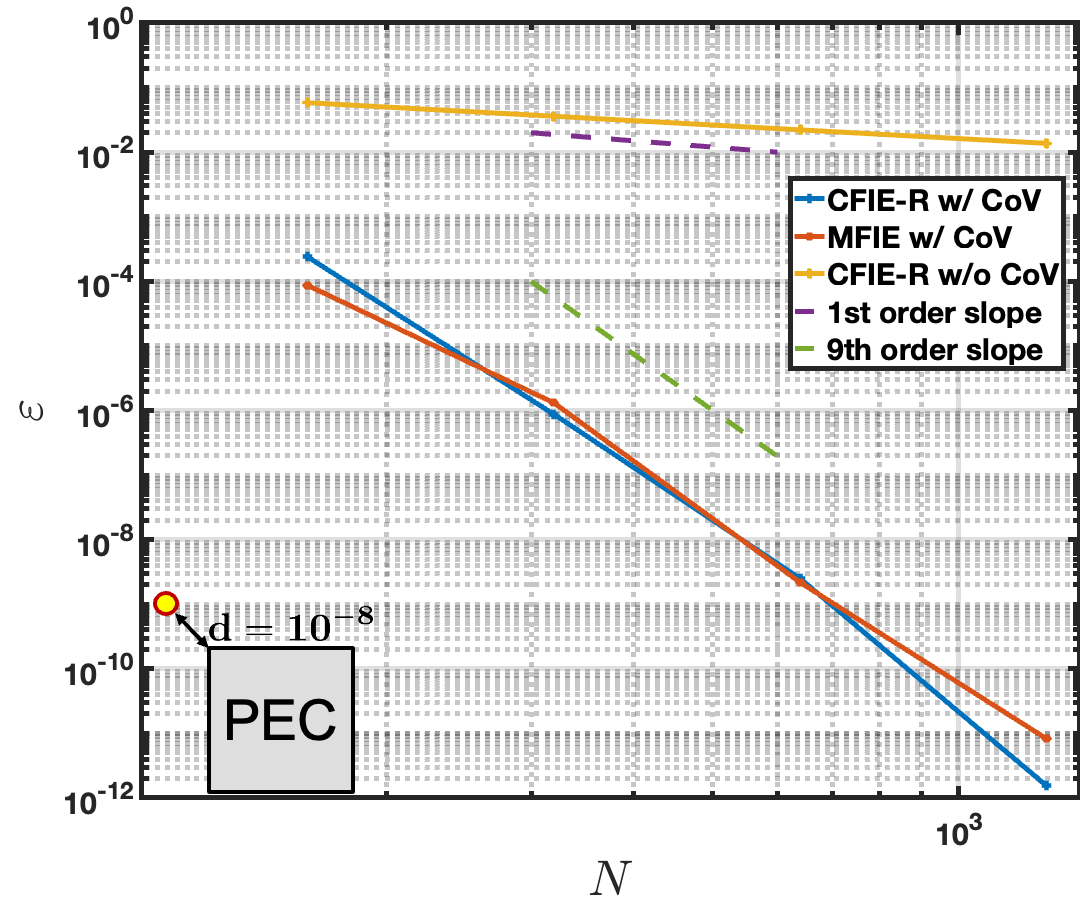}\hspace{0.5cm}
        \includegraphics[width=0.4\textwidth]
        {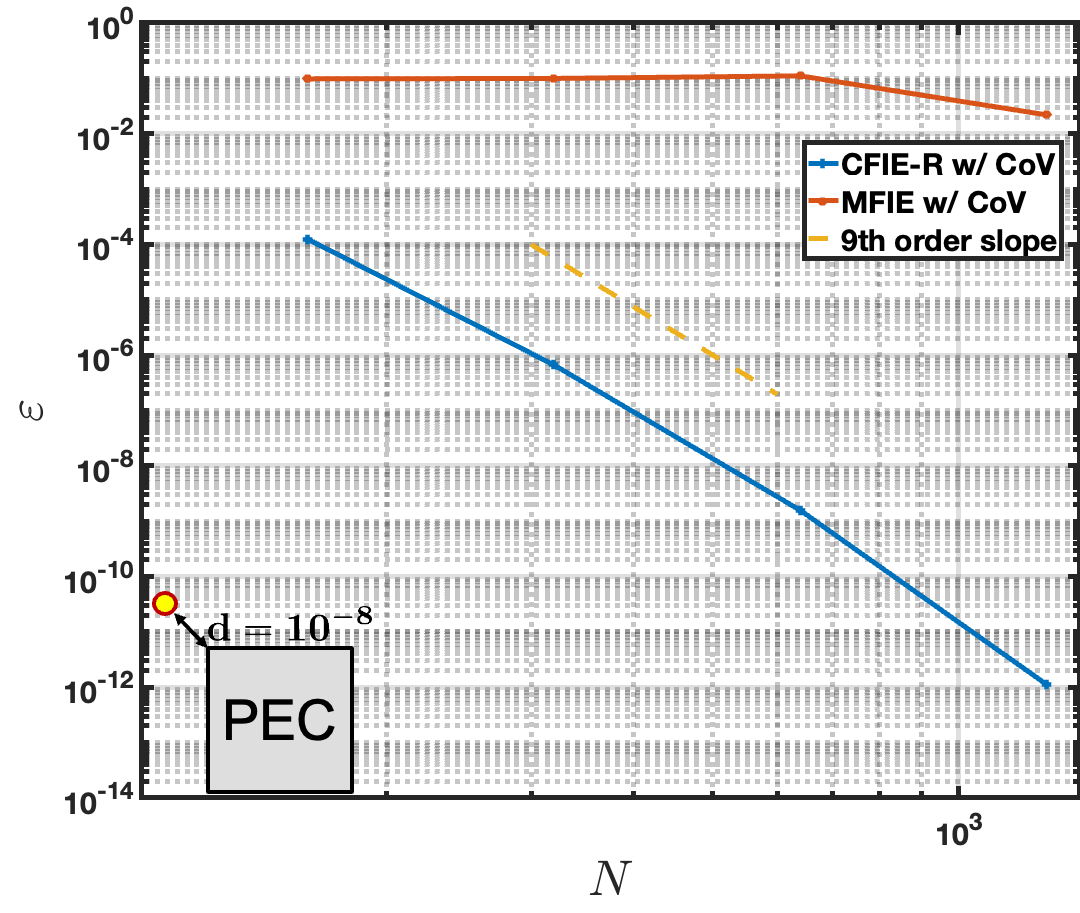}%      <-- No space here
        \caption{Convergence of the field $u$ for PEC cylinder of
          square cross-section for the non-resonant frequency $k = 10$
          (left) and the resonant frequency $k = 9.5548$ (right) at a
          distance $10^{-8}$ away from a corner.}
\label{fig:results}
\end{figure}

For the next example, the scattering from the same square
cross-section cylinder is considered with the planewave incident
excitation illustrated in Figure~\ref{fig:excitation}. This example
compares the performance of the corner-regularized MFIE and CFIE-R
formulations at wavenumbers $k=10$ and $k=9.5548$, the second of which
is very close to a resonance of the MFIE. A CoV of order $p=6$ is used
in all cases except for a curve ``w/o'' CoV. Figure ~\ref{fig:results}
left illustrates the convergence behavior of the field produced by the
corner-regularized MFIE and CFIE-R formulations at a point
$d = 10^{-8}$ diagonally offset from the top-left corner of the
scatterer, with a wavenumber $k = 10$. Both the MFIE and CFIE-R
produce similar accuracies when the wavenumber is not near an MFIE
resonance, achieving a relative error of $10^{-12}$ and 9-th order
convergence. For comparison, the convergence of the CFIE-R without a
change of variables (CoV) is also plotted, demonstrating the expected
first-order convergence caused by the corner singularity, with errors
that remain larger than $10^{-2}$ even for discretizations containing
over 1000 unknowns. Figure ~\ref{fig:results} right plots the
convergence of the field evaluation at the same point for the
wavenumber is $k=9.5548$, which is very close to a resonance of the
MFIE operator. In this scenario, the CFIE-R exhibits convergence
nearly identical to that observed for $k=10$. However, the MFIE fails
entirely, producing inaccurate results regardless of the number of
unknowns used.

\begin{figure}[!ht] % example dataset
    \centering% <-- superfluous in this example as commented by Zarko
        \includegraphics[width=0.4\textwidth]
        {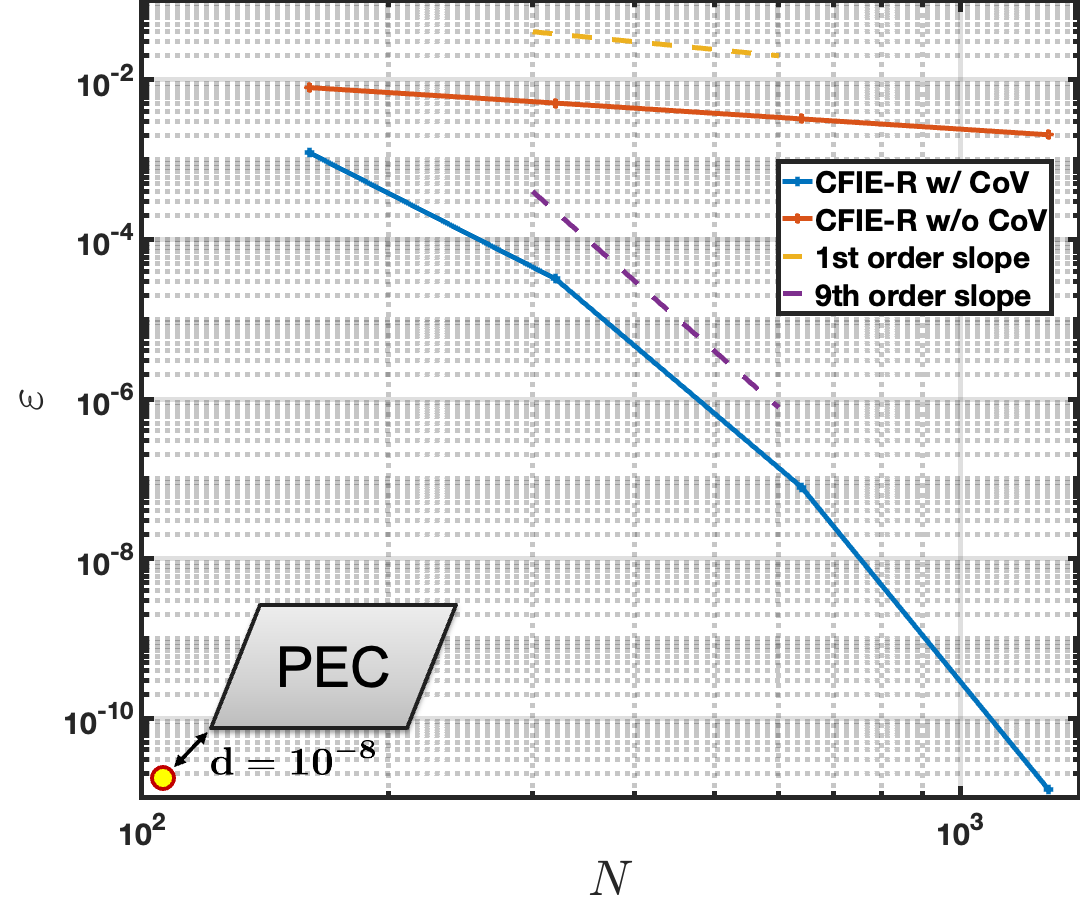}
    \hspace{0.5cm}
        \includegraphics[width=0.4\textwidth]
        {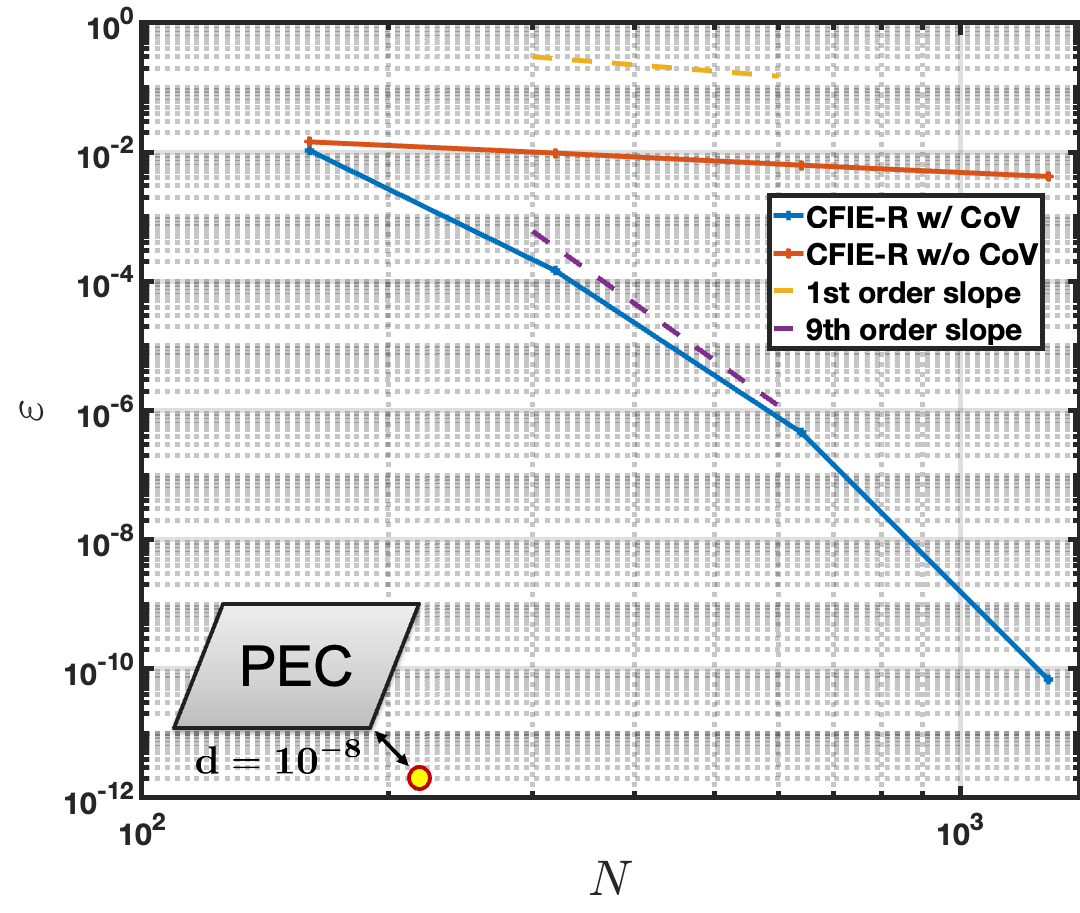}%      <-- No space here
      \caption{Convergence of the field $u$ for PEC cylinder with a
        parellogram-shaped cross-section at a point $10^{-8}$ away
        from a corner. Left: Acute angle. Right: Obtuse angle.}
\label{fig:results_parallogram}
\end{figure}

\begin{figure}[!ht] % example dataset
    \centering% <-- superfluous in this example as commented by Zarko
        \includegraphics[width=0.4\textwidth]
        {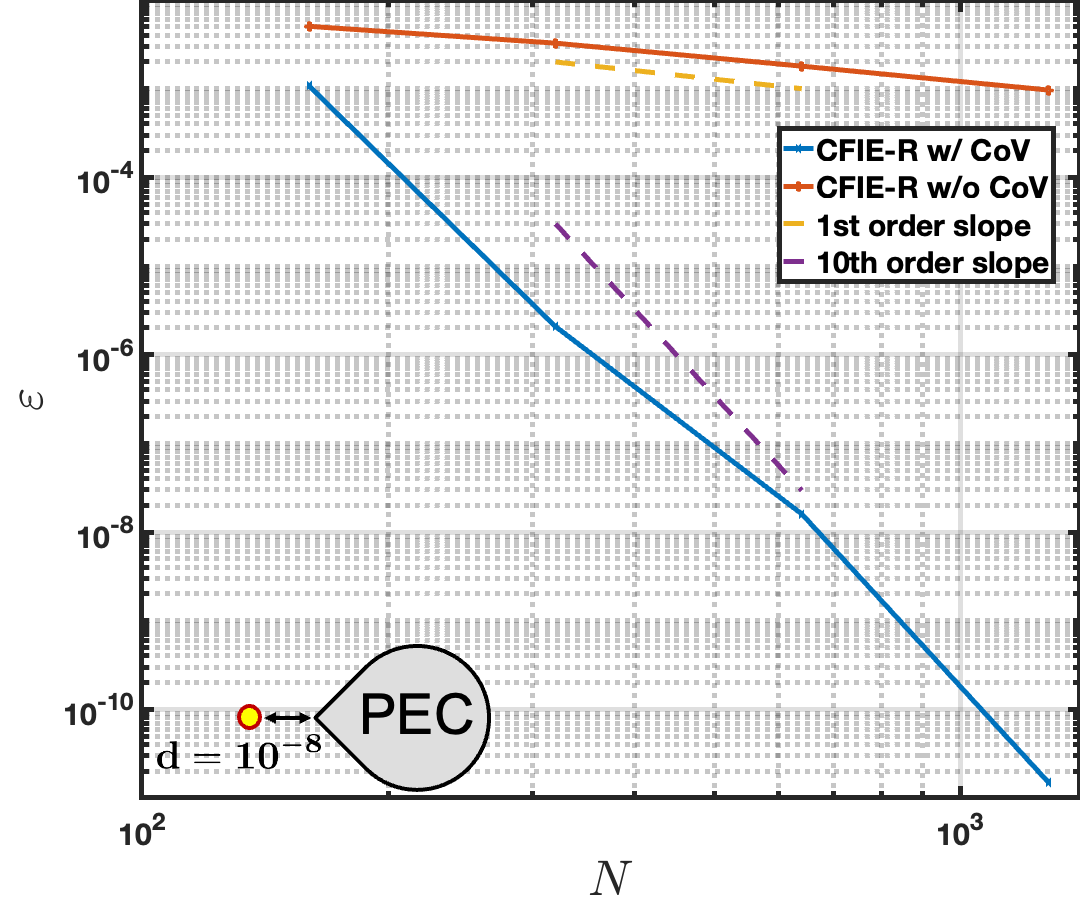}
    \hspace{0.5cm}
        \includegraphics[width=0.4\textwidth]
        {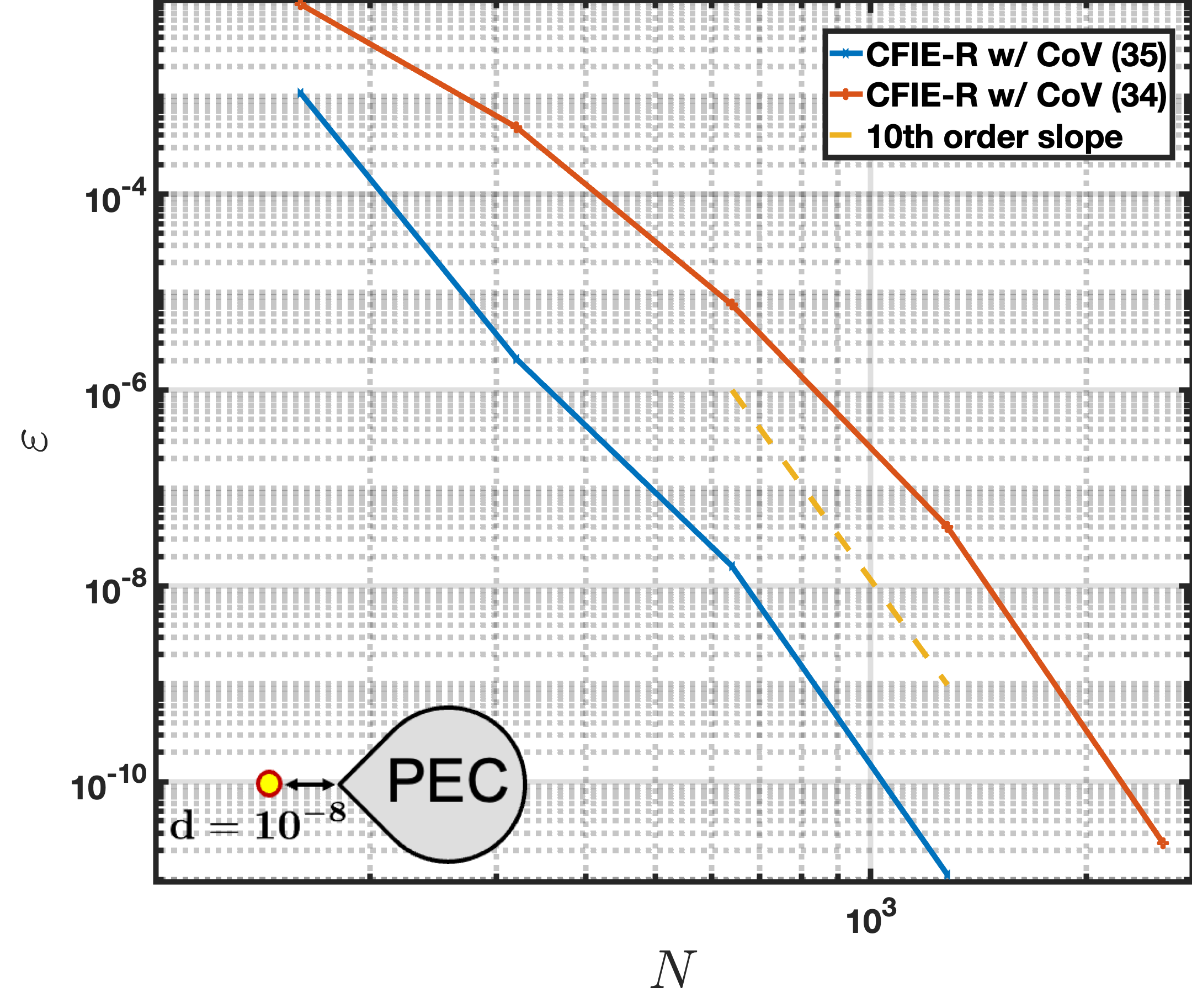}%      <-- No space here
      \caption{Left: Convergence of the field $u$ for PEC cylinder
        with a teardrop-shaped cross-section at a point $10^{-8}$ away
        from the corner. Right: comparison of the algorithm's
        performance for the two CoVs considered.}
\label{fig:results_teardrop}
\end{figure}

\begin{figure}[!ht] % example dataset
    \centering% <-- superfluous in this example as commented by Zarko
        \includegraphics[width=0.4\textwidth]
        {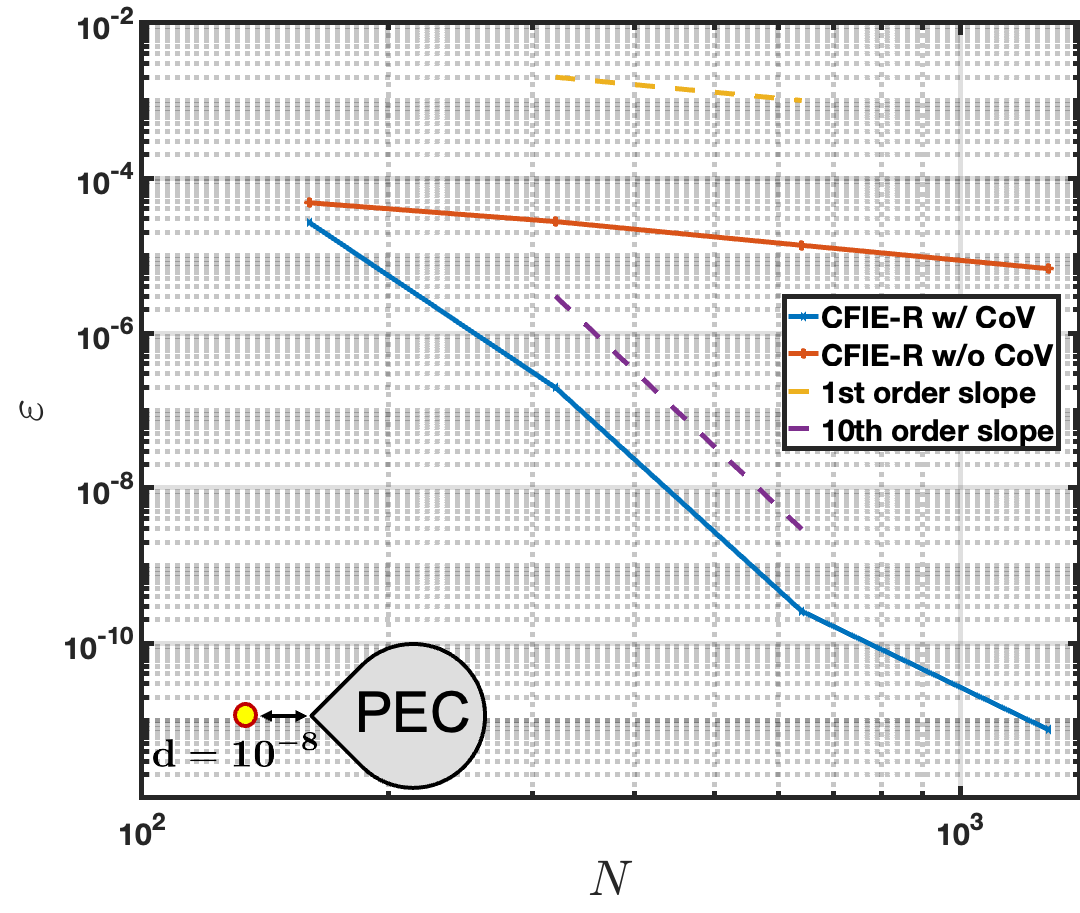}
    \hspace{0.5cm}
        \includegraphics[width=0.4\textwidth]
        {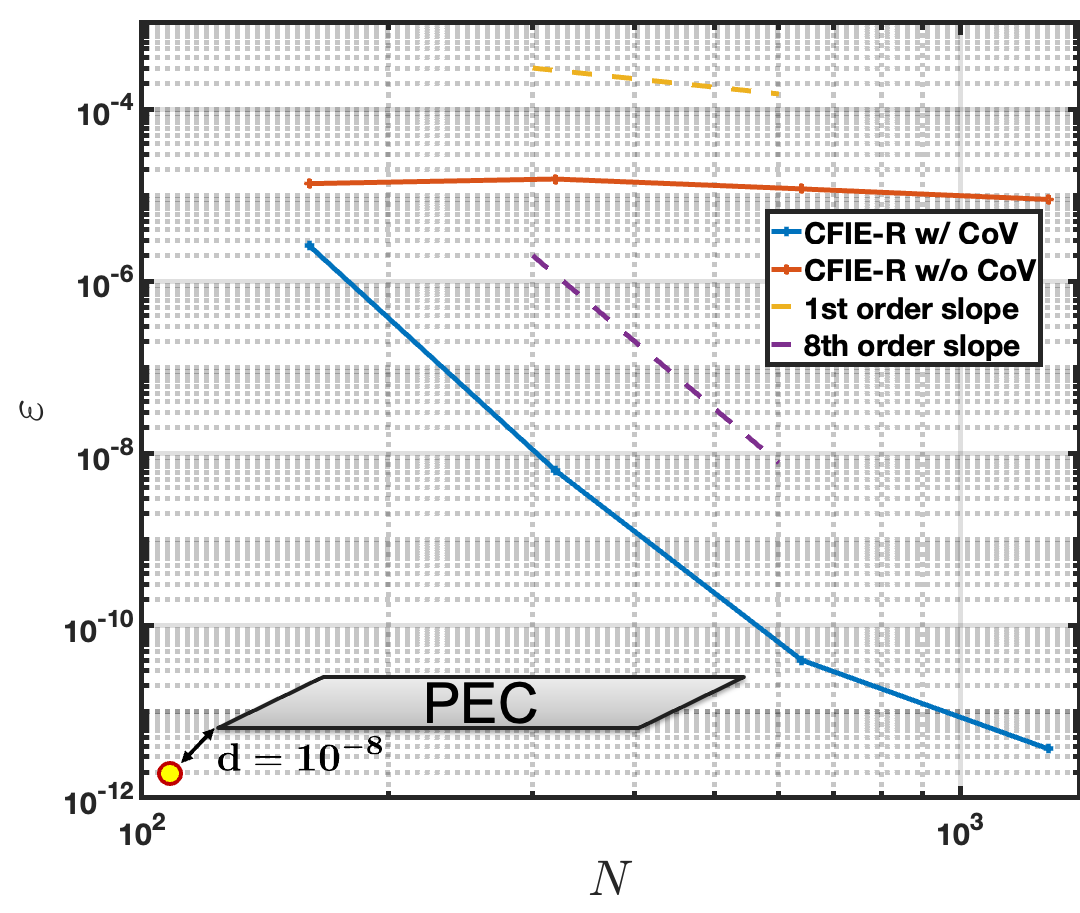}%      <-- No space here
        \caption{Convergence of the field $u$ for PEC cylinder with
          teardrop-shaped (left) and parellogram-shaped ``needle''
          (right) at points $10^{-8}$ away from a corner. (For the right image the observation point lies at a distance $d=10^{-8}$ of the needle corner point with acute angle $0.57^\circ(=0.01$ radians).)}
\label{fig:extreme_results}
\end{figure}
\begin{figure}[!ht] % example dataset
 \begin{center}% <-- superfluous in this example as commented by Zarko
        \includegraphics[width=0.3\textwidth]
        {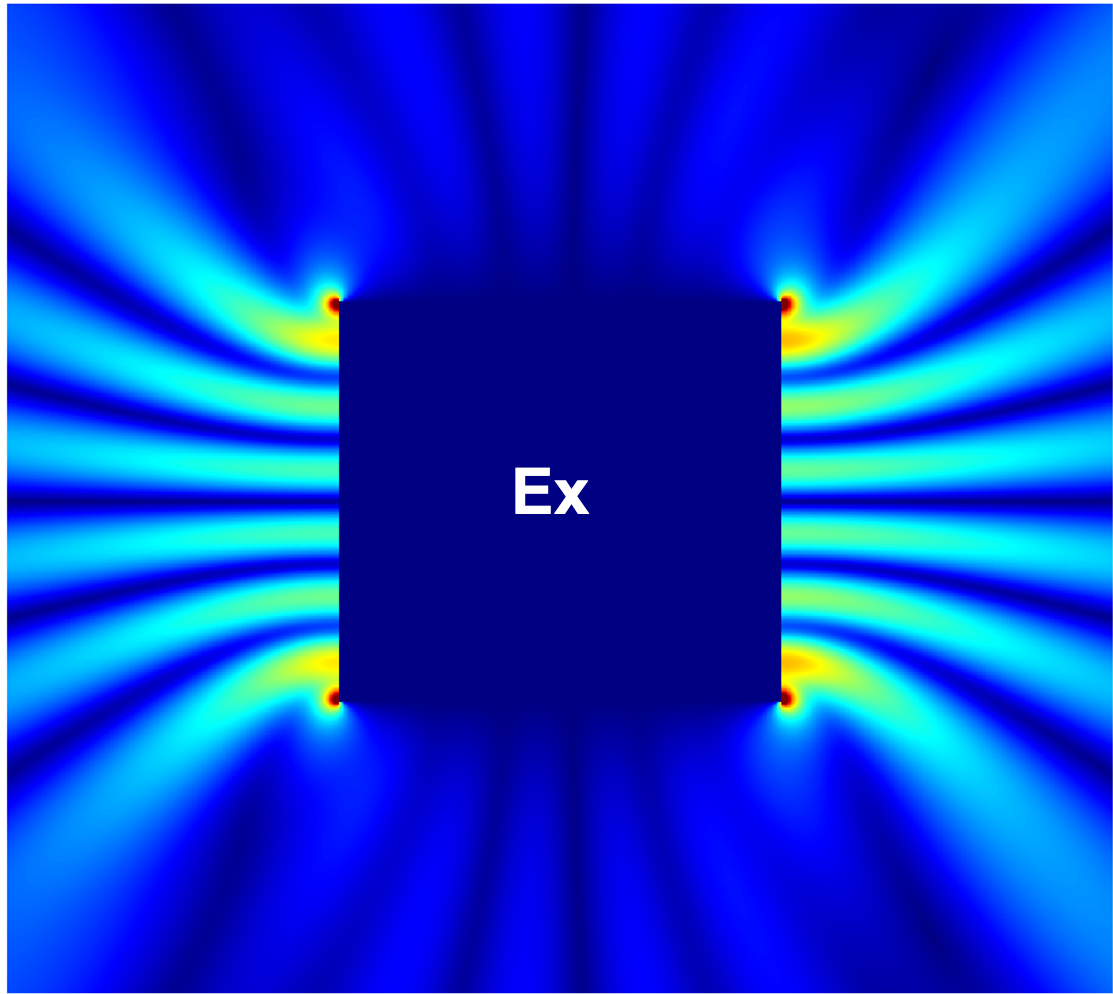}\hspace{0.05cm}
        \includegraphics[width=0.3\textwidth]
        {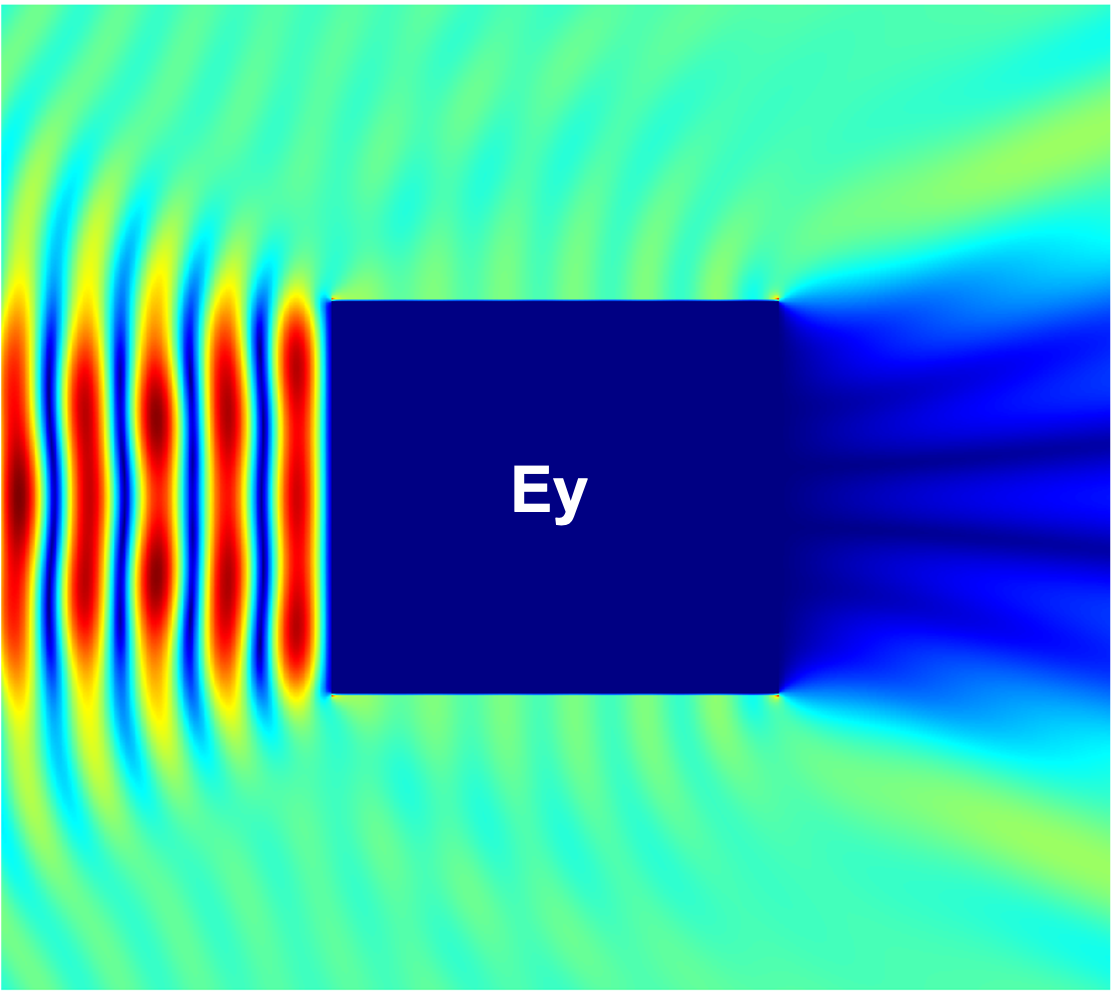}\hspace{0.05cm}
         \includegraphics[width=0.3\textwidth]
         {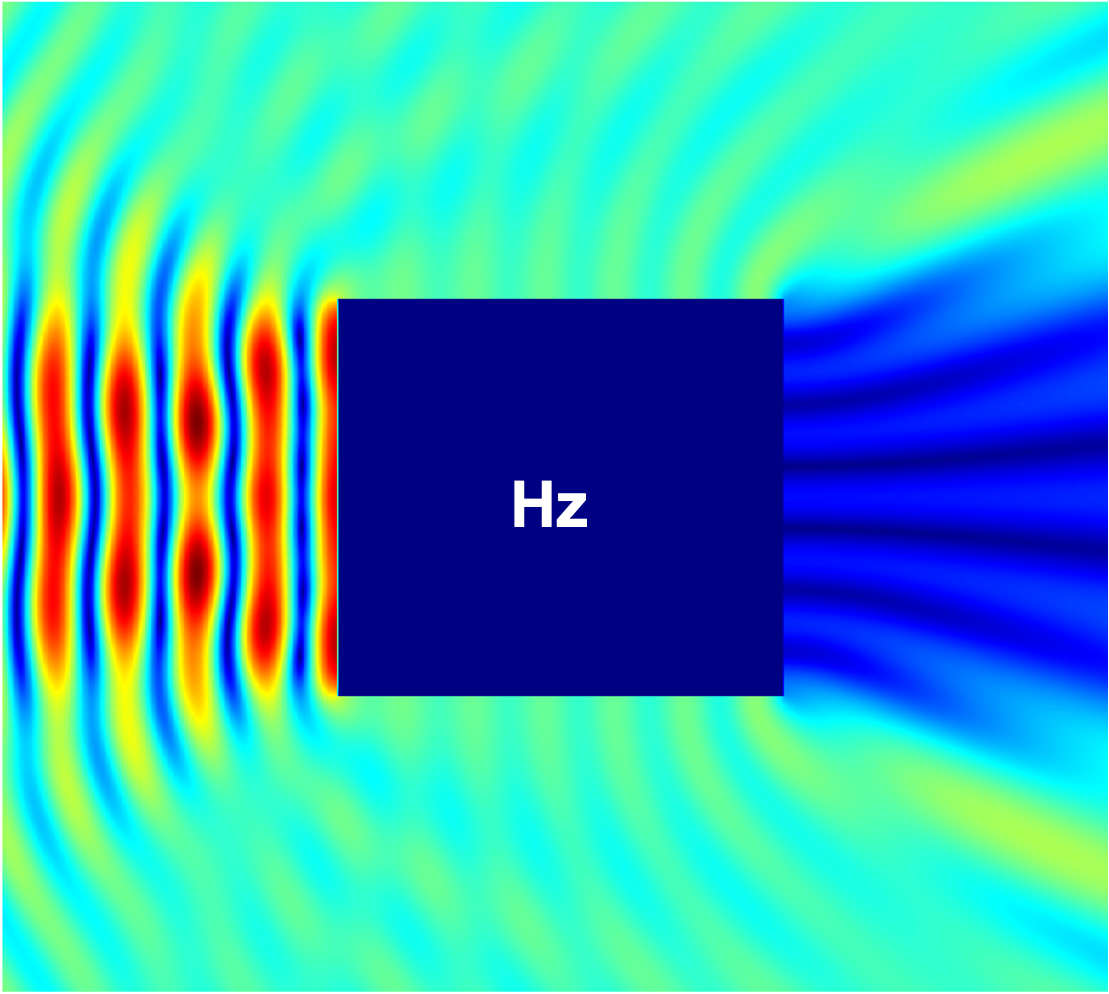}
        \includegraphics[width=0.3\textwidth]
        {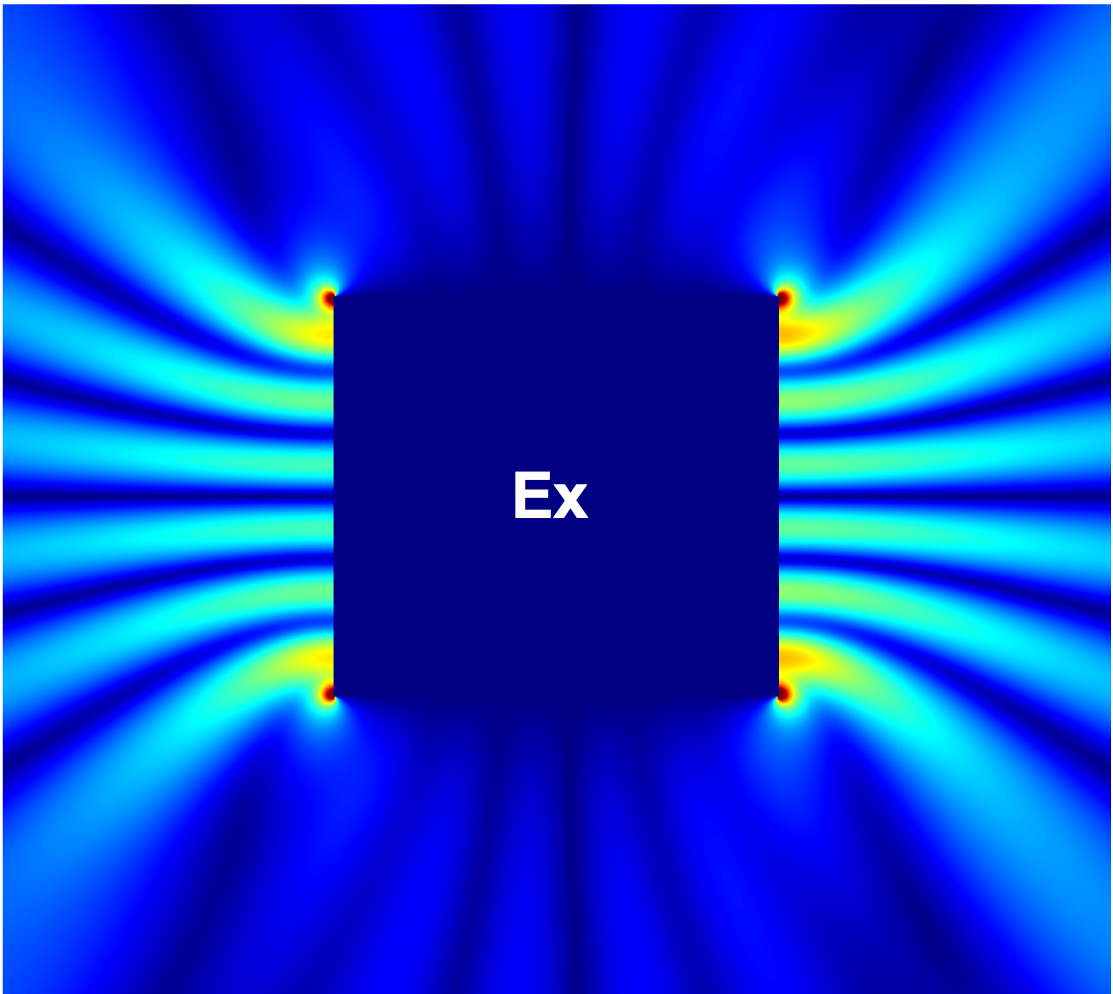}\hspace{0.05cm}
        \includegraphics[width=0.3\textwidth]
        {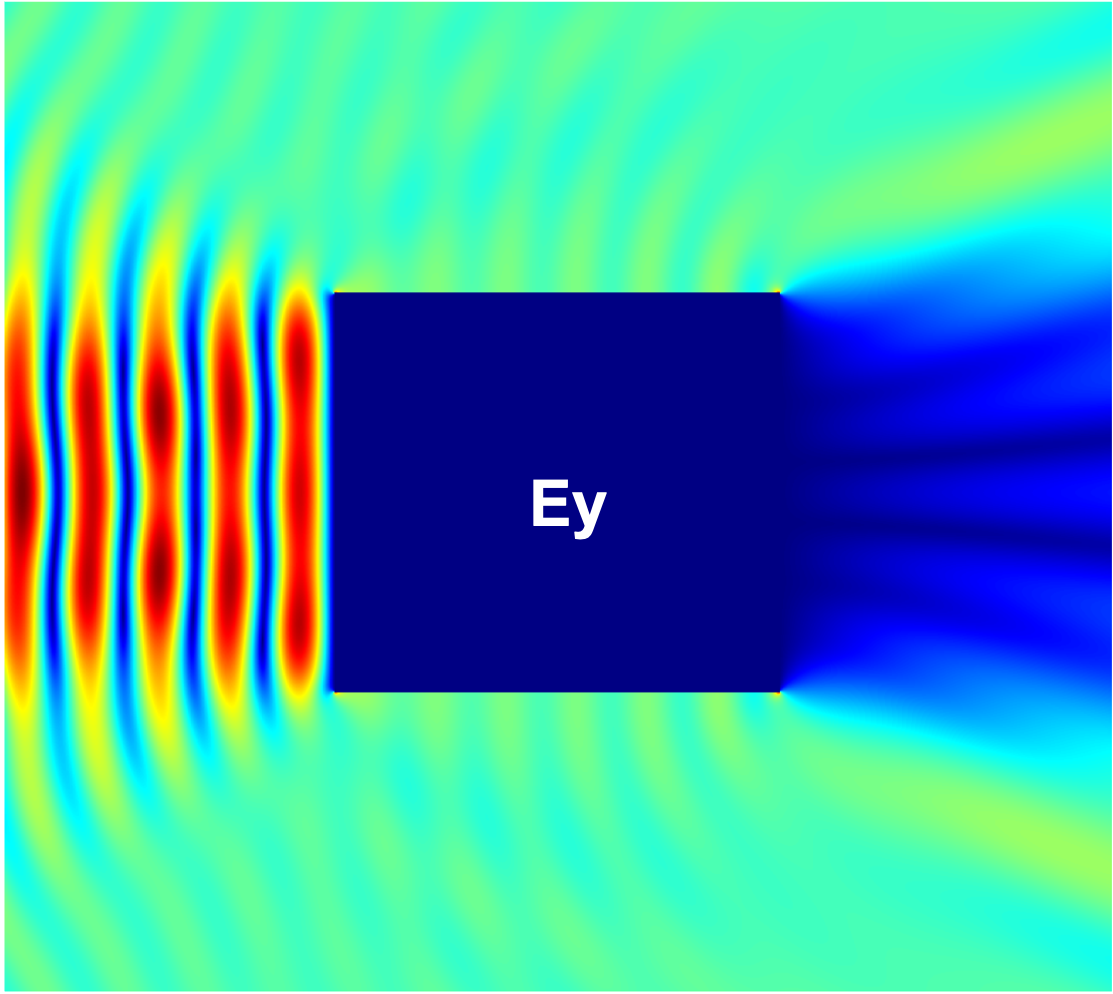}\hspace{0.05cm}
        \includegraphics[width=0.3\textwidth]
        {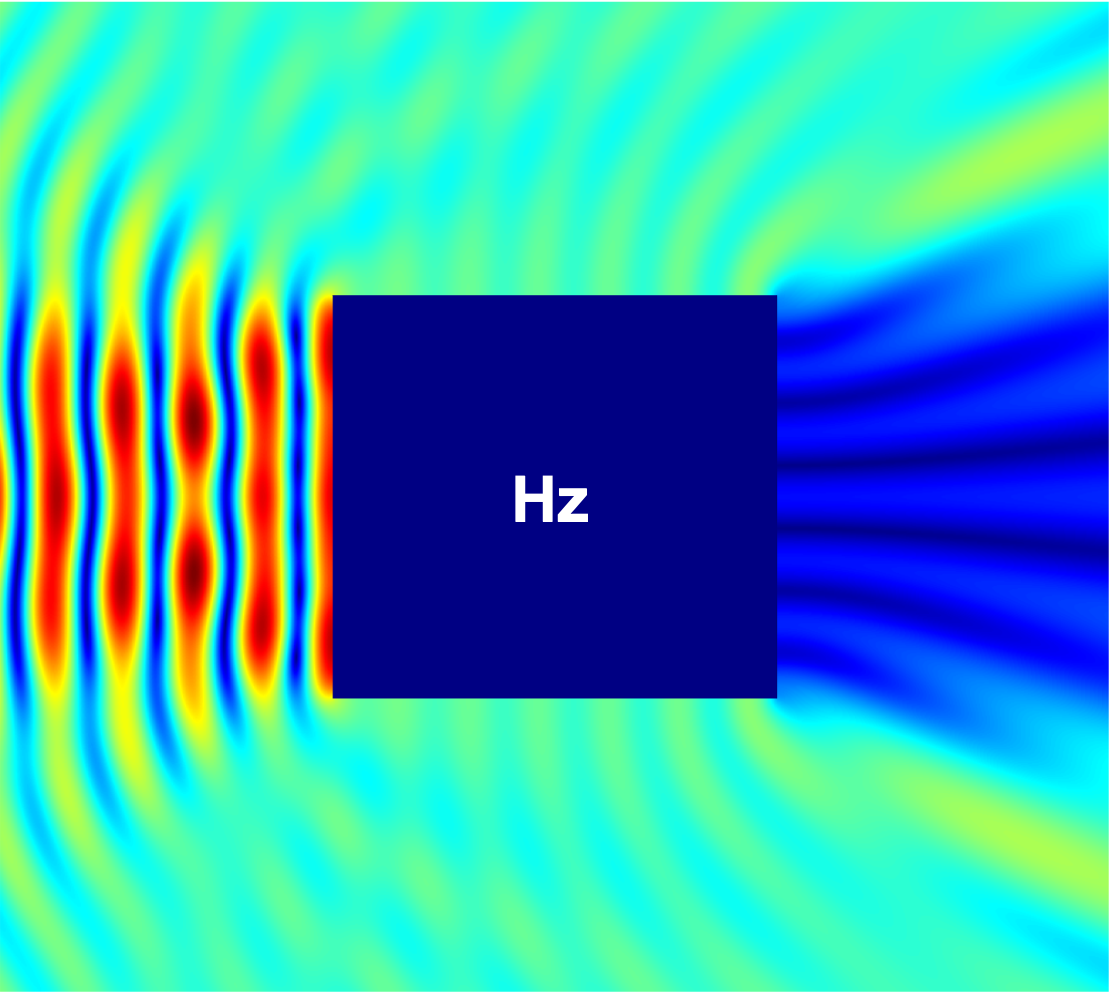}%      <-- No space here
        \end{center}
        \caption{Total field magnitude plots for $E_x,E_y$ and $H_z$
          field components for a cylinder with square
          cross-section. Upper row: Fine-mesh results produced by the
          commercial solver Lumerical. Lower row: Results produced by
          the proposed solver.}
\label{fig:lumericalvsbie}
\end{figure}

The next two examples examine the field convergence of the proposed
method for scatterers different from the square cross-section
cylinders considered up to this point, namely, cylinders with
parallelogram ($45^\circ$ acute angle) and teardrop-shaped
cross-sections. (The parallelogram has sides of lengths $2$ and
$2\sqrt{2}$.) The field convergence is evaluated at a point
$d=10^{-8}$ diagonally offset from the bottom-left and bottom-right
corners in Figure~\ref{fig:results_parallogram} left and right, which
have interior angles of $45^\circ$ and $135^\circ$, respectively. The
results indicate that the corner-regularized CFIE-R with a CoV order
of $p=6$ achieves similar convergence near both corners, despite their
differing interior angles, attaining better than 9th-order convergence
and a relative error below $10^{-11}$. This performance resembles
closely those observed for the square cross-section examples. In
contrast, the original rectangular polar method without CoV provides
first-order convergence only, and it fails to achieve even two digits
of accuracy.

The ``teardrop'' shaped geometry
$$x(z) = \left(2\sin{\frac{z}{2}},-\beta \sin{z} \right), \quad \beta = \tan{\frac{\alpha\pi}{2}, \quad 0 < \alpha < 1}, \quad 0\le z \le2\pi $$
is considered next where the parameter $\alpha$ controls the angle of
the corner; $\alpha$ values resulting interior angles of $90^{\circ}$
and $9^\circ$ are considered in what
follows. Figure~\ref{fig:results_teardrop} left illustrates the
convergence of the field evaluated at a point located a distance
$d=10^{-8}$ diagonally from the teardrop corner, showing that the CoV
order of $p=6$ yields 10th-order convergence with a relative error
below $10^{-10}$, demonstrating once again higher than expected
convergence on account of over-resolution of the corner. In contrast,
the original RP method~\eqref{absorbed_ricfie} exhibits only linear
convergence, achieving fewer than three digits of accuracy even with
more than 1000 unknowns.  Figure~\ref{fig:results_teardrop} right
compares, for the teardrop scatterer, the performance resulting from
use CoVs~\eqref{eq:CoVCK} and~\eqref{eq:power_cov}. While both CoV
functions achieve similar best-case accuracy and convergence rates in
field evaluation, the CoV~\eqref{eq:CoVCK} proves significantly more
efficient, delivering a solution that is up to three digits more
accurate for the same number of unknowns. This superior performance is
attributed to the uneven distribution of discretization points it
produces along the domain boundary, as discussed in
Section~\ref{corn_reg}.

The proposed methods consistently maintain accuracy across a wide
range of corner angles, as illustrated in what follows for a $9^\circ$
teardrop and a ``needle-like'' parallelogram with an acute angle of
approximately $0.57^\circ(=0.01$ radians). In both cases, the relative
error in the field evaluation is plotted for a point at a distance
$d=10^{-8}$ diagonally offset from the corner, as shown in
Figure~\ref{fig:extreme_results}. The results demonstrate that the
corner-regularized CFIE-R with a CoV order of $p=6$ achieves a
relative error of approximately $10^{-11}$ in both cases, with
10th-order convergence for the small-angle teardrop and 8th-order
convergence for the needle-like parallelogram. In contrast, the
original RP method for solving the CFIE-R system fails to achieve more
than three digits of accuracy and exhibits only linear convergence in
both cases.

\begin{figure}[!ht] % example dataset
    \centering% <-- superfluous in this example as commented by Zarko
        \includegraphics[width=0.3\textwidth]
        {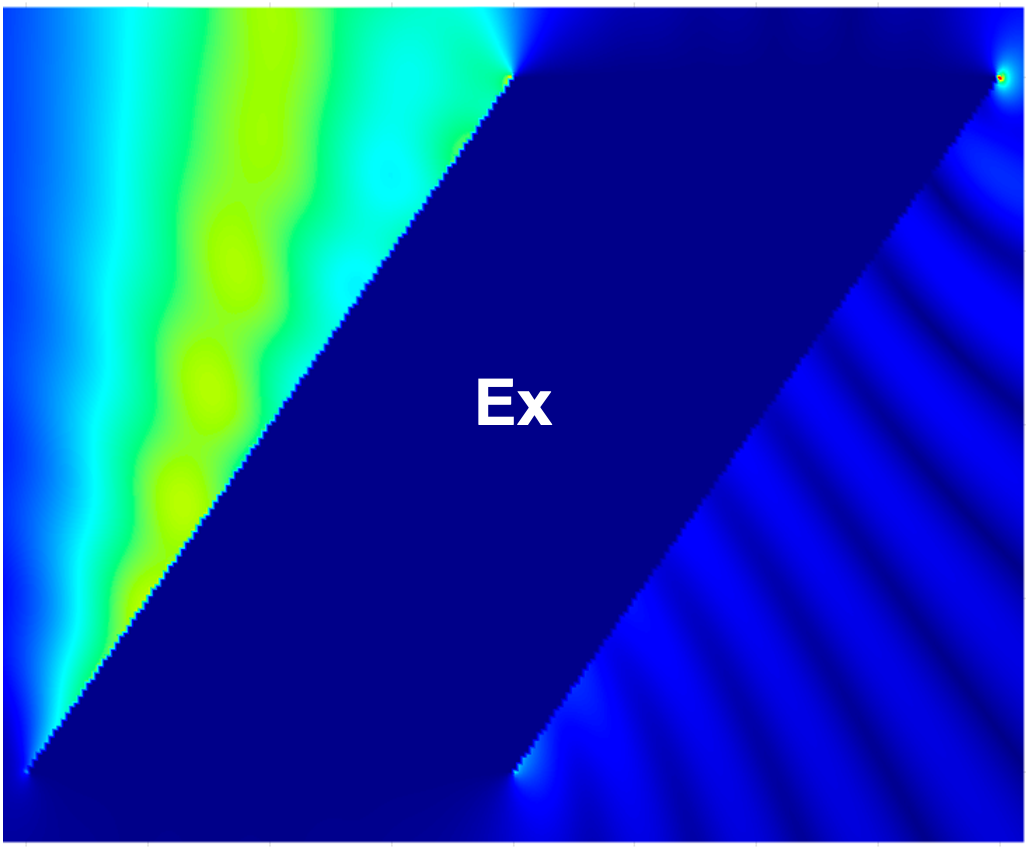}
   \hspace{0.05cm}
        \includegraphics[width=0.3\textwidth]
        {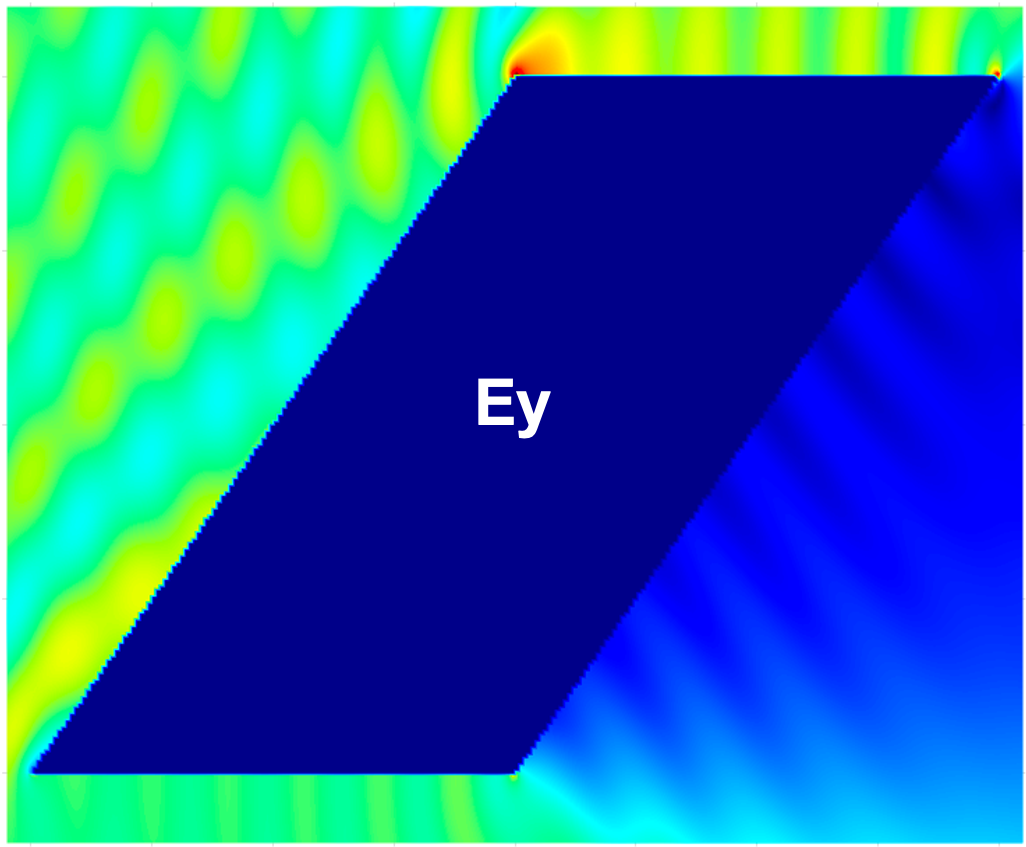}
    \hspace{0.05cm}
        \includegraphics[width=0.3\textwidth]
        {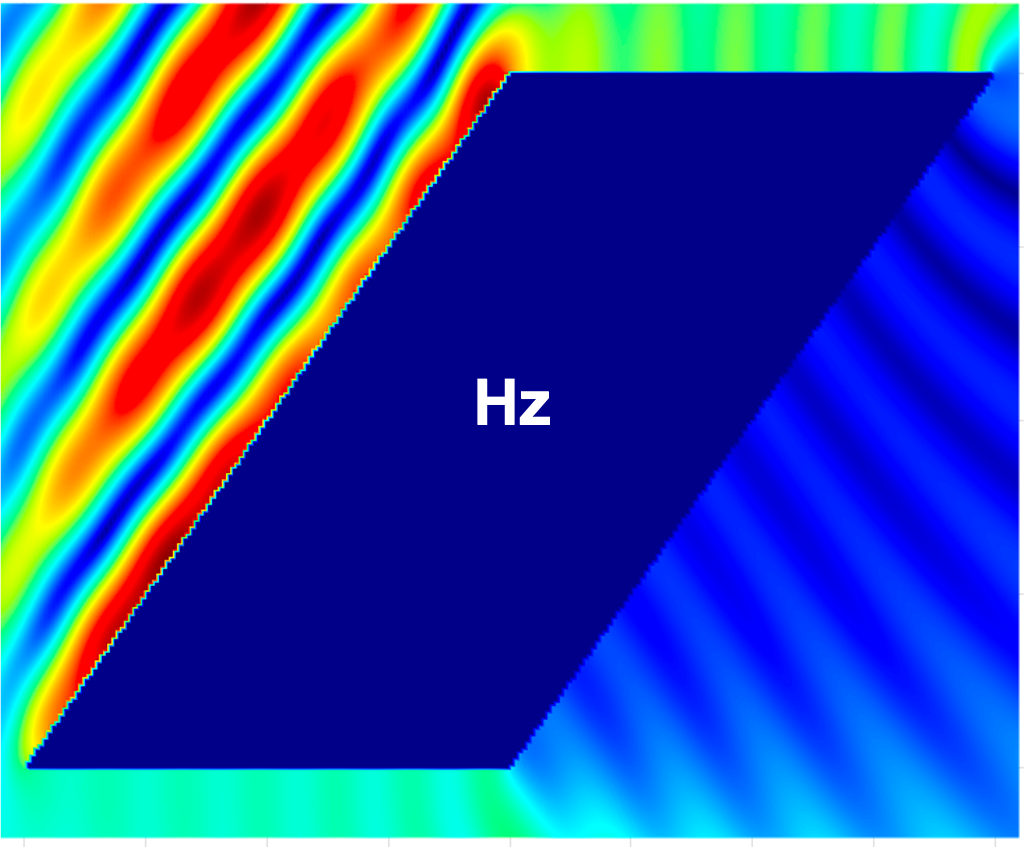}
        \includegraphics[width=0.3\textwidth]
        {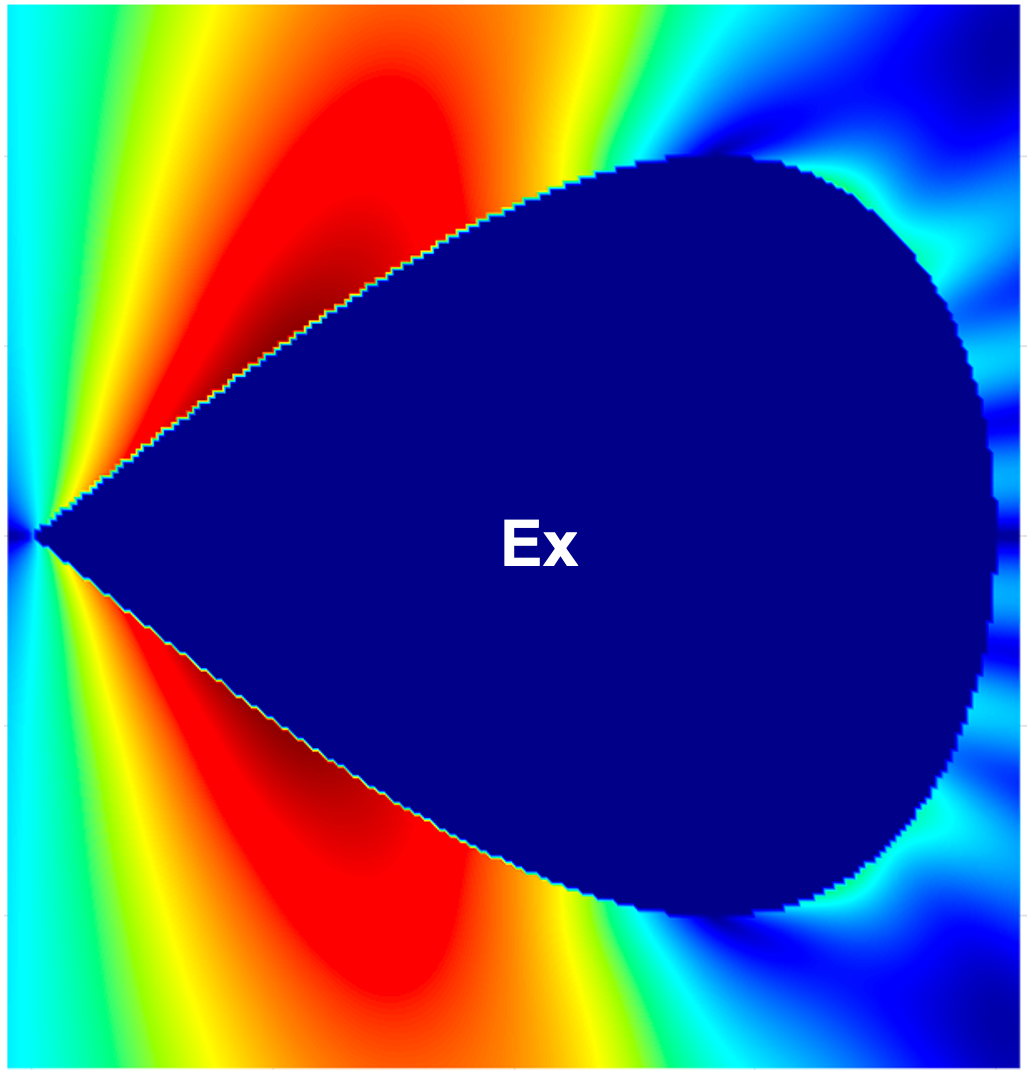}
    \hspace{0.05cm}
        \includegraphics[width=0.3\textwidth]
        {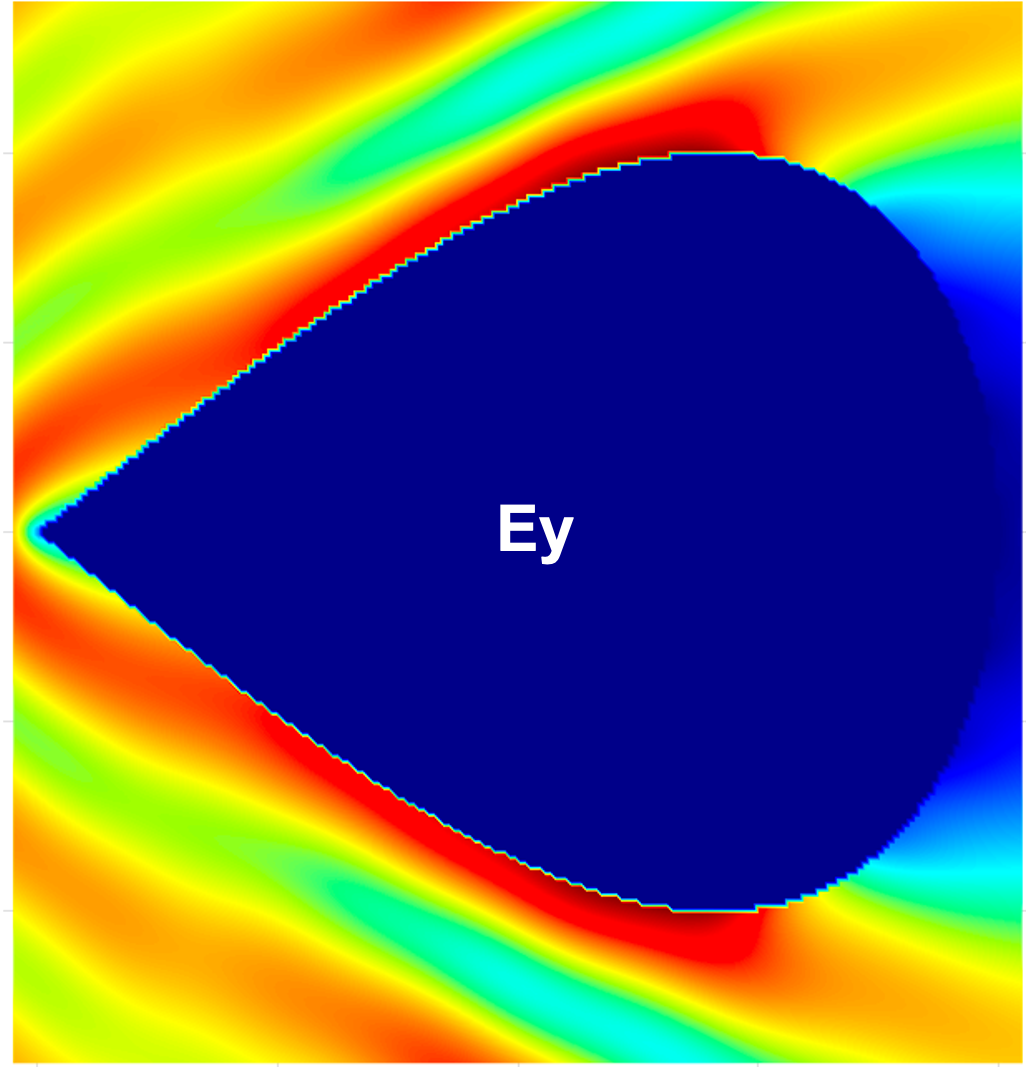}
    \hspace{0.05cm}
        \includegraphics[width=0.3\textwidth]
        {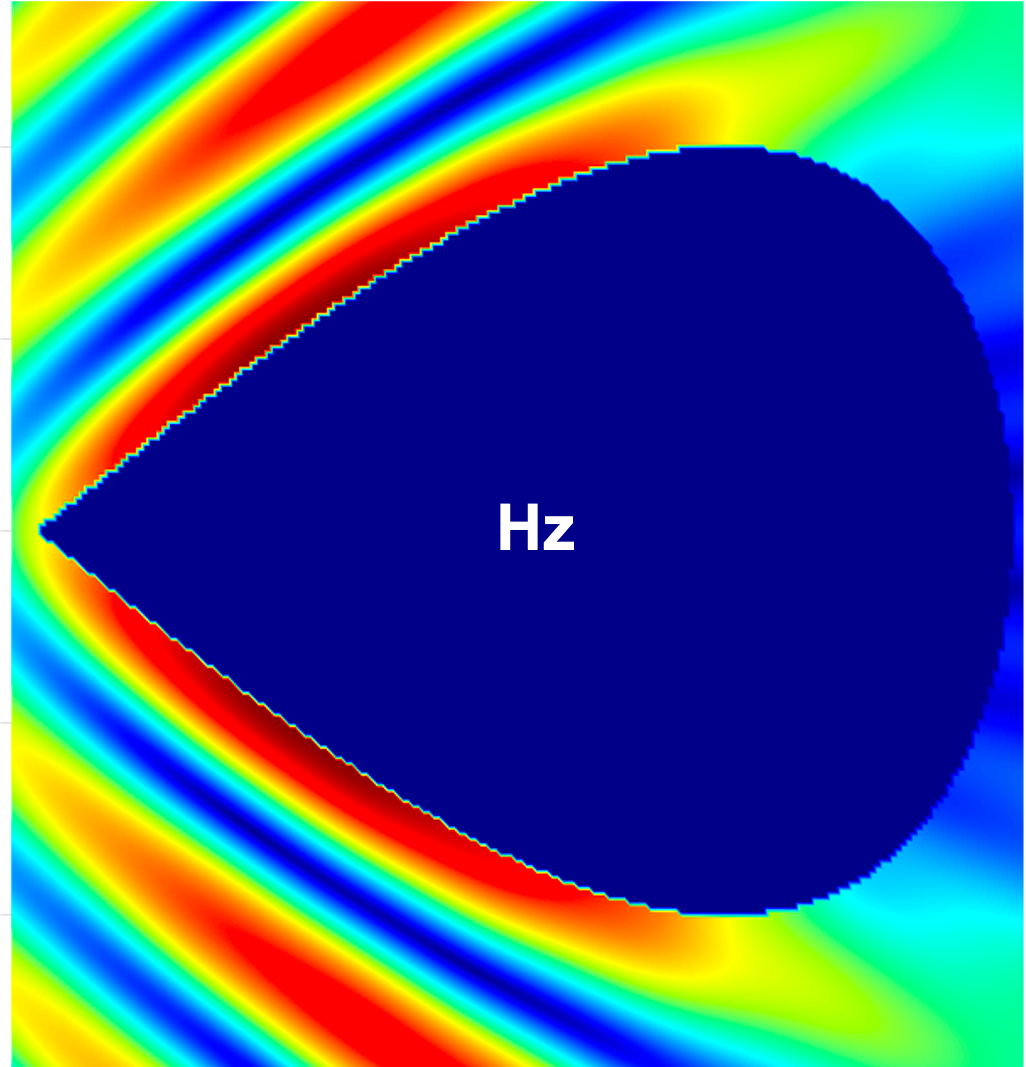}
\caption{Total field magnitude plots for $E_x,E_y$ and $H_z$ field components for cylinders with a parallelogram and teardrop shaped cross-sections}
\label{fig:scatteredfields}
\end{figure}

Figure~\ref{fig:lumericalvsbie} compares the magnitudes of the total
field components \(E_x\), \(E_y\), and \(H_z\) within a
\([-2.5, 2.5] \times [-2.5, 2.5]\) region for the PEC cylinder with a
square cross-section excited with a planewave with $k=10$, computed
using the proposed method and the commercial FDTD solver Lumerical
with an extremely fine mesh ($\frac{\lambda}{125}$ maximum grid step size). The results show excellent
agreement, validating the correctness of the proposed
method. igure~\ref{fig:scatteredfields} additionally presents the
magnitudes of the total field components \(E_x\), \(E_y\), and \(H_z\)
for the parallelogram with a \(45^\circ\) interior angle and the
teardrop with a \(90^\circ\) interior angle examples considered
above in this section.

% \begin{figure}[!ht] % example dataset
%     \centering% <-- superfluous in this example as commented by Zarko
%     \subfigure{% <-- No space here
%         \includegraphics[width=0.45\textwidth]
%         {ex_lumerical_ppt.png}}
%     \hfill
%     \subfigure{%   <-- No space here
%         \includegraphics[width=0.45\textwidth]
%         {matlab_ex_final.png}}%      <-- No space here
%     \hfill
%     \subfigure{%   <-- No space here
%         \includegraphics[width=0.45\textwidth]
%         {ey_lumerical_ppt.png}}%      <-- No space here
%     \hfill
%     \subfigure{% <-- No space here
%         \includegraphics[width=0.45\textwidth]
%         {matlab_ey_final.png}}
%     \hfill
%     \subfigure{%   <-- No space here
%         \includegraphics[width=0.45\textwidth]
%         {hz_lumerical_ppt.png}}%      <-- No space here
%     \hfill
%     \subfigure{%   <-- No space here
%         \includegraphics[width=0.45\textwidth]
%         {matlab_hz_final.png}}%      <-- No space here
% \caption{Total field magnitude plots for $E_x,E_y$ and $H_z$ field components for a cylinder with square cross-section. Results from Lumerical (in the left column) vs results from the proposed solver (in the right column).}
% \label{fig:lumericalvsbie}
% \end{figure}

\subsection{Density errors and density corner-exponent
  evaluation}\label{exponent}
\begin{figure}[!ht] % example dataset
    \centering% <-- superfluous in this example as commented by Zarko
        \includegraphics[width=0.4\textwidth]
        {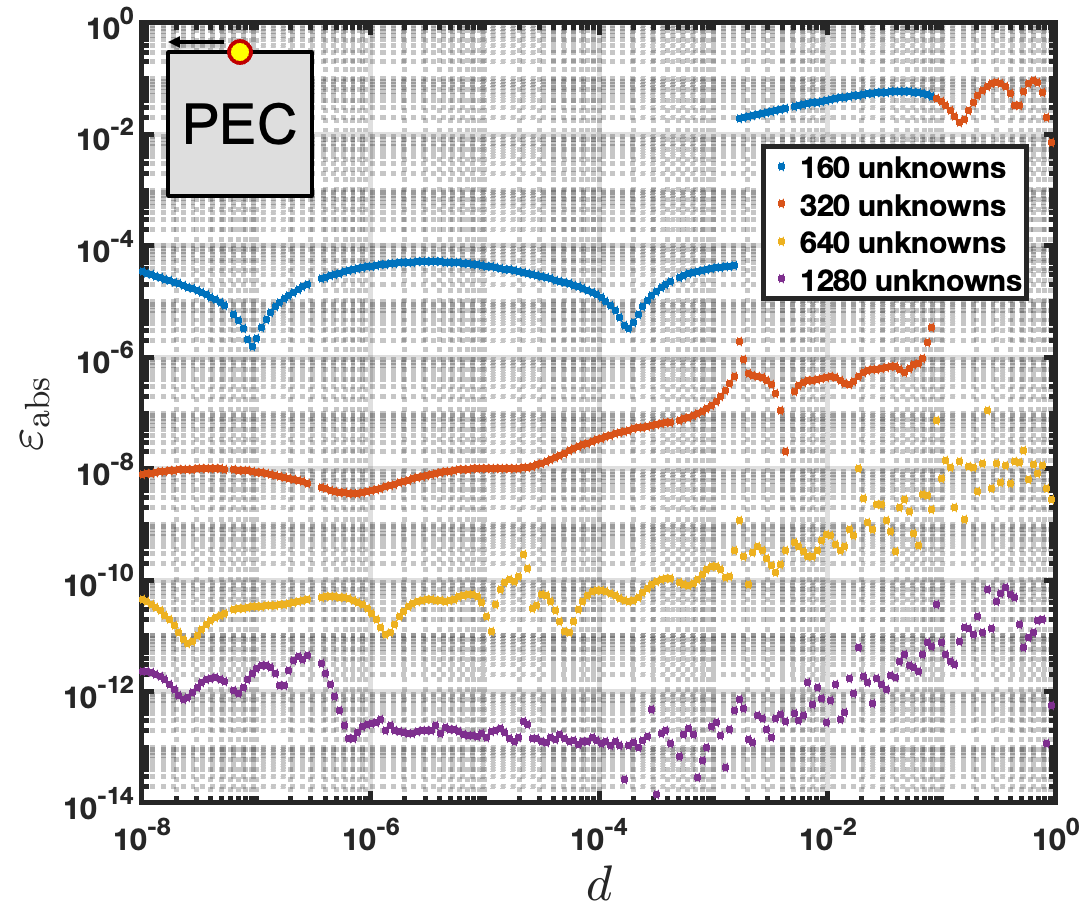}
    \hspace{0.5cm}
        \includegraphics[width=0.4\textwidth]
        {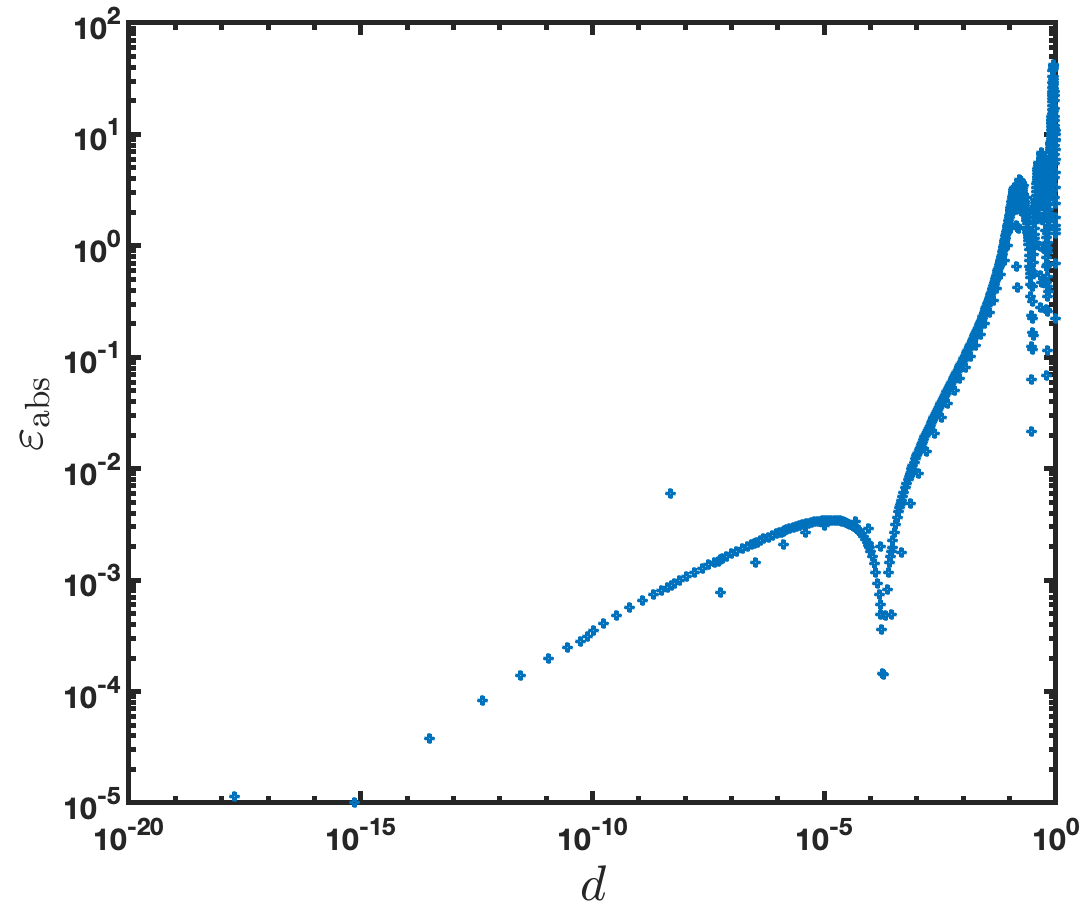}%      <-- No space here
      \caption{Left: Absolute value $\varepsilon_\mathrm{abs}$ of the
        error in the evaluation of the unknown $\psi$ for a PEC
        scatterer with a square cross-section. Right: Error in the
        numerical asymptotics of the corner exponent vs. an exact
        result near a corner point.}
\label{fig:square_errors}
\end{figure}

This section considers the accuracy obtained for the density solution
$\psi =(\psi_q)_{q=1}^M$ introduced in~\eqref{eq:ch_of_unk} at points
both near and away from a corner. The square cross-section PEC
cylinder with sidelength $a=2$ and planewave excitation as in
Figure~\ref{fig:excitation} are used for this
illustration. Figure~\ref{fig:square_errors} left displays the
absolute error in $\psi$ vs. distance from the corner for the top
horizontal boundary of the square. As expected, the error decreases as
the number of unknowns used to discretize the geometry is
increased. As discussed in Section~\ref{comparisonagainstRP}, for the
coarser discretizations the error near the corner is significantly
lower than away from the corner.  As dicussed in that section, the
errors decrease more slowly near the corner as the discretizations are
refined, which implies that the overall solution error is limited by
that near the center of the boundary for the coarser
discretizations---which explains the near 10-th order rate of
convergence observed in some cases in spite of the lower order CoV
used. The method achieves density absolute errors lower than
$10^{-11}$ for evaluation points $10^{-8}$ away from a corner.

It is interesting to additionally consider the asymptotic behavior of
the the CFIE-R density $\phi$ in~\eqref{operator_regularized} that may
be obtained numerically by division from the density $\phi$
in~\eqref{eq:ch_of_unk}. For a $90^{\circ}$ interior angle, the
density $\phi$ at a Cartesian distance $d$ from the corner, which is
denoted by $\\varphi (d)$ in what follows, tends to infinity like
$d^{\nu}$, with $\nu = -1/3$, as $d\to 0$~\cite{6}. The exponent $\nu$
can be extracted from the numerical solution as the $d\to 0$ limit of
the slope $\widetilde{\nu}(d)$ of a log-log plot of the density
$\varphi(d)$ at the distance $d$ from the corner, as a function of
$d$. The slope $\widetilde\nu$ is given by
\begin{equation}\label{asymptote_evaluation}
  \widetilde{\nu}(d) = % \frac{d \log \phi(d)}{d \log d},
  \frac{\frac{d}{dt}\varphi(e^t)}{\varphi(e^t)}\big |_{t=\log(d)}.
\end{equation}
Figure~\ref{fig:square_errors}(b), which plots
$\left|\widetilde\nu(d)-(-1/3)\right|$ vs. the distance from the
corner $d$, shows that indeed, $\widetilde\nu(d)\to -1/3$ as $d\to 0$:
the precision of the numerical solution $\psi$ suffices to produce the
singular exponent $\nu = \widetilde\nu(0)$ in the density $\phi$ with
an accuracy of at least 8 digits.

\section{Conclusions\label{conc}}
This contribution introduced and validated a novel formulation and
discretization strategy for the solution of problems of
electromagnetic scattering by objects containing geometric
singularities. The proposed method achieves high-order convergence,
comparable to that obtainable by related Nystr\"om methods for smooth
geometries, even for field evaluations at points extremely close to
corners, and even for extremely sharp, needle-like corner
angles. Building on the well-conditioned, resonance-free regularized
CFIE-R formulation~\cite{1}, the approach incorporates a
corner-regularization technique that combines a polynomial-like change
of variables (CoV) with an associated change of unknown. By
incorporating the vanishing line element from the CoV directly into
the unknown this reformulation ensures smooth behavior of the new
unknown near corners. Using additionally several precision-preserving
techniques to address the challenges posed by CoV-related clustering
of points near the corners, the method results in fast convergence and
high accuracies in the restuling scattered field, at distances both
close and far from the scatterer's corners.

Extensive numerical results confirm the unique solvability,
resonance-free character, and high-order convergence of the proposed
approach for various geometries, including corners with very small
interior angles. In most cases, the method achieves a relative error
near $10^{-12}$ with faster than 9th-order convergence for field
evaluations at arbitrary distances from a corner. Total field
magnitudes outside the scatterers are plotted for several geometries
demonstrating excellent agreement with results obtained from a
commercial FDTD solver.

In conclusion, the proposed methodology represents a significant
advancement in the numerical modeling of electromagnetic scattering
problems. It offers a robust, efficient, and highly accurate solution
framework, particularly for scatterers with geometric corner
singularities. This approach establishes a strong foundation for
future extensions to three-dimensional Maxwell scattering problems for
scatterers containing corners and edges.

\section*{Acknowledgments\label{ack}}
The authors gratefully acknowledge support by the Air Force Office of Scientific Research (FA9550-20-1-0087, FA9550-25-1-0020, FA9550-21-1-0373 and FA9550-25-1-0015) and the National Science Foundation (CCF-2047433 and DMS-2109831).
%% References with BibTeX database:
%\bibliographystyle{cas-model2-names}
\bibliographystyle{acm}
\bibliography{2dbie-refs}

%% \bibliographystyle{elsarticle-num}
%% \bibliography{<your-bib-database>}

%% Authors are advised to use a BibTeX database file for their reference list.
%% The provided style file elsarticle-num.bst formats references in the required Procedia style

\end{document}